\documentclass{dis}
\usepackage{amssymb,amsmath,amsthm,graphicx}
\usepackage{diagrams}
\usepackage{bbm}
\usepackage[titles]{tocloft}

\theoremstyle{definition}
\newtheorem{defi}{Definition}[section]
\newtheorem*{claim}{Claim}
\theoremstyle{plain}
\newtheorem{thr}{Theorem}
\newtheorem{lem}[defi]{Lemma}
\newtheorem{cor}[defi]{Corollary}
\newtheorem{pro}[defi]{Proposition}

\newtheorem{problem}[defi]{Problem}
\theoremstyle{remark}
\newtheorem{rem}[defi]{Remark}

\newtheorem{exa}[defi]{Example}

\newcommand{\A}{\mathcal{A}}
\newcommand{\F}{\mathcal{F}}

\renewcommand{\H}{\mathcal{H}}
\newcommand{\M}{\mathcal{M}}
\newcommand{\Simple}{\mathcal{S}}
\newcommand{\W}{\mathcal{W}}
\renewcommand{\P}{\mathcal{P}}

\newcommand{\Algebra}{\mathfrak A}
\newcommand{\U}{\mathfrak{U}}
\renewcommand{\P}{\mathfrak{P}}

\newcommand{\T}{\mathbb{T}}
\newcommand{\Z}{\mathbb{Z}}
\newcommand{\N}{\mathbb{N}}
\newcommand{\R}{\mathbb{R}}
\newcommand{\C}{\mathbb{C}}
\newcommand{\Q}{\mathbb{Q}}

\newcommand{\Eins}{\mathbbm{1}}
\newcommand{\eps}{\varepsilon}
\renewcommand{\phi}{\varphi}

\newcommand{\Hartset}{\mathfrak{H}}
\newcommand{\Riemann}{\mathcal{I}}
\newcommand{\Riemanntop}{\mathcal{R}}

\newcommand{\limn}{\lim_{n \to \infty}}

\newcommand{\e}{\varepsilon}

\newcommand{\tm}{\subseteq}

\begin{document}
\mathclass{Primary 43A60; Secondary 26A42.}
\thanks{The authors would like to thank the Austrian Science Fund (FWF) for
financial support through grants S8312, S9612 and Y328.}
\abbrevauthors{G. Maresch and R. Winkler}
\abbrevtitle{Hartman functions}

\title{Compactifications, Hartman functions and (weak) almost periodicity}

\author{Gabriel Maresch}
\address{Institute of Discrete Mathematics and Geometry\\ Vienna University of Technology\\
Wiedner Hauptstra\ss e 8-10/104\\ 1040 Vienna, Austria\\
E-mail: gabriel.maresch@tuwien.ac.at}

\author{Reinhard Winkler}
\address{Institute of Discrete Mathematics and Geometry\\ Vienna University of Technology\\
Wiedner Hauptstra\ss e 8-10/104\\ 1040 Vienna, Austria\\
E-mail: reinhard.winkler@tuwien.ac.at}

\maketitledis

\begin{abstract}
In this paper we investigate Hartman functions on a topological group $G$. Recall that $(\iota, C)$ is a group compactification of $G$ if $C$ is a compact group, $\iota: G\to C$ is a continuous group homomorphism and $\iota(G)\tm C$ is dense. A bounded function $f: G \mapsto \C$ is a Hartman function if there exists a group compactification $(\iota, C)$ and $F: C \to \C$ such that $f= F\circ\iota$ and $F$ is Riemann integrable, i.e.\  the set of discontinuities of $F$ is a null set w.r.t.\ the Haar measure. In particular we determine how large a compactification for a given group $G$ and a Hartman function $f: G\to \C$ must be, to admit a Riemann integrable representation of $f$. The connection to (weakly) almost periodic functions is investigated.

\noindent In order to give a systematic presentation which is self-contained to a reasonable extent, we include several separate sections on the underlying concepts such as finitely additive measures on Boolean set algebras, means on algebras of functions, integration on compact spaces, compactifications of groups and semigroups, the Riemann integral on abstract spaces, invariance of measures and means, continuous extensions of transformations and operations to compactifications, etc.
\end{abstract}

\makeabstract

\tableofcontents

%
%

\chapter{Introduction}\label{CHintro}
\section{Motivation}
By a topological dynamical system $(X,T)$ we mean a continuous transformation $T:X\to X$ acting on a compact space $X$ (which in many cases is supposed to be metrizable). Symbolic dynamics is concerned with the special case $X=A^{\N}$ or $X=A^{\Z}$ with a finite set $A$, called the alphabet. Here the transformation is the shift $T=\sigma: (a_n)\in X \mapsto (a_{n+1})\in X$. The importance of this special case is due to the fact that, for a suitable finite partition (Markov partition) $X=X_1\cup\ldots\cup X_n$ of a metrizable space $X$ and the alphabet $A=\{1,\ldots,n\}$, most information of the original system $(X,T)$ is contained in the associated symbolic system which is defined below.

Consider the coding $F: X\to A$, $F(x)=i$ if $x\in X_i$. Let $\phi: X\to A^{\N}$, $x\mapsto (T^n x)_{n\in\N}$ or, if $T$ is bijective, $\phi: X\to A^{\Z}$, $x\mapsto (T^n x)_{n\in\Z}$. The case of bijective $T$ applies for the major part of the exposition. The associated dynamical system $(Y,\sigma)$ with $Y=\overline{\phi(X)}$ is a subshift, i.e.\ $Y$ is a closed and $\sigma$-invariant subset of $A^{\Z}$. The connection between $(X,T)$ and $(Y,\sigma)$ is expressed by the commuting diagram:
\begin{diagram}
X            &\rTo^{T}&X\\
\dTo>{\phi}&         &\dTo>{\phi} \\
Y            &\rTo^{\sigma} &Y.
\end{diagram}
If $\phi$ is continuous this means that $(Y,\sigma)$ is a factor of $(X,T)$. However, this can be guaranteed only if the $X_i$ are clopen subsets of $X$ which, for instance for connected $X$, is impossible. The classical way of avoiding this disadvantage is to choose the partition in such a way that $\phi$ is injective and $\phi^{-1}$ has a continuous extension $\psi$ such that $(X,T)$ is a factor of $(Y,T)$:
\begin{diagram}
Y            &\rTo^{\sigma}&Y\\
\dTo>{\psi}&         &\dTo>{\psi} \\
X            &\rTo^{T} &X.
\end{diagram}
In order to apply results from ergodic theory (such as Birkhoff's Theorem) one looks for invariant measures. Assume that $\mu$ is such a $\sigma$-invariant measure on $Y$, i.e.\ $\mu(\sigma^{-1}[B])=\mu(B)$ for all Borel sets $B\subseteq Y$. Then $\mu_T(M):=\mu(T^{-1}[M])$ defines a $T$-invariant measure $\mu_T$ on $X$.

The situation is particularly nice if $T$ is uniquely ergodic, i.e.\ if there is a unique $T$-invariant Borel measure. In this case the limit relation
\begin{align}
\lim_{N\to \infty} \frac 1N \sum_{n=0}^{n-1} f(T^nx)=\int_X f d\mu_T
\end{align}
does hold not only up to a set of zero $\mu_T$-measure, but even uniformly for all $x\in X$ whenever $f: X\to \R$ is continuous and bounded. By obvious approximation this statement extends to all bounded
$f: X\to \R$ with
\begin{align}\label{riemann}
\forall \eps >0 \quad \exists f_1,f_2: X\to \R \mbox{ continuous}, \quad f_1\le f \le f_2, \quad \int_X(f_2-f_1)d\mu_T < \eps.
\end{align}
In the case $X=[0,1]$, equipped with the Lebesgue measure, \eqref{riemann} is equivalent with the requirement $\mu_T(\mbox{disc})=0$, i.e.\ that the set $\mbox{disc}(f)$ of discontinuity points of $f$ is a null set. In other words, $f$ is Riemann integrable. If $f$ takes only finitely many values $r_1,\ldots,r_s$ this condition is equivalent with $\mu_T(\partial X_i)=0$ for the topological boundary of $X_i:=f^{-1}[\{r_i\}]$, $i=1,\ldots, s$. Indeed, this condition is usually assumed for partitions in the context of symbolic dynamics. In this paper we allow $F: X\to \C$ to have infinitely many values, but, motivated by the above considerations, assume that $F$ is Riemann integrable.

A very important class of uniquely ergodic systems are group rotations, i.e.\ $T: C\to C$, $x\mapsto x+g$ where $g\in C$ is a topological generator of the compact (abelian) group $C$, meaning that the cyclic group generated by $g$ is dense in $C$. The unique invariant measure for the transformation $T$ is given by the Haar measure $\mu_C$ on $C$. The induced coding sequences $(a_n)_{n\in \Z}$ are given by $a_n=F(x+ng)$ and may be used to form a factor of $(\iota,C)$ . Indeed, if we consider the mapping $\iota: \Z\to C$, $n\mapsto ng$, we have $a=F\circ\iota$. $(\iota, C)$ is a group compactification of $\Z$ since $\iota$ is a (trivial) continuous group homomorphism with image $\iota(\Z)$ dense in $C$. Allowing $\Z$ to be replaced by an arbitrary topological group, we finally arrive at the definition of Hartman functions, the main objects of our paper:

A function $f: G\to \C$ on a topological group $G$ is called a Hartman function if there is a group compactification $(\iota, C)$ of $G$ and a function $F: C\to \C$ which is Riemann integrable w.r.t.\ the Haar measure and satisfies $f=F\circ\iota$. $F$ is called a representation of $(\iota, C)$.

In particular, almost periodic functions (defined by continuous $F$) are Hartman functions. The name \emph{Hartman} function refers to the Polish mathematician Stanis{\l}aw Hartman who was, up to our knowledge, the first to consider these objects in the 1960s in his work in harmonic analysis \cite{Hart61,Hart63,Hart64}. He focused on the Bohr compactification $(\iota_b, bG)$ of the group $G$. It is not difficult to see that our definition is equivalent with the requirement $(\iota, C)=(\iota_b, bG)$. The question whether for a given Hartman function $f$, there are small compactifications with a representation $f=F\circ\iota$ is one of our major topics.

Additionally we investigate the connection of Hartman functions and weak almost periodicity. Recall that a function is weakly almost periodic if it has a continuous representation in a semitopological semigroup compactification, or, equivalently, in the weak almost periodic compactification $(\iota_w,wG)$. While every almost periodic function is Hartman, this is not true in the weak case. A more systematic overview of the content of this paper is given at the end of this Section.

\section{Recent results on Hartman sets, sequences and functions}
 For an extended survey on recent research on Hartman sets, Hartman sequences and Hartman functions we refer to \cite{Wink07}. Here we only give a very brief summary.

The series of papers we report on was initiated by investigations of M.\ Pa\v{s}teka and R.F.\ Tichy \cite{Past96,Past97,PaTi94} on the distribution of sequences induced by the algebraic structure in commutative rings $R$. The authors used the completion $\overline{R}$ w.r.t.\ a natural metric structure such that $R$ is compact and thus carries a Haar measure $\mu$. The restriction of $\mu$ to the $\mu$-continuity sets $M$, i.e.\ to those sets with $\mu(\partial M)=0$ has been pulled back in order to obtain a natural concept of uniform distribution in the original structure $R$.

One easily observes that the measure theoretic part of the construction depends only on the additive group structure of $R$. Thus the natural framework for a systematic investigation is that of group compactifications $(\iota, C)$ of a topological group $G$ and of the finitely additive measure $\mu_{(\iota,C)}$ on $G$ defined for $\iota$-preimages of $\mu$-continuity sets as follows
 \begin{align}
  \mu_{(\iota,C)}\left(\iota^{-1}[M]\right):=\mu(M),\quad M\subseteq C.
 \end{align}
 This has been studied in \cite{FPTW99}. Results for the special case $G=\Z$ are presented in \cite{Sand02,ScSW00}: Hartman sets $\iota^{-1}[M]\subseteq \Z$ are identified with the function $\Eins_{\iota^{-1}[M]}: \Z \to \{0,1\}$ and
called Hartman sequences. The relation to Beatty resp.\ Sturmian sequences and continued fractions expansion is described. It is shown that the system of Hartman sequences is generated by the system of Beatty sequences by means of Boolean combinations and approximation in measure.

The connection to ergodic theory already mentioned in \cite{ScSW00} is stressed further in \cite{Wink02}: Hartman sequences can be considered as symbolic coding sequences of group rotations (as described in the previous section). The problem to identify the underlying dynamical system turns out to be equivalent to the identification of the group compactification $(\iota, C)$ of $\Z$ inducing the Hartman set $\iota^{-1}[M]\tm \Z$. As an alternative to classical methods such as spectral analysis of the dynamical system, a purely topological method has been presented. Each Hartman set $\iota^{-1}[M]\tm \Z$ defines in a natural way a filter on $\Z$. Under rather mild assumptions this filter is the $\iota$-preimage of the neighborhood filter $\mathfrak{U}(0_C)$ of the identity in $C$ and contains all necessary information about $(\iota, C)$.

These methods have been applied to questions from number theory in \cite{Beig07} and generalized to the setting of topological groups in \cite{BeSW06}.

The aspect of symbolic dynamics has been studied further in \cite{StWi05} by investigation of subword complexity of Hartman sequences. Recall that the subword complexity $p_a: \N \to \N$ induced by the sequence $a\in \{0,1\}^{\Z}$ is a function associating to each $n\in\N$ the number of different 0-1 blocks of length $n$ occurring in $a$. Clearly $1\le p_a(n)\le 2^n$. The main facts in this context are:
\begin{enumerate}
\item $\lim_{n\to\infty} \tfrac 1n \log p_a(n) = 0$, corresponding to the fact that group rotations have entropy $0$.
\item Whenever $\lim_{n\to\infty} \tfrac 1n \log p_n = 0$ for a sequence $p_n$ with $1\le p_n \le 2^n$, then there is a Hartman sequence $a$ with $p_n(a)\ge p_n$ for every $n\in\N$.
\item\label{projectionbody} The Hartman sequence  $a=\Eins_{\iota^{-1}[M]}$, where $M\tm \T^s$ is an $s$-dimensional cube, satisfies $p_a(n)\sim c_M\cdot n^s$ with an explicit constant $c_M>0$ (we omit the number theoretic assumptions).
\end{enumerate}

An amazing geometric interpretation of the constant $c_M$ was recently given in \cite{Stein08}, where statement \ref{projectionbody} has been generalized to convex polygons $M$ and $c_M$ corresponds to the volume of the projection body of $M$.

The investigation of Hartman functions has been started in \cite{Mare05} where, for instance, results from \cite{Wink02} on Hartman sequences have been generalized. In the present paper we continue these investigations and include a systematic and considerably self-contained treatment of the topological and measure theoretic background.

\section{Content of the paper}
 Chapter \ref{CHprel} presents measure theoretic and topological preliminaries. Section \ref{Salgebras} fixes notation concerning (Boolean) set algebras and related algebras of functions. In Section \ref{Sfam} we investigate finitely additive measures on set algebras and the integration of functions from corresponding function algebras. Consequently we present the connection between measures and means. One of the most fundamental phenomena in analysis is that compactness is used to obtain $\sigma$-additivity of measures and thus makes Lebesgue's integration theory work. Riesz' Representation Theorem plays a crucial r\^ole in this context, which we recall in Section \ref{Scompintegral}. For the case that compactness is absent one can try to force compactness by considering compactifications. In Section \ref{Scompcont} we construct compactifications in such a way that a given set of bounded functions admits continuous extensions. We touch the classical representation theorems of Gelfand and Stone. Among all compactifications of a given (completely regular) topological space $X$ there is a, in a natural sense, maximal compactification, the Stone-\v{C}ech compactification $(\iota_{\beta}, \beta X)$. In Section \ref{Sstonecech} we collect important properties. Having presented the basics concerning compactifications, measures, means and the Riemann integral, we put these concepts together in Section \ref{Scompriemann}. Section \ref{Smeans}, the last one in Chapter \ref{CHprel}, presents the interpretation of the Stone-\v{C}ech compactification $(\iota_{\beta}, \beta X)$ of a discrete space $X$ as the set of all multiplicative means on $X$.  This motivates us to investigate means with more restrictive properties, such as invariance.

Chapter \ref{CHinvariance} is concerned with invariance of measures and means under transformations and operations. In particular we investigate in Section \ref{Ssingle} questions of existence and uniqueness. For a transformation $T: X\to X$ invariance is closely related to the behavior of Ces\'aro means along $T$-orbits,  a concept which leads to the notion of Banach-density. In Section \ref{Sapplications-means} we treat several examples and applications: finite $X$, $X=\Z$ and $T:x\mapsto x+1$, compact $X$ and continuous $T: X\to X$, shift spaces and symbolic dynamics, the free group generated by two elements. In Section \ref{Strafocomp} we consider compactifications under the additional aspect of extending transformations and (semi)group actions in a continuous way. For binary or, more generally, $n$-ary, operations continuous extensions do not always exist. The arising problems are treated in Section \ref{Sjoint}. In particular we show that $n$-ary operations on $X$, $n\ge 2$, can be continuously extended to $(\iota_{\beta},\beta X)$ only in very special cases. Nevertheless it is useful to formulate a general framework in order to unify the most interesting classical situations: topological and semitopological group and semigroup compactifications. This is done in Section \ref{Sopcomp}. In Section \ref{Sinvgroups} these constructions are discussed in the context of invariant means and measures. We mention the notion of weak almost periodicity and touch very briefly amenable groups and semigroups.

Chapter \ref{CHhartman} develops the basic theory of Hartman functions. Section \ref{Shartman} presents several equivalent conditions describing the connection with almost periodicity and the Bohr compactification, i.e.\ the maximal group compactification. Replacing group compactifications by semitopological semigroup compactifications one obtains the weak almost periodic compactification, weak almost periodic functions and weak Hartman functions. This is presented in Section \ref{Sweakhartman}. The category of all group compactifications of a topological group $G$ is particularly well understood if $G$ is abelian and carries a locally compact group topology. The key ingredient is Pontryagin's Duality Theorem. We recall this situation in Section \ref{Slcacomp}. One of the most interesting questions concerning a Hartman function $f: G\to \C$ is how small a group compactification $(\iota,C)$ can be taken if one asks for a Riemann integrable representation of $f$. This question is treated in Section \ref{Srealize}.
We give an answer for LCA groups in terms of the minimal cardinality of a dense subgroup in the Pontryagin dual $\hat G$ of $G$.

Chapter \ref{CHwap} is devoted to the comparison of Hartman functions and weakly almost periodic functions. It turns out that a generalization of what is called a jump-discontinuity in basic analysis plays an important r\^ole. Generalized jump discontinuities are established in Section \ref{Sgjd} and used in Section \ref{Shartmanwap} to give necessary conditions of weak almost periodicity of Hartman functions. This leads to the investigation of Hartman functions without such generalized jumps in Section \ref{Shartmannogjd}. Hartman functions with small support are treated in \ref{Shartmansmall}. Finally, Section \ref{ShartmanZ} discusses particular examples of Hartman functions on the integers which are neither almost periodic nor converge to $0$. The results use the Fourier-Stieltjes transform of measures.

Finally a short summary is given, including a diagram which illustrates the relation between several spaces of functions which are interesting in our context.


\chapter{Measure theoretic and topological preliminaries}\label{CHprel}
\section{Set algebras $\Algebra$ and $\Algebra$-functions}\label{Salgebras}

We start with fixing notation which is suitable to imitate the
construction of the Riemann integral in the slightly more general
context which will be ours.

\begin{defi}\label{Dsetalgebra}
A (boolean) \textbf{set algebra} $\Algebra$ (on a set $X$) is a
system of subsets of $X$ with $\emptyset,X \in \Algebra$ for
which $A,B \in \Algebra$ implies
$A \cup B, A \cap B, X \setminus A \in \Algebra$.
\end{defi}

\begin{exa}\label{Eborel01}
Let $X=[0,1] \tm \R$ be the unit interval and $\Algebra = \Algebra([0,1])$ the system of all finite unions of subintervals
$I \tm [0,1]$ (open, closed and one-sided closed, also
including singletons and the empty set). This is the most classical
situation. But it is worth to note that we might replace
$[0,1]$ by any totally ordered $X$,
for instance by any $D \tm [0,1]$ dense in $[0,1]$ (as $D = \Q\cap [0,1]$).
\end{exa}

We are interested in the integration of complex valued functions
on $X$:

\begin{defi}\label{Dfunctionalgebra}
Let $\A$ be a set of functions $f: X \to \C$. We call the subset
$\A_{\R}$ of all $f \in \A$ with $f(X) \tm \R$ the \textbf{real part}
of $\A$. $\A$ is called real if $\A_{\R}=\A$. If $\A$ is a vector
space or an algebra over $\R$ (or $\C$) we call $\A$ a real (or
complex) \textbf{space} resp.\ a real (or complex) \textbf{algebra of
functions}. For any $A \tm X$ let $\Eins_A(x)=1$ for $x \in A$ and
$\Eins_A(x)=0$ for $x \in X \setminus A$. For an algebra $\A$ we
always assume $\Eins_X \in \A$. A complex space or algebra $\A$ of
functions is called a \textbf{$*$-space} resp.\ a \textbf{$*$-algebra} if
$f\in \A$ implies $\overline f \in \A$ for the complex conjugate
$\overline f$ of the function $f$. We write $B(X)$ for the set of
all bounded $f: X \to \C$, $B_{\R}(X):=B(X)_{\R}$ for its real part.
A *-algebra $\A$ on $X$ which is complete with respect to the
topology of uniform convergence on $X$ is called a \textbf{$C^*$-algebra}.
\end{defi}

Note that whenever $\A$ is a real space we can form the
complexification $\A_{\C} = \{f_1+if_2:\ f_1,f_2 \in \A\}$ which is
a complex vector space, and a $*$-algebra whenever $\A$ is a real
algebra. For any complex linear space or algebra $\A$, to be a
$*$-space resp.\ a $*$-algebra is equivalent with the following
property: Whenever $f = f_1+if_2$ is the decomposition of $f$ into
real part $f_1$ and imaginary part $f_2$, then $f \in \A$ if and
only if $f_1,f_2 \in \A_{\R}$. Thus for the investigation of
$*$-algebras $\A$ it suffices to investigate the real part $\A_{\R}$
whenever convenient. Furthermore any $*$-algebra of functions is
closed under taking absolute values: $|f|=\sqrt{f\overline{f}}$; a fact which can be seen
by approximating the square-root by polynomials.

\begin{defi}\label{Dsimplefunction}
Let $\Algebra$ be a set algebra on $X$. A function $f: X \to \C$ is
called \textbf{$\Algebra$-simple} if it has a representation
\[
f = \sum_{i=1}^n c_i \Eins_{A_i}
\]
with $A_i \in \Algebra$ and $c_i \in \C$.
The set of all $\A$-simple $f$ is denoted by $\Simple_{\Algebra}$.
We denote the uniform closure  $\overline{\Simple_{\Algebra}}$
of $\Simple_{\Algebra}$ by $B(\Algebra)$.
Members of $B(\Algebra)$ are also called \textbf{$\Algebra$-functions}.
\end{defi}

More explicitly, for a set algebra $\Algebra$ on $X$ the
function $f: X \to \C$ lies in $B(\Algebra)$ if and only if
for all $\e>0$ there is a $f' \in \Simple_{\Algebra}$ with
$|f(x)-f'(x)|<\e$ for all $x \in X$.

\begin{pro}\label{Palgebrainclusion1}
All of the sets $\Simple_{\Algebra} \tm B(\Algebra) \tm B(X)$ are
$*$-algebras. In general the inclusions cannot be replaced by
equality.
\end{pro}
\begin{proof}
It is clear that $\Simple_{\Algebra}$, $B(\Algebra)$ and $B(X)$ are
$*$-algebras satisfying the stated inclusions. Thus it suffices to
show that $\Simple_{\Algebra} \neq B(\Algebra) \neq B(X)$ if one takes
$X=[0,1]$ and $\Algebra = \Algebra([0,1])$
the set algebra of all finite unions of
subintervals of $[0,1]$. Then $f \in C(X) \tm B(\Algebra)$
but $f \notin \Simple_{\Algebra}$ if we take $f(x)=x$, hence
$\Simple_{\Algebra} \neq B(\Algebra)$. On the other side all
$f \in B(\Algebra)$ are Riemann integrable in the classical sense
which is not the case for arbitrary $f \in B(X)$.
\end{proof}

For every set algebra $\Algebra$, $B(\Algebra)$ is a
$C^*$-algebra. But not every $C^*$-algebra $\A$ can be written as
$\A=B(\Algebra)$ for an appropriate $\Algebra$. The situation is
explained by the following facts.

\begin{pro}\label{Psetsversusfunctions}
For a set $\A$ of complex valued functions $f: X \to \C$ let
$\Algebra_{\A}:= \{A \tm X:\ \Eins_A \in \A\}$. Then:
\begin{enumerate}
\item $\Algebra_{\A}$ is a set algebra whenever $\A$ is an algebra.
\item Every set algebra $\Algebra$ on $X$ satisfies
   $\Algebra = \Algebra_{B(\Algebra)}$.
\item For every uniformly closed algebra $\A$ one has
  $B(\Algebra_{\A}) \tm \A$ while the converse inclusion does not hold
  in general.
\end{enumerate}
\end{pro}
\begin{proof}
\hfill\par\vspace{-2ex}
\begin{enumerate}
\item Follows from $\Eins_X \in \A$,
$\Eins_{A_1 \cap A_2} = \Eins_{A_1} \cdot \Eins_{A_2}$,
$\Eins_{X \setminus A} = \Eins_X - \Eins_A$ and the identity
$A_1 \cup A_2 = X \setminus ((X\! \setminus\! A_1) \cap (X\! \setminus\! A_2))$.
\item
The inclusion $\Algebra \tm \Algebra_{B(\Algebra)}$ is obvious.
For the converse assume
$A \in \Algebra_{B(\Algebra)}$, i.e.\ $\Eins_A \in B(\Algebra)$.
Then there are $f_n \in \Simple_{\Algebra}$ uniformly converging to $\Eins_A$.
There are representations $f_n = \sum_{i=1}^{k_n}\alpha_{i,n}\Eins_{A_{n,i}}$
such that for each $n$ the $A_{n,i} \in \Algebra$, $i=1,\ldots,k_n$,
are pairwise disjoint. For sufficiently large fixed $n$,
each $x \in X$ satisfies
either $|f_n(x)-1| < \tfrac 12$ (if $x \in A$) or $|f_n(x)| < \tfrac 12$
(if $x \notin A$). This shows that $A_{n,i} \tm A$ or
$A_{n,i} \tm X \setminus A$ for any such fixed $n$ and all
$i=1,\ldots,k_n$, hence $A = \bigcup_{i: A_{n,i} \tm A} A_{n,i} \in \Algebra$.
\item
The stated inclusion is obvious. The example $\A = C([0,1])$,
$\Algebra_{\A} = \{\emptyset, X\}$,
$B(\Algebra_{\A}) = \{c\Eins_X:\ c \in \C\}$,
shows that the inclusion might be strict. \qedhere
\end{enumerate}
\end{proof}

\section{Finitely additive measures and means}\label{Sfam}

\begin{defi}\label{Dfam}
Let $\Algebra$ be a set algebra on $X$. A function $p: \Algebra \to [0,\infty]$ with $p(\emptyset)=0$ is called a \textbf{finitely additive measure}, briefly \textbf{fam} (on $X$ or, more precisely, on $\Algebra$) if it is finitely additive, i.e.\ if $p(A_1 \cup A_2) = p(A_1)+p(A_2)$ whenever $A_1 \cap A_2 = \emptyset$. $p$ is called a \textbf{finitely additive probability measure}, briefly \textbf{fapm},
if furthermore $p(X)=1$.
\end{defi}

\begin{exa}\label{Eborel02}
Continuing Example \ref{Eborel01}, for $X=[0,1]$ and $\Algebra=\Algebra([0,1])$, the
system of all finite unions of intervals, one takes $p(I) = b-a$
for $I=[a,b]$ with $0 \le a \le b \le 1$. This definition uniquely
extends to a fapm on the set algebra $\Algebra([0,1])$
of all finite unions of intervals.
We will refer to this $p$ as the natural measure.
The construction does not depend on the completeness
(compactness) of $[0,1]$ and hence can be done as well for dense
subsets $D \subset X$. For instance one could consider (finite
unions of) intervals of rationals.
\end{exa}

\begin{defi}\label{Dmean}
Let $\A$ be a linear space of functions on a set $X$.
Then a \textbf{mean} $m$ on $\A$ is a linear functional $m: \A \to \C$
which is \textbf{positive}, i.e.\ $f \ge 0$ implies $m(f) \ge 0$, and satisfies
$m(\Eins_X)=1$.
\end{defi}

Note that whenever $\A$ is real and $m$ is a mean on $\A$ then
$\overline m(f_1+if_2):=m(f_1)+im(f_2)$ for $f_1,f_2 \in \A$ is the unique
extension $\overline m$ to the complexification $\A_{\C}$ of $\A$.
Very often we simply write $m$ for $\overline m$.

For real functions $f$ every mean $m$, by positivity,
satisfies $\inf f \le m(f) \le \sup f$. As a consequence we have:

\begin{pro}\label{Pmeancontinuous}
Every mean $m$ on $\A$ is continuous with respect to the norm
$||f||_{\infty}:= \sup_{x \in X}|f(x)|$ and thus has a unique
extension to the uniform closure $\overline{\A}$ of $\A$.
\end{pro}

Every mean induces a further notion of closure:

\begin{defi}\label{Dmclosure}
Let $m$ be a mean on a linear space $\A$ of functions on $X$. Then
the \textbf{real $m$-closure} $\overline{\A_{\R}}^{(m)}$ of $\A$ is the
set of all $f: X \to \R$ such that for all $\e>0$ there are $f_1,f_2
\in \A_{\R}$ with $f_1 \le f \le f_2$ and $m(f_2-f_1) < \e$. For $f
\in \overline{\A_{\R}}^{(m)}$, $\overline{m}(f)$ is defined to be
the \emph{unique} value $\alpha \in \R$ with $m(f_1) \le \alpha \le
m(f_2)$ for all $f_1,f_2 \in \A$ with $f_1 \le f \le f_2$. The
(complex) \textbf{$m$-closure} $\overline{\A}^{(m)}$ is the set of all
$f=f_1+if_2$ with $f_1,f_2 \in \overline{\A}^{(m)}_{\R}$.
Furthermore we define $\overline{m}(f):=
\overline{m}(f_1)+i\overline{m}(f_2)$ for such $f=f_1+if_2$.
$\overline{m}$ is called the \textbf{completion} of $m$, sometimes also
simply denoted by $m$. In the case $\overline{\A}^{(m)} = \A$ we
call $m$ \textbf{complete} and $\A$ \textbf{$m$-closed}.
\end{defi}

\begin{rem}
Distinguish the $m$-closure from the completion
with respect to the pseudo-metric $d_m(f,g):= m(|f-g|)$.
By definition ($m$ is continuous w.r.t.\ $d_m$) the
$m$-closure is always contained in the $d_m$-completion:
$\overline{\A}^{(m)}\subseteq \overline{\A}^{(d_m)}$.
The closure w.r.t.\ $m$ corresponds to the integral in the sense
of Riemann, the completion w.r.t.\ $d_m$ to that of Lebesgue (modulo null-sets).
\end{rem}

Every fapm $p$ defined on a set algebra $\Algebra$ on a set $X$
induces a linear functional $m_p$ in the natural way.
Standard arguments (using that $\Algebra$ is closed under intersections
and that $p$ is finitely additive) show that for an $\Algebra$-simple
$f= \sum_{i=1}^n c_i\Eins_{A_i} \in \Simple_{\Algebra}$ the value
\[
m_p(f) = m_p\left(\sum_{i=1}^n c_i \Eins_{A_i} \right):=
  \sum_{i=1}^n c_i p(A_i)
\]
does not depend on this particular representation of $f$ as
a linear combination. Obviously this $m_p$
is a mean on $\Simple_{\Algebra}$ and thus,
by Proposition \ref{Pmeancontinuous}, has a unique extension to the
algebra $B(\Algebra) = \overline{\Simple_{\Algebra}}$ as well
as to $\overline{\Simple_{\Algebra}}^{(m_p)}$.

We want to extend the domain of $m_p$ from $\Simple_{\Algebra}$ to the space
$\Riemann_p$ defined as follows.

\begin{defi}\label{Driemannp}
For a given fapm $p$ on $\Algebra$ let
$\Riemann_p:= \overline{\Simple}_\Algebra^{(m_p)}$.
The members $f \in \Riemann_p$ are called \textbf{integrable} (w.r.t.\ $p$).
The extension of $m_p$ to $\Riemann_p$, usually also denoted
by $m_p$, is called the mean induced by $p$.
\end{defi}

We leave the proof of the following easy properties to the reader:

\begin{pro}\label{Palgebrainclusion2}
Let $\Algebra$ be a set algebra on the set $X$, $p$ a fapm defined
on $\Algebra$. Then $B(\Algebra) \tm \Riemann_p \tm B(X)$, $\Riemann_p$ is $m_p$-closed  and $m_p$ is a mean on $\Riemann_p$.
\end{pro}

\begin{rem}\label{Rcstaralgebra}
$\Riemann_p$ is uniformly closed. In particular $\Riemann_p$ is a $C^*$-algebra.
Indeed, let $f_n \to f$ uniformly where $f_n \in \Riemann_p$. For given $\eps>0$ there
exists $f_n$ such that $\|f-f_n\|_{\infty}\le \tfrac{\eps}{4}$ and  $f_{n,1}, f_{n,2} \in B(\Algebra)$ such that $f_{n,1}\le f_n\le f_{n,2}$ and $m_p(f_{n,2}-f_{n,1})\le\tfrac{\eps}{2}$. Observe that
\[f_{n,1} - \tfrac{\eps}{4} \le f_n-\tfrac{\eps}{4} \le f \le f_n+\tfrac{\eps}{4} \le f_{n,2} + \tfrac{\eps}{4}\]
and thus $m_p\left( (f_{n,2} + \tfrac{\eps}{4})-(f_{n,1} - \tfrac{\eps}{4})\right)\le \eps$ shows $f\in\Riemann_p$.
\end{rem}
The inclusions stated in Proposition \ref{Palgebrainclusion2} are in general strict as the following example shows.

\begin{exa}\label{Eclosureinclusions}
Let again $\Algebra=\Algebra([0,1])$
be the set algebra of all finite
unions of subintervals of $X=[0,1]$, $p$ the natural measure on $\Algebra$.
Then $\Riemann_p$ is the set of all $f:[0,1] \to \C$ which are
integrable in the classical Riemann sense, thus a proper subset of
$B(X)$. Consider $f:= \Eins_C$ where
$C = \left\{ \sum_{n=1}^{\infty} \tfrac{a_n}{3^n}:\ a_n \in \{0,2\} \right\}$
is Cantor's middle third set. Then $f \in \Riemann_p$, but
$f \notin B(\Algebra)$: $f \in B(\Algebra)$ would yield the existence of
$f_1 = \sum_{i=1}^n c_i \Eins_{A_i} \in \overline{\Simple_{\Algebra}}$
with $A_i \in \Algebra$, $c_i \in \C$ and $||f-f_1||_{\infty} < \tfrac 12$.
We may assume that the $A_i$ are pairwise disjoint.
Consider $A:= \bigcup_{i: |c_i-1| < \frac 12} A_i \in \Algebra$ and
$f_2:= \Eins_A \in \Simple_{\Algebra}$. Then
$||f-f_2||_{\infty} < \tfrac 12$ which, since $f$ and $f_2$ only take
the values 0 and 1, implies $f=f_2$ and $C=A$, a finite union
of intervals, contradiction.
\end{exa}

We have seen that each fapm $p$ on a set algebra in a natural way
induces a mean $m$ on the $C^*$-algebra $\Riemann_p$. Recall from
the first statement in Proposition \ref{Psetsversusfunctions} that
$\Algebra_{\A} := \{A \subset X:\ \Eins_A \in \A\}$ is a set algebra
whenever $\A$ is an algebra of functions. Given a mean $m$ on $\A$,
$p_m(A):= m(\Eins_A)$ clearly defines a fapm on
$\Algebra_m:=\Algebra_{\A}$. We ask whether the constructions
$\varphi: (\Algebra,p) \mapsto (\Riemann_p,m_p)$ and $\psi: (\A,m)
\mapsto (\Algebra_m,p_m)$ are inverse to each other. In general this
is not the case.

\begin{exa}\label{Ealgebrasets}
Consider any algebra $\A$ of
continuous functions on a nontrivial
connected space $X$ (for instance $X=[0,1]$)
containing functions which are not constant,
and any nontrivial mean $m$ on $\A$.
Then $\Algebra_{\A} = \{\emptyset,X\}$ and hence $\Riemann_{p_m}$
only contains the constant functions and does not coincide with $\A$.
\end{exa}

However, this is not surprising if we note that $\A$ in the above example
is not $m$-closed, while $\Riemann_p$ is $m_p$-closed. Thus we have
to assume this property for all function algebras and means,
and to use the analogue property for fapm's.

\begin{defi}\label{Dcompletesetalgebra}
Consider a fapm $p$ on a set algebra $\Algebra$ on the set $X$. Then the
\textbf{$p$-completion} $\overline{\Algebra}^{(p)}$ of $\Algebra$
is defined as the set of all $A \tm X$ with the following property:
For each $\e>0$ there are $A_1,A_2 \in \Algebra$ with $A_1 \tm A \tm A_2$
and $p(A_2 \setminus A_1) < \e$. For $A\in \overline{\Algebra}^{(p)}$ we define
$p(A)$ to be the unique $\alpha$ with $p(A_1)\le \alpha\le p(A_2)$ for all $A_1, A_2 \in \Algebra$
with $A_1\tm A \tm A_2$. In this way we canonically extend $p$ to all of $\overline{\Algebra}^{(p)}$.
In the case
$\overline{\Algebra}^{(p)} = \Algebra$ we call $p$ \textbf{complete} and
$\Algebra$ \textbf{$p$-closed}.
\end{defi}

It is clear that the $p$-completion of a set algebra is again a set algebra.
Note furthermore that in the case that $p$ is $\sigma$-additive the notion
coincides with the usual concept of a complete measure.

\begin{pro}\label{Pcompletesetalgebra}
Let $\Algebra$ be a set algebra on $X$ and $p$ a fapm on $\Algebra$.
\begin{enumerate}
\item
$\Algebra \tm \Algebra_{m_p}$ and $p_{m_p}(A)=p(A)$ whenever $A \in \Algebra$.
\item $\overline{\Algebra}^{(p)}=\Algebra_{m_p}$. In particular the equality $\Algebra = \Algebra_{m_p}$ holds if and only if
$\Algebra$ is $p$-closed.
\end{enumerate}
\end{pro}
\begin{proof}
The first statement is obvious. To prove the second statement assume first
that $A\in \Algebra_{m_p}$ and pick any $\eps>0$. Then we have, by definition of $\Algebra_{m_p}$, that  $\Eins_A \in \Riemann_p$. By definition of $\Riemann_p$ this means that there are $f_1,f_2 \in \Simple_{\Algebra}$ such that $f_1 \le \Eins_A \le f_2$ and $m_p (f_2-f_1) < \e$. There is a representation $f_2-f_1 = \sum_{i=1}^n c_i \Eins_{A_i}$ such that the $A_i$ are nonempty, pairwise disjoint and both  $f_1$ and $f_2$ are constant on each $A_i$. $f_2-f_1 \ge 0$ implies $c_i \ge 0$ for all $i$.
  Consider the partition of $\{1,\ldots,n\}$
  into three sets $I_1,I_2,I_3$ in such a way that $A_i \tm A$
  for $i \in I_1$ and $A_i \cap A = \emptyset$ for $i \in I_2$.
  For $i \in I_3$ we require that $A_i$ intersects $A$ as well as
  $X \setminus A$. We define $B_1:= \bigcup_{i \in I_1} A_i$
  and $B_2:= B_1 \cup \bigcup_{i \in I_3} A_i$, hence $B_1 \tm A \tm B_2$
  and $B_1,B_2 \in \Algebra$. Note that $f_1 \le \Eins_A \le f_2$
  together with the fact that the $f_1$ and $f_2$ are constant on each
  $A_i$ implies that for $i \in I_3$ we have $f_1 \le 0$ and $f_2 \ge 1$,
  therefore $c_i \ge 1$. We conclude
  \begin{eqnarray*}
       p(B_2 \setminus B_1) = \sum_{i \in I_3} p(A_i)
       &\le&\sum_{i \in I_3} c_i p(A_i) = m_p
       \left(\sum_{i \in I_3} c_i \Eins_{A_i}\right)\\
            & \le&  m_p \left(\sum_{i=1}^n c_i
       \Eins_{A_i}\right)=m_p(f_2-f_1) < \e.
  \end{eqnarray*}
Since $\e>0$ was arbitrary this implies $A \in \overline{\Algebra}^{(p)}$.

If on the other hand we are given a set $A \in \overline{\Algebra}^{(p)}$ and $\eps>0$, then there exist
$B_1,B_2\in \Algebra$ such that $B_1\tm A \tm B_2$ and $p(B_2\setminus B_1)<\eps$. Passing to the indicator functions
$\Eins_{B_1}\le \Eins_A \Eins_{B_2}$ and noting $\Eins_{B_1}, \Eins_{B_2} \in \Simple_{\Algebra}$ we see that $A\in\Algebra_{m_p}$.
\end{proof}

The analogue statement for the converse construction says that,
given a mean $m$ on a $C^*$-algebra $\A$,
$\A = \Riemann_{p_m}$ if and only if $\A$ is $m$-closed.
Later we will use topological constructions for a proof
of this fact, see Proposition \ref{Pmclosed}.

\section{Integration on compact spaces}\label{Scompintegral}
Throughout this text the notion of compactness always includes the Hausdorff separation axiom.
In this section we assume that $X$ is a compact space.
If $\mu$ is a Borel probability measure on $X$ then
$m=m_{\mu}: f \mapsto \int_X f d\mu$ defines a mean on $\A=C(X)$,
the $C^*$-algebra of all continuous $f: X \to \C$.
One of the main reasons that integration theory is particularly
successful on (locally) compact spaces is that also a converse is true:
Positive functionals induce $\sigma$-additive measures.
This is the content of the celebrated Riesz' Representation Theorem,
which we use in the following version:

\begin{pro}[Riesz]\label{riesz}
Let $X$ be compact and $m$ a mean on $C(X)$.
Then there is a unique regular probability measure
$\mu = \mu_m$ which is the completion of its restriction to
the $\sigma$-algebra of Borel sets on $X$ and such that
$m(f) = \int_X f d\mu$ for all $f \in C(X)$. (Recall that \textbf{ regular}
means that for every $\mu$-measurable $A$ and all $\e>0$ there are closed $F$ and open $G$ with $F \tm A \tm G$ and
$\mu(G \setminus F) < \e$.)
\end{pro}
A proof can be found for instance in Rudin's book \cite{Rudi87}.

On the compact unit interval $X=[0,1]$ the classical Riemann integral
can be taken as a mean $m$ on $C(X)$. Then
the measure $\mu_m$ according to Riesz' Representation Theorem is
the Lebesgue measure on $[0,1]$. Note that in this case
$\A = C(X)$ is not $m_{\mu}$-closed, since all Riemann integrable
functions (essentially by the very definition of the Riemann integral)
are members of the $m$-closure of $\A$ but not necessarily continuous.
Sets $A$ with topological boundary $\partial A$
of measure 0 play an important r\^ole.

\begin{defi}\label{DJordan}
Let $\mu$ be a complete Borel measure on $X$.
A set $A \tm X$ is called \textbf{$\mu$-Jordan measurable} or a
\textbf{$\mu$-continuity set} if the topological boundary $\partial A$ of $A$
satisfies $\mu(\partial A) = 0$.
The system of all $\mu$-continuity sets (which forms a set algebra on $X$)
is denoted by $\mathfrak{C}_{\mu}(X)$.
\end{defi}

In the classical case $X = [0,1]$, $\mu$ the Lebesgue measure,
the continuity sets $A$ are exactly those $A \tm [0,1]$ for which $\Eins_A$
is integrable in the Riemann sense. The uniform closure of the
linear span of such $\Eins_A$ coincides with the Riemann integrable functions.
In order to treat the Riemann integral in the context of compactifications
we fix well-known characterizations of classical Riemann integrability
in our somewhat more general context.

For a function $f$, defined on the topological space $X$, we will denote by $\mbox{disc}(f)$ the set of discontinuity points of $f$.

\begin{pro}\label{Priemannchar}
Let $X$ be compact,
$\mu$ a finite complete regular Borel measure
on $X$ and $f: X \to \R$ bounded.
Then the following conditions are equivalent:
\begin{enumerate}
\item\label{discnullset} $\mbox{disc}(f)$ is
$\mu$-measurable and a $\mu$-null set.

\item $f \in \overline{S_{\mathfrak{C}_{\mu}}}=B(\mathfrak{C}_{\mu})$, i.e.\
$f$ can be approximated by simple $\mathfrak{C}_{\mu}$-functions w.r.t\ uniform convergence.

\item\label{contapprox} $f \in \overline{C(X)}^{\,m_{\mu}}$, i.e.\
for every $\eps>0$ there exist $f_1,f_2 \in C(X)$
such that $f_1\le f \le f_2$ and
$\int_X (f_2-f_1) d\mu< \eps$.
\end{enumerate}
If one (and hence all) of these conditions are satisfied, then
$f$ is $\mu$-measurable.
\end{pro}

\begin{proof}
First we prove that condition \ref{discnullset}. implies that $f$ is measurable.
By regularity there is a decreasing sequence of open sets $O_n$, $n \in \N$,
of measure $\mu(O_n) < \tfrac 1n$ with $\mbox{disc}(f) \tm O_n$.
Let $f_n$ be the restriction of $f$ to $X \setminus O_n$. For any
Borel set $B \tm \R$ we have
$f^{-1}[B] = \bigcup_{n \in \N} f_n^{-1}[B] \cup N$ with
$N \tm D:= \bigcap_{n \in \N} O_n$, $\mu(D)=0$. By the completeness of
$\mu$ we conclude that $N$ and thus $f^{-1}[B]$
and finally $f$ is measurable. Now we start with the cyclic proof
of the equivalences.\\

$1 \Rightarrow 2$:
Assume that $\mu(\mbox{disc}(f))=0$ and, w.l.o.g.\, that
$f(X) \tm [0,1]$. We introduce the level-sets
$M_t:=[0\le f<t]$  which are measurable by the first part
of the proof, and the function
\[\phi_f(t):= \mu(M_t).\]
Since $\phi_f$ is increasing, it has at most countably many
points of discontinuity.
Consider
$ \mu(\{x:f(x)=t\})\le\phi_f(r)-\phi_f(s)$ for $s<t<r$.
If $\phi_f$ is continuous at $t$ this implies
\[\sup_{s<t} \phi_f(s)=f(t)=\inf_{r>t} \phi_f(r),\] and so
$\{x:f(x)=t\}$ is a $\mu-$null set for $t\notin \mbox{disc}(\phi_f)$.
Now let $x\in \partial M_t$. If $f$ is continuous at $x$ we clearly
have $f(x)=t$. So
\[\partial M_t \subseteq \mbox{disc}(f)\cup \{x: f(x)=t\}.\]
The first set on the right-hand side is a $\mu-$null set
by our assumption and the second one is a $\mu$-null set at least for
each continuity point $t$ of $\phi_f$. So for all but at most countably
many $t$ the set $M_t$ is a $\mu$-continuity set. In particular
the set $N_f:=\{t: \mu(\partial M_t)=0\}\subseteq [0,1]$ is dense.

Now we approximate $f$ uniformly by members of
$S_{\mathfrak{C}_{\mu}}$: Given $\eps>0$, pick $n\in \N$ such that $n>\frac{1}{\eps}$ and pick real
numbers $\{t_i\}_{i=0}^{n}\subset N_f$ with
\[
t_0=0<t_1<\frac{1}{n}<\ldots<t_{i}<\frac{i}{n}<t_{i+1}<
\ldots<\frac{n-1}{n}<t_{n} = \|f\|_{\infty}\le 1.
\]
Let $A_i:=M_{t_{i}}\!\setminus\! M_{t_{i-1}}$. Then
$|f(x)-\tfrac{i-1}{n}|<\eps$ on $A_i$,  $i=1,\ldots,n$. Since
$X= M_1\setminus M_0=\bigcup_{i=1}^n A_i$ we conclude
\[ \Big|\sum_{i=1}^n \frac{i}{n}
\Eins_{A_i}(x)-f(x)\Big|<\eps.\]

$2 \Rightarrow 3$:
Let $\A_0$ denote the set of all bounded $g: X \to \R$ satisfying
Condition \ref{contapprox}, i.e.\ such that for each $\e>0$ there are $g_1,g_2 \in C(X)$ with $g_1 \le g \le g_2$ and $\int_X( g_2-g_1)\, d\mu < \e$.
It is a routine check that $\A_0$ is a linear space and uniformly
closed. Thus it suffices to show that $\Eins_A \in \A_0$ whenever
$A \in \mathfrak{C}_{\mu}$. For such an $A$ and any given $\e > 0$
we use the regularity of $\mu$ to get an open set $O$ with
$\partial A \tm O$ and $\mu(O) < \e$. Since compact spaces are normal
we can find closed sets $A_1,A_2$ and open sets $O_1,O_2$ with
$$A \setminus O \tm O_1 \tm A_1 \tm A^o \tm \overline A \tm O_2 \tm A_2
\tm A \cup O.$$ Take continuous Urysohn functions $f_1$ for $A \setminus O$ and
$X \setminus O_1$, $f_2$ for $A_2$ and $X \setminus (A \cup O)$,
i.e.\ $$\Eins_{A\setminus O}\le f_1 \le f \le f_2\le\Eins_{A\cup O}.$$
Then $\int_X( f_2 - f_1)\, d\mu \le \mu(O) < \e$.\\

$3 \Rightarrow 1$:
Define the oscillation $\mbox{Os}_f(x)$ of $f$ at a point $x$ by
\[
\mbox{Os}_f(x):=\limsup_{y\to x} f(y)- \liminf_{y \to x} f(y).
\]
Let $A_k:=[\mbox{Os}_f(x) \ge \tfrac 1k]$ be the set of all $x \in X$ where the oscillation of $f$ is at least $\tfrac 1k$. Pick any $\e > 0$ and $k \in \N$. By Condition \ref{contapprox} there are
continuous $f_1^{\e},f_2^{\e}$ with $f_1^{\e} \le f \le f_2^{\e}$
and $\int_X (f_2^{\e} - f_1^{\e})\, d\mu < \tfrac{\e}k$. Note that $A_k
\tm B_k^{\e}:= \{x \in X:\, f_2^{\e}(x)-f_1^{\e}(x) \ge \tfrac 1k\}$ and
$\mu(B_k^{\e}) < 2\e$.  Since $\e>0$ was arbitrary we have $\mu(A_k)=0$.
Since $\mbox{disc}(f) = \bigcup_{k \in \N} A_k$ this proves Condition \ref{discnullset}.
\end{proof}

The equivalence of \ref{discnullset} and \ref{contapprox} can also be found in \cite{Tala82}.

\begin{defi}\label{Driemannmu}
Let $\mu$ be a finite, complete and regular Borel measure on the
compact space $X$ and
$f:X \to \C$ be a bounded function with decomposition
$f=f_1+if_2$ into real and imaginary part. Then $f$ is called \textbf{$\mu$-Riemann integrable} if both $f_1$ and $f_2$ satisfy the
equivalent conditions of Proposition \ref{Priemannchar}. We denote
the set of all $\mu$-Riemann integrable $f$ by
$\Riemanntop_{\mu}(X)$ or $\Riemanntop_{\mu}$.
\end{defi}

The three conditions in Proposition \ref{Priemannchar}
immediately transfer to complex valued functions.

\begin{cor}\label{Criemanncharcomplex}
Let $\mu$ be a finite, complete and regular Borel measure on the
compact space $X$.
For a bounded $f: X \to \C$ the following conditions are equivalent.
\begin{enumerate}
\item $f \in \Riemanntop_{\mu}$, i.e.\ $f$ is $\mu$-Riemann integrable.
\item $\mu(\mbox{disc}(f))=0$.
\item $f \in B(\Simple_{\mathfrak{C}_{\mu}})$.
\end{enumerate}
In particular $\Riemanntop_{\mu} = \Riemann_{p}$ if $p(A):= \mu(A)$
for $A \in \mathfrak{C}_{\mu}$.
\end{cor}

Every $f \in \Riemanntop_{\mu}$ is $\mu$-measurable and the set
$\mbox{disc}(f)$ of discontinuities of a Riemann integrable $f$ is
small not only in the measure theoretic but also in the topological
sense.

\begin{pro}\label{Pfsigma}
Let $X$ be compact and $\mu$ a finite regular
Borel measure, $\mbox{supp}(\mu)=X$ .Let $f\in \Riemanntop_{\mu}(X)$ be Riemann integrable. Then $\mbox{disc}(f)$ is a
meager $\mu$-null set, in particular the set of continuity points of $f$ is dense in $X$.
\end{pro}
\begin{proof}
We may assume that $f \in R_{\mu}(X)$ is real-valued.
It suffices to show that $\mbox{disc}(f)$ is meager.
As in the proof of Proposition \ref{Priemannchar}
let us denote the oscillation of $f$ at $x$ by
$\mbox{Os}_f(x)$. A folkloristic argument shows that the sets
$A_n:=[\mbox{Os}_f \ge \tfrac{1}{n}],\, n>1$ are closed.
The sets $A_n$ are all $\mu$-null sets since $A_n \tm \mbox{disc}(f)$. Using that $\mu$ has full support, this implies that all sets $A_n$ are
nowhere dense, i.e.\ $\mbox{disc}(f)=\bigcup_{n>0 0} A_n$ is a meager
$F_{\sigma}$-set of zero $\mu-$measure.
\end{proof}

We want to illustrate the r\^ole of
the regularity assumption on $\mu$ in Proposition \ref{Priemannchar}.
For this we use the example of a non regular Borel measure occurring
in Rudin's book \cite[Exercise 2.18]{Rudi87}.

\begin{exa}\label{Eregularity}
Let $X=[0,\omega_1]$
be the set of all ordinals up to the first uncountable one
equipped with the order topology. Thus $X$ is a compact space.

We need the fact that every (at most) countable family of
uncountable compact subsets $K_n \tm X$ has an uncountable
intersection $K$. To see this consider any increasing sequence $x_0
< x_1 < x_2 < \ldots \in X$ which meets every $K_n$ infinitely many
times. 
It follows that $\alpha_0 := \sup_n x_n < \omega_1$ is in the closure of all $K_n$, hence in $K$. Since we may require $x_0
> x$ for any given $x < \omega_1$ the same construction can be
repeated in order to obtain an $\alpha_1 \in K$ with $\alpha_1 >
\alpha_0$. Transfinite induction with the limit step
$\alpha_{\lambda}:= \sup_{\nu < \lambda} \alpha_{\nu}$ generates the
closed and thus compact subset of all $\alpha_{\nu}$, $\nu <
\omega_1$, which is contained in $K$.

Easy consequences: We call a set $S \tm X$ of type 1 if $S \cup
\{\omega_1\}$ contains an uncountable compact $K$. If $S$ is of type
1 the complement of $S$ must not have the same property. Call $S \tm
X$ of type 0 if $(X \setminus S) \cup \{\omega_1\}$ contains an
uncountable compact $K$. The system of all
sets of either type 0 or type 1, forms a $\sigma$-algebra $\Algebra$
containing all Borel sets.

Letting $\mu(S)=i$ if $S$ is of Type $i=0,1$, $\mu$ is a complete measure
defined on $\Algebra$. Note that every countable set is a $\mu$-null
set. The set $\{\omega_1\}$ has measure 0 and is a counterexample
for outer regularity: The function $\Eins_{\{\omega_1\}}$ obviously
satisfies conditions 1 and 2 in Proposition \ref{Priemannchar}, but
not condition 3. To see this last assertion consider any continuous $f: X \to \C$ and take
$\beta_n$ such that $|f(x)-f(\omega_1)| < \tfrac 1n$ for all $x \ge
\beta_n$. Then $\beta:= \sup_n \beta_n < \omega_1$ has the property
that $f(x)=f(\omega_1)$ for all $x \ge \beta$. It follows that
$\int_X f\, d\mu = f(\omega_1)$ for all $f \in C(X)$. In particular
$g \le \Eins_{\omega_1} \le h$, $g,h \in C(X)$ implies $\int_X(h-g)\,
d\mu \ge 1$, contradicting condition 3.

Nevertheless we might apply Riesz' Representation Theorem \ref{riesz} to the
functional $m(f):= \int_X f\, d\mu$. A quick inspection shows that
$\mu_m = \delta_{\omega_1}$, i.e.\ the associated unique regular
Borel measure is the point measure concentrated at the point $\omega_1$.
As a complete measure, this $\mu_m$ is defined on the whole
power set of $X$. Finally we observe that
$\Eins_{\omega_1} \notin \Riemann_{\mu_m}$.
\end{exa}

\section{Compactifications and continuity}\label{Scompcont}

The previous Section has illustrated that compactness plays an
important r\^ole in integration theory. This motivates us to
investigate compactifications, the topic of this
purely topological Section. Let $X$ be a, possible discrete, topological space.

We will interpret functions $f:X \to \C$ as
restrictions of functions $F: K \to \C$ on compact spaces $K$.
For our needs the following setting is appropriate.

\begin{defi}\label{Dcomp}
A pair $(\iota,K)$, $K$ compact, $\iota: X \to K$ a continuous mapping, is called a \textbf{compactification}
of $X$ whenever $\overline{\iota(X)}=K$, i.e.\
whenever the image of $X$ under
$\iota$ is dense in $K$.
The function $F: K \to \C$ is called a \textbf{representation} of $f: X \to \C$
whenever $f = F \circ \iota$, i.e.\ whenever the diagram
\begin{diagram}
        &               &K\\
        &\ruTo^{\iota}  &\dTo>{F} \\
X       &\rTo^f         &\C
\end{diagram}
commutes. In this case we also say that $f$ can be represented
in $(\iota,K)$. If $F\in C(K)$ we say that $F$ is a continuous representation.
\end{defi}
Note that in the definition of a compactification $\iota$ is neither required to be a homeomorphic embedding nor to be injective. If there is a continuous representation $F$ of $f$ in $(\iota,K)$, then
this $F$ is uniquely determined by continuity and the fact that $\iota(X)$
is dense in $K$. Furthermore $f=F \circ \iota$ is continuous as well.
In this Section we are therefore mainly interested in continuous $f$.
Let us consider first a rather trivial example.

\begin{exa}\label{Enaturalcomp}\hfill\par\vspace{-2ex}
\begin{itemize}
\item
Let $f: X \to \C$ be bounded and continuous.
Surely $K_f:= \overline{f(X)}$ is compact. Define
$\iota_f: x \mapsto f(x)$ and let $F_f: K_f \to \C$ be the inclusion mapping. Then $(\iota_f,K_f)$ is a compactification of $X$ and $F_f$ is a continuous representation of $f$ in $(\iota_f,K_f)$. We call $F_f$ the \textbf{natural continuous representation of $f$}.
\item
Let $f: X\to \C$ be merely bounded. If we impose the discrete topology on $X$, $f$ is continuous and the associated compactification
$(\iota_f,K_f)$ is then a compactification of the discrete space $X_{\mbox{\tiny dis}}$.
\end{itemize}

\end{exa}

One observes the following minimality property of the natural continuous representation:
If $F: K \to \C$ is any continuous representation of $f$ in any
compactification $(K,\iota)$ of $X$, then $\pi: K \to K_f = \overline{f(X)}$,
$\pi(k):= F(k)$, is continuous, onto and satisfies $\pi \circ \iota = \iota_f$.
This motivates the following definition.

\begin{defi}\label{Dcompvergleich}
Let $(\iota_1,K_1)$ and $(\iota_2,K_2)$ be two compactifications of $X$.
Then we write $(\iota_1,K_1) \le (\iota_2,K_2)$ (via $\pi$) and
say that $(\iota_1,K_1)$ is smaller than $(\iota_2,K_2)$
or, equivalently, $(\iota_2,K_2)$ is bigger than $(\iota_1,K_1)$,
if $\pi: K_2 \to K_1$ is continuous satisfying $\iota_1 = \pi \circ \iota_2$,
i.e.\ making the diagram
\begin{diagram}
        &               &       K_2\\
        &\ruTo^{\iota_{2}}&\dTo>{\pi} \\
X       & \rTo^{\iota_{1}}&        K_1\\
\end{diagram}
commutative. For the case that $\pi$ is a homeomorphism we
say that $(\iota_1,K_1)$ and $(\iota_2,K_2)$ are
equivalent via $\pi$ and write $(\iota_1,K_1) \cong (\iota_2,K_2)$.
\end{defi}

A consequence of the continuity of the involved maps and of the
fact that the images $\iota_i(X)$ are dense is that $\pi$ as in
Definition \ref{Dcompvergleich} is unique. By compactness, $\pi$
is onto as well. If $\pi$ happens to
be injective it is a homeomorphism, i.e.\ $(\iota_1,K_1)$ and
$(\iota_2,K_2)$ are equivalent. Furthermore one easily sees that,
whenever $(\iota_1,K_1) \le (\iota_2,K_2)$ via $\pi_1$ and
$(\iota_2,K_2) \le (\iota_1,K_1)$ via $\pi_2$ then
$\pi_2 \circ \pi_1$ is the identity on $K_1$ and $\pi_1 \circ \pi_2$
is the identity on $K_2$, hence $\pi_2 = \pi_1^{-1}$, $\pi_1$ and $\pi_2$
are isomorphisms and both compactifications are equivalent.

\begin{pro}\label{Pcompequiv}
$(\iota_1,K_1) \cong (\iota_2,K_2)$ if and only if both,
$(\iota_1,K_1) \le (\iota_2,K_2)$ and $(\iota_2,K_2) \le (\iota_1,K_1)$.
\end{pro}

Note that maps $\pi_1,\pi_2$ as in Definition \ref{Dcompvergleich}
may as well be considered to be the morphisms in a
category whose objects are all compactifications of $X$.
Other related categories arise if one allows only continuous representations
of one fixed $f: X \to \C$.
In this terms the minimality property of the natural compactification
asserts that $(\iota_f,K_f)$ is a universal object and thus
unique up to equivalence.

\begin{pro}\label{Pcomprepr}
Let $F_1$ be a representation of $f: X \to \C$ in a compactification
$(\iota_1,K_1)$ of $X$, and suppose $(\iota_1,K_1) \le (\iota_2,K_2)$
via $\pi$. Then $F_2:= F_1 \circ \pi$ is a representation of $f$ in
$(\iota_2,K_2)$ which is continuous whenever $F_1$ is continuous.
\end{pro}

Given a family of compactifications $(\iota_i,K_i)$, $i \in I$, of $X$, we get a common upper bound by taking products:
Let $\iota(x) := (\iota_i(x))_{i \in I} \in P:= \prod_{i \in I} K_i$ and
$K:= \overline{\iota(X)} \tm P$.
Then one obtains a compactification $(\iota,K)$
which, by the projections $\pi_{i_0}: K \to K_{i_0}$, $i_0 \in I$,
$(k_i)_{i \in I} \mapsto k_{i_0}$, indeed satisfies
$(\iota_i,K_i) \le (\iota,K)$ for all $i \in I$. Sometimes we use
the notation $\bigvee_{i\in I} (\iota_i,K_i)$ for $(\iota,K)$.

\begin{defi}\label{Dprodcomp}
For compactifications $(\iota_i,K_i)$ of $X$, $i \in I$, the
compactification $(\iota,K)$, $\iota: x \mapsto (\iota_i(x))_{i_\in I}$,
$K:= \overline{\iota(X)} \tm \prod_{i \in I}K_i$, is called the
\textbf{product compactification} of all $(\iota_i,K_i)$, $i \in I$.
\end{defi}

\begin{pro}\label{Pproduktcompisminimal}
For compactifications $(\iota_i,K_i)$ of $X$, $i \in I$, the supremum
$\sup_{i\in I} (\iota_i,K_i)$ is equivalent to the product compactification
$(\iota,K)$ of all $(\iota_i,K_i)$, $i \in I$.
\end{pro}
\begin{proof}
We have already seen that $\sup_{i\in I} (\iota_i,K_i) \le (\iota,K)$. Let $(\iota', K')$ be another
compactification of $X$ such that $(\iota_i,K_i)\le (\iota', K'),\, i\in I$. Denote by $\pi_i: K' \to K_i$ the $i$-th projection. Define a mapping $\pi: K' \to K$ via $k'\mapsto (\pi_i(k'))_{i\in I}$. Note that
$\pi \circ \iota' = \iota$,
hence $\pi(\iota'(X)) \tm K$ and
\[
\pi(K') = \pi(\overline{\iota'(X)}) \tm \overline{\pi(\iota'(X))} \tm K.
\]
It is immediate to check that $\pi$ is continuous and satisfies $\pi\circ\iota'=\iota$; thus $(\iota,K)\le (\tilde\iota,\tilde K)$.
\end{proof}

Analogously the product compactification can be used to obtain
a minimal compactification where all functions from an arbitrary
given family have a continuous representation: Let $f_i: X \to \C$, $i \in I$, be bounded and continuous functions on $X$. We consider the natural continuous representations of the $f_i$, i.e.\ $(\iota_i,K_i):= (\iota_{f_i},K_{f_i})$ and $F_i: K_i \to \C$, the inclusion mappings. Let $(\iota,K)$ be the product of all $(\iota_i,K_i)$, $i \in I$.

\begin{defi}\label{Dprodnaturalcomp}
Let us denote the $C^*$-algebras of bounded resp.\ continuous resp.\ bounded and continuous $f:X \to \C$ by $B(X)$, $C(X)$ resp.\ $C_b(X)$. For a given family of $f_i \in C_b(X)$, $i \in I$, the compactification $(\iota,K)$, constructed as above is called the \textbf{natural compactification} for the family of all
$f_i$, $i \in I$. If $\A = \{f_i:\ i \in I\}$ we also
write $(\iota_{\A},K_{\A})$ for $(\iota,K)$.
\end{defi}

\begin{pro}\label{Palgebracontinuous}
Let $\A\subseteq C_b(X)$, then the following holds
\begin{enumerate}
\item Every $f \in \A$ has a continuous representation in
  the natural compactification $(\iota_{\A},K_{\A})$ of $\A$.
\item Suppose that $(\iota,K)$ is any compactification of $X$
  where every $f \in \A$ has a continuous representation.
  Then $(\iota_{\A},K_{\A}) \le (\iota,K)$, i.e. $(\iota_{\A},K_{\A})$ is minimal among the compactifications with this property .
\item $\{F\circ\iota: F \in C(K_{\A})\}$ is a $C^*$-algebra and the $*$-algebra generated by $\A$ is dense in this $C^*$-algebra. In particular, if $\A$ is $C^*$-algebra, then $\A$ contains exactly those $f$ which have a continuous
 representation in $(\iota_{\A},K_{\A})$.
\end{enumerate}
\end{pro}
\begin{proof}
\begin{enumerate}
\item
For each $i_0 \in I$, $G_{i_0}: K \to \C$, $(c_i)_{i \in I} \mapsto c_{i_0}$,
is continuous and satisfies $f_{i_0} = G_{i_0} \circ \iota$ for
each $i_0 \in I$. Thus all $f_i$ can be continuously represented
in $(\iota_{\A},K_{\A})$.
\item
Let $(\iota',K')$ be an arbitrary compactification
of $X$ where continuous representations $G'_i: K' \to \C$ of
$f_i = G'_i \circ \iota'$ exist. As in Proposition \ref{Pproduktcompisminimal} we define
$\pi: k' \mapsto (G'_i(k'))_{i \in I} \in \prod_{i \in I} \overline{f_i(X)}$.
$\pi$ is continuous because all components are.
Again we have $\pi(K') \tm K$, $\pi: K' \to K$ and $(\iota,K) \le (\iota',K')$.
Furthermore $G'_i = G_i \circ \pi$ for all $i \in I$, since
the mappings on both sides are continuous and coincide on the
dense set $\iota'(X)$.
\item It is clear that the mapping $F\mapsto F\circ\iota$ maps the $C^*$-algebra $C(K_{\A})$ again
on a $C^*$-algebra and that this map is a continuous homomorphism between $C^*$-algebras.

For the rest of the proof we can assume w.l.o.g. that $\A$ is a $*$-algebra. It remains to prove that the $*$-algebra $\A':=\{F\in C(K_{\A}): F\circ\iota\in \A\}$ is dense in $C(K_{\A})$. We employ the
Stone-Weierstra{\ss} theorem. Obviously $\A'$ is a
$*$-algebra containing all constant functions. We are done
if $\A'$ is point separating. Pick $c \neq c' \in K_{\A}$. Recall that
the points in $K_{\A}$ are of the form $c=(c_f)_{f \in \A}$ and
$c'=(c'_f)_{f \in \A}$ with $c_f,c'_f \in \C$. Hence there is some
$f_0 \in \A$ such that $c_{f_0} \neq c'_{f_0}$. By definition, $K_{\A}$
is the closure of the set of all $(f(x))_{f \in \A}$, $x \in X$. It
follows that there are $x,x' \in X$ with $f_0(x)$ arbitrary close to
$c_{f_0}$, $f_0(x')$ to $c'_{f_0}$, hence $f_0(x) \neq f_0(x')$. Let
$F_0 = \pi_{f_0} \in C_b(X)$ implying $f_0 = F_0 \circ \iota$ and $F_0
\in \A'$ with $F_0(c) = f_0(x) \neq f_0(x') =F_0(c')$. Thus $\A'$ is
indeed point separating, which completes the proof.
\end{enumerate}
\end{proof}

\begin{pro}\label{Palgsepa}
Let $\A$ be a $C^*$-algebra on $X$. Then $\A$ separates points of $X$ if and only if in the natural compactification $(\iota_{\A},K_{\A})$ the map $\iota_{\A}: X\to K_{\A}$ is one-one.
\end{pro}
\begin{proof}
Recall that $\iota_{\A}(x):=(f(x))_{x\in \A}$. $\A$ separates points of $X$ if and only if for all $x_1,x_2 \in X$ with $x_1\neq x_2$ there exists $f\in \A$ such that $f(x_1)\neq f(x_2)$, i.e. $\iota_{\A}(x_1)\neq\iota_{\A}(x_2)$.
\end{proof}

\begin{cor}[Gelfand]\label{Calgebracomp}
The mapping $\A \mapsto (\iota_{\A},K_{\A})$ is (modulo equivalence
of compactifications) a bijective and order-preserving
correspondence between compactifications of $X$ and
$C^*$-subalgebras of $C_b(X)$ which contain $\Eins_X$. In particular $\A$ and $C(K_{\A})$ are isomorphic as
$C^*$-algebras.
\end{cor}

\begin{rem} Note that Corollary \ref{Calgebracomp} applies to $C^*$-subalgebras of $B(X)$ as well. All one has to do is to identify $B(X)$ with $C_b(X_{\mbox{\tiny dis}})$. Thus $B(X)$ is an $C^*$-algebra of continuous functions.
\end{rem}

\begin{exa}\label{EStonespace}
 Let us consider the special case that $\A = B(\Algebra_{\A})\subseteq C_b(X_{\mbox{\tiny dis}})$
and write $\Algebra = \Algebra_{\A}$. We consider the set
$\A_1:= \{\Eins_A:\ A \in \Algebra\}$, the corresponding compactification
$(\iota_1,K_1):= (\iota_{\A_1},K_{\A_1})$ and the commutative diagram

\begin{diagram}
        &               &       K_{\A}\\
        &\ruTo^{\iota_{\A}}&\dTo>{\pi} \\
X       & \rTo^{\iota_{\A_1}}&        K_1\\
\end{diagram}
with  $\pi: (c_f)_{f \in \A} \mapsto (c_f)_{f \in \A_1} \in \{0,1\}^{\A_1}$.
We claim that $\pi$ is injective. Suppose first $\Eins_A(x) = \Eins_A(y)$
for all $A \in \Algebra$. Then $f(x)=f(y)$ for all $f \in \Simple_{\Algebra}$
and hence for all $f$ from the closure $B(\Algebra)=\A$. Suppose now
that $c=(c_f)_{f \in \A_1}=\pi(a)=\pi(b) \in \overline{\iota_{\A_1}(X)}$
with $a=(a_f)_{f \in \A}$ and $b=(b_f)_{f \in \A}$.
Then $a_f=b_f=c_f$ for all $f \in \A_1$.
There is a net $(x_{\nu})_{\nu \in N}$, $N$ a directed set, such that
$\iota_{\A_1}(x_{\nu}) \to c$.
Define $\iota_{\A_1}(x_{\nu}) = (c_f^{\nu})_{f \in \A}$.
Note that $c_f^{\nu} = f(x_{\nu})$. Thus we have $f(x_{\nu}) \to c_f=a_f=b_f$
for all $f \in \A_1$, hence, by linearity, for all $f \in \Simple_{\Algebra}$
and, by uniform closure, for all $f \in B(\Algebra) = \A$.
Therefore we conclude that
$\iota_{\A}(x_{\nu}) = (f(x_{\nu}))_{f \in \A} \to a
=(a_f)_{f \in \A} = (b_f)_{f \in \A} = b$, proving that $\pi$ is injective.
Thus $(\iota_{\A},K_{\A}) \cong (\iota_{\A_1},K_{\A_1})$.
A clopen subbasis of $K_1$ is given by all sets
$A_0':= \{(c_A)_{A \in \Algebra}:\ c_{A_0}=1\}$, $A_0 \in \Algebra$.
\end{exa}

\begin{cor}[Stone]\label{CStonespace}
If $\A = B(\Algebra_{\A})$ then the compact space $K_{\A}$ is totally disconnected.
\end{cor}

Note that the natural context of our discussion are classical theorems
due to Gelfand, Banach and Stone. Without going into formal details
these results are as follows. Gelfand's representation theorem
states that every (abstract) commutative unital $C^*$-algebra $\A$
(meaning that complex conjugation is replaced by an abstract
operation with corresponding properties)
is isometrically isomorphic to some $C(K)$ where $K$ is a suitable
compact space. In this context $K$ is also called
the structure space or Gelfand compactum for $\A$. By the Banach-Stone
theorem, two compact spaces $K_1$ and $K_2$ are homeomorphic
if and only if $C(K_1) \cong C(K_2)$ as unital Banach algebras.
Furthermore, by Stone's theorem, for every Boolean algebra $B$ there is a
totally disconnected compact space $K$, the so called Stone space
associated to $B$, such that for the
systems $\mbox{Cl}(K)$ of all clopen subsets of $K$ we have $B \cong \mbox{Cl}(K)$
as Boolean algebras. Two such spaces $K_1$ and $K_2$ are homeomorphic
if and only if $\mbox{Cl}(K_1) \cong \mbox{Cl}(K_2)$. Finally, the Stone space
of a Boolean set algebra $\Algebra$ is homeomorphic to the Gelfand
compactum for $B(\Algebra)$. For the interested reader we refer
to \cite{Davi96} and \cite{DuSc88}.

\section{The Stone-{\v C}ech compactification $\beta X$}\label{Sstonecech}

We now apply the construction of the natural compactification for an algebra $\A$, to the case $\A = C_b(X)$, i.e.\ to the algebra of all bounded and continuous $f: X \to \C$.

\begin{defi}\label{Dstonecech}
The maximal compactification $(\iota_{\beta},\beta X)$ of
a topological space $X$, corresponding to the algebra $C_b(X)$ in
the sense on Corollary \ref{Calgebracomp}, is denoted by
$(\iota_{\beta},\beta X)$ and is called the \textbf{Stone-{\v C}ech
compactification} of $X$.
\end{defi}

$(\iota_{\beta},\beta X)$ is characterized uniquely up to
equivalence by the universal property that for every continuous
$\varphi: X \to K$, $K$ compact, there is a (unique) continuous
$\psi: \beta X \to K$ with $\varphi = \psi \circ \iota_{\beta}$.
To see this we may w.l.o.g.\ assume $K = \overline{\varphi(X)}$
such that $(\varphi,K)$ is a compactification of $X$.
By the maximality of $(\iota_{\beta},\beta X)$ and Corollary
\ref{Calgebracomp} this just means that there is a $\psi$ as claimed.
For uniqueness assume that $(\iota,K)$ is another compactification
of $X$ with this universal property. Every $f \in C_b(X)$ has a range
contained in a compact set $K_0 \tm \C$. By the universal property
there is a continuous $\psi: K \to K_0$ with $\psi \circ \iota = f$.
Hence, again by Corollary \ref{Calgebracomp} the algebra corresponding
to $(\iota,K)$ contains $C_b(X)$. Thus $(\iota,K)$ has to be maximal,
i.e.\ equivalent $(\iota_{\beta}, \beta X)$.

Nevertheless, in order to obtain an interesting and rich structure one needs
sufficiently many bounded and continuous functions.

\begin{defi}\label{Dcompletelyregular}
$X$ is called \textbf{completely regular} if it
fulfills the following separation property: For every closed $A \tm X$ and $x \in X \setminus A$ there is a continuous $f: X \to [0,1]$ with $f(x)=1$ and $f(a)=0$ for all $a \in A$. Such an $f$ is called Urysohn function for $A$ and $x$.
\end{defi}

Under this assumption every Urysohn function gives rise to a compactification
separating two points $x \neq y \in X$, yielding that
$\iota_{\beta}$ is injective. $\iota_{\beta}$ is even a homeomorphic
embedding of $X$ into $\beta X$. To see this, it suffices to show
that for $x \in O \tm X$, $O$ open, there is an open set
$O_{\beta} \tm \beta X$ containing $\iota_{\beta}(x)$ such that
$\iota_{\beta}(O) \supseteq O_{\beta} \cap \iota_{\beta}(X)$.
Take a Urysohn function $f_0$ for $x$ and $A:= X \setminus O$,
recall that $\iota_{\beta}: x \mapsto (f(x))_{f \in C_b(X)}$ and
observe that $O_{\beta}:= \{(c_f)_{f \in C_b(X)}:\ c_{f_0}>0\}$
has the required properties.

Let us now consider the case of discrete $X$, i.e.\ $C_b(X)=B(X)$. Then each $\Eins_A \in B(X)$, $A \tm X$, has a continuous representation
in $(\iota_{\beta},\beta X)$ which must be of the form $\Eins_{A^*}$
with some clopen $A^* = \overline{\iota_{\beta}(A)} \tm \beta X$.
(Therefore the usual notation $A^* = \overline A$ as a closure,
though not rigorously correct in our setting,
does not lead to contradictions.)

Conversely, every clopen set $B \tm \beta X$ can be written as $B=A^*$ with $A:= \iota_{\beta}^{-1}[B]$.
Furthermore such sets form a basis for the topology in $\beta X$: Let $O \tm \beta X$ be open and $x \in O$. Then, by the separation
properties of compact spaces, there is an open set $O_x$ such that
$x \in O_x \tm \overline{O_x} \tm O$. For $A_x:= \iota_{\beta}^{-1}[O_x]$
we obtain $x \in A_x^* \tm O$. This shows that $O = \bigcup_{x \in O}A_x^*$ can be written as a union of clopen sets.

Let $A=\{a\}$, $a \in X$, be a singleton and $x \neq \iota_{\beta}(a)$.
There is an open neighborhood $O$ of $X$ not containing $\iota_{\beta}(a)$.
Thus the continuous representation of $\Eins_A$ in $(\iota_{\beta}, \beta X)$
has to take the constant value $0$ on $O$, hence
$\Eins_{A^*} = \Eins_{\{\iota_{\beta}(a)\}}$. By continuity this shows
that $\{\iota_{\beta}(a)\}$ is open, i.e.\ $\iota_{\beta}(a)$ is
an isolated point in $\beta X$. A further consequence is that
$A^* \cap \iota_{\beta}(X \setminus A) = \emptyset$ and
$A^* \cap (X \setminus A)^* = \emptyset$. Since
\[
\beta X = \overline{\iota_{\beta}(X)} =   \overline{\iota_{\beta}(A)
  \cup \iota_{\beta}(X \setminus A)} =  \overline{\iota_{\beta}(A)}
  \cup \overline{\iota_{\beta}(X \setminus A)} = A^* \cup (X \setminus A)^*
\]
we conclude that $\Phi: A \mapsto A^*$ is an isomorphism of
Boolean set algebras between  $\P(X)$, the powerset of $X$, and $\mbox{Cl}(\beta X)$, the system of all clopen sets in $\beta X$.

Consider $\F_x:= \{\iota_{\beta}^{-1}[O]:\ x \in O \tm \beta X,
O \ \mbox{open}\}$. Obviously $\F_x$ is a filter on $X$. For arbitrary
$A \tm X$, by $A^* \cup (X \setminus A)^* = \beta X$, we have
either $x \in A^*$ or $x \in (X \setminus A)^*$. In the first case
this means $A = \iota_{\beta}^{-1}[A^*] \in \F_x$, in the second case
$X \setminus A \in \F_x$. Thus $\F_x$ is an ultrafilter.
Conversely every ultrafilter $\F$ on $X$ induces an ultrafilter
$\F_{\beta}$ on $\beta X$ consisting of all $F_{\beta} \tm \beta X$
which contain $\iota_{\beta}(F)$ for at least one $F \in \F$.
The compactness of $\beta X$ guarantees that $\F_{\beta}$ converges
to some $x \in \beta X$ which is possible only if $\F = \F_x$.
This shows that the points in $\beta X$ are in a natural bijective
correspondence with the ultrafilters on $X$.

We summarize the collected facts about $\beta X$.

\begin{pro}\label{Pstonecechsummary}
Let $X$ be a completely regular topological space.
Then the Stone-{\v C}ech compactification $(\iota_{\beta},\beta X)$ of $X$
has the following properties.
\begin{enumerate}
\item
For every continuous $f: X \to K$, $K$ compact, there is a continuous
$\varphi: \beta X \to K$ with $f = \varphi \circ \iota_{\beta}$,
i.e.\ making the diagram
\begin{diagram}
        &                    &       \beta X\\
        &\ruTo^{\iota_{\beta}}&\dTo>{\phi} \\
X       & \rTo^{f}           &        K\\
\end{diagram}
commutative.
\item
$\iota_{\beta}: X \to \iota_{\beta}(X) \tm \beta X$ is a homeomorphism.
\item
Assume that $X$ is discrete.
  \begin{enumerate}
  \item The mapping $A \mapsto A^*:= \overline{\iota_{\beta}(A)}$
    is an isomorphism of Boolean set algebras between  $\P(X)$, the powerset of $X$, and $\mbox{Cl}(\beta X)$, the system of all clopen sets in $\beta X$.
  \item The clopen subsets of $\beta X$ form a topological basis in $\beta X$.
  \item The isolated points in $\beta X$ are exactly those of the form
    $\iota_{\beta}(x)$, $x \in X$.
  \item $\beta X$ can be represented as the set of all ultrafilters
    on $X$ where $\iota_{\beta}(x) = \F_x:= \{F \tm X:\ x \in F\}$
    for all $x \in X$. Then $A^*$ consists of those ultrafilters $\F$
    on $X$ with $A \in \F$.
\end{enumerate}
\end{enumerate}
\end{pro}

\section{Compactifications, measures, means and Riemann integral}
\label{Scompriemann}

We are now going to consider compactifications $(\iota,K)$
of a set (or a topological space) $X$ in connection with complete
Borel probability measures $\mu$ on $K$.

\begin{defi}\label{Dmeasurecomp}
Let $(\iota,K)$ be a compactification of $X$, $\mu$ a complete and regular Borel probability
measure on $K$ and $\A$ a set of
complex valued $\mu$-measurable functions on $K$. Then we call the
quadruple $(\iota,K,\mu,\A)$ \textbf{admissible} if the following condition is satisfied: Whenever $F_1 \circ \iota = F_2 \circ \iota$
for $F_1,F_2 \in \A$ then
\[\int_K F_1\, d\mu = \int_K F_2\, d\mu.\]
For arbitrary $\A$ define $\A^*:= \iota^*(\A) = \{F \circ \iota:\ F \in \A\}$.
For admissible $(\iota,K,\mu,\A)$ we define
\[m_{\mu}:\; f=F \circ \iota \mapsto \int_K F\, d\mu, \quad F \in \A.\]
\end{defi}

Note that $m_{\mu}$ is well defined on $\A^*$ and a bounded linear functional whenever $\A$ is a linear space, called the \textbf{mean induced by $(\iota,K,\mu,\A)$}.

It is clear that for all compactifications $(\iota,K)$ of $X$ and
all $\mu$ we get an admissible quadruple if we take $\A:= C(K)$.
In this case $F_1 \circ \iota = F_2 \circ \iota$ with $F_1,F_2 \in \A$ is possible only for $F_1=F_2$ (recall that $\iota(X)$ is dense in $K$).
For our subsequent investigations the following similar statement
for $\A = \Riemanntop_{\mu}$ is fundamental.

\begin{pro}\label{Priemannmean}
For every compactification $(\iota,K)$ of $X$ and every complete and regular Borel probability
measure on $\mu$ on $K$ the quadruple
$(\iota,K,\mu,\Riemanntop_{\mu})$ is admissible. Hence
\[m(F \circ \iota):= \int_K F\, d\mu \]
is a well defined mean on the algebra $\Riemanntop_{\mu}^*$.
\end{pro}
\begin{proof}
W.l.o.g.\ we may assume that $\mu$ has full support,
i.e.\ all nonempty open sets in $K$ have positive measure.
Let $f=F_1\circ\iota=F_2\circ\iota$ with $F_i\in \Riemanntop_{\mu}$.
First we assume that $F_i = \Eins_{A_i}$ for certain $\mu$-continuity sets
$A_i \in \mathfrak{C}_{\mu}$. The symmetric difference
$A:= A_1 \bigtriangleup A_2$ is a $\mu$-continuity set with empty
interior. We conclude that $\partial A \tm \overline A$ has zero measure
and hence $\int \Eins_{A_1} d\mu = \int \Eins_{A_2} d\mu$.
By linearity this property extends to functions $F_i \in \Simple_{\mathfrak{C}_{\mu}}$
and, using a standard approximation argument,
to arbitrary $F_i \in \Riemanntop_{\mu}$.
\end{proof}

We have to compare compactifications also in a measure theoretic sense.
For this reason we fix further notation.

\begin{defi}\label{Dcompmeasure}
Suppose that $\mu_i$ is a complete Borel probability measure on $K_i$, where
$(\iota_i,K_i)$ is a compactification of $X$, $i=1,2$. Then we write
$(\iota_1,K_1,\mu_1) \le (\iota_2,K_2,\mu_2)$ if
$(\iota_1,K_1) \le (\iota_2,K_2)$
via $\pi: K_2 \to K_1$
which, in addition to being continuous is also measure preserving, i.e.\ whenever $A_1 \tm K_1$
is $\mu_1$-measurable then its pre-image $A_2:=\pi^{-1}[A_1]$ is
$\mu_2$-measurable with $\mu_2(A_2)=\mu_1(A_1)$.
\end{defi}

\begin{rem} In the above situation we also could have \emph{defined}
the measure on $K_1$ via $\mu_1(A_1):=\mu_2(\pi^{-1}[A_1])$. This construction is called
pullback, and $\mu_1$ is often denoted by $\pi\circ \mu_2$.
\end{rem}

We know by Proposition \ref{Pcomprepr}
that every $f: X \to \C$ which has a continuous representation
$F_1: K_1 \to \C$ in the compactification $(\iota_1,K_1)$ has a continuous
representation $F_2:= F_1 \circ \pi$ in $(\iota_2,K_2)$ whenever
$(\iota_1,K_1) \le (\iota_2,K_2)$ via $\pi: K_2 \to K_1$.
We get a similar assertion if we replace continuity by Riemann integrability.

\begin{pro}\label{Pcompmeasure}
Suppose that $f: X \to \C$ has a $\mu_1$-Riemann integrable representation
$F_1: K_1 \to \C$ in the compactification $(\iota_1,K_1,\mu_1)$. Whenever $(\iota_1,K_1,\mu_1) \le (\iota_2,K_2,\mu_2)$ via $\pi$ then $F_2:= F_1 \circ \pi$ is a $\mu_2$-Riemann integrable representation
of $f$ in $(\iota_2,K_2,\mu_2)$.
\end{pro}

\begin{proof}
It is clear that $F_2:=F_1\circ\pi$ is a realization of $f$ whenever $F_1$ is. It is immediate to check that
$\mbox{disc}(F_2\circ\pi)\subseteq \pi^{-1}[\mbox{disc}(F_1)]$.
Thus one obtains
\[
{\mu_2}(\mbox{disc}(F_2)\le {\mu_2}(\pi^{-1}[\mbox{disc}(F_1)])=
{\mu_1} (\mbox{disc}(F_1))=0.\]
Thus $F_2$ is $\mu_2$-Riemann integrable whenever
$F_1$ is $\mu_1$-Riemann integrable.
\end{proof}

Proposition \ref{Pcompmeasure} shows that
$(\iota_1,K_1,\mu_1) \le (\iota_2,K_2,\mu_2)$ implies
$\{F_1\circ\iota: F_1\in \Riemanntop_{\mu_1}\}\subseteq
\{F_2\circ\iota: F_2\in \Riemanntop_{\mu_2}\}$. This observation is of particular interest if there exists a maximal $(\iota, K, \mu)$.

Conversely, assume that a $C^*$-algebra $\A$ of bounded functions $f: X
\to \C$ and a mean $m$ on $\A$ is given. Let $(\iota_{\A},K_{\A})$ be the natural compactification for $\A$. By Proposition \ref{Palgebracontinuous} the mapping $\iota^*: F \to F \circ \iota$ is a bijection between $C(K_{\A})$ and $\A$. Thus $m'(F):=m(F \circ \iota)$ is well defined and a mean on $C(K_{\A})$. By Riesz' Representation Theorem \ref{riesz} $m'$ induces a Borel probability measure $\mu$ on $K_{\A}$
with $m'(F) = \int F d\mu$ for all $F\in C(K_{\A})$ which is unique on the
$\sigma$-algebra of Borel sets and its $\mu$-completion. So it is
not surprising that the $m$-closure of $\A$ contains all $f = F \circ \iota$ with $F \in \Riemanntop_{\mu}$.

\begin{pro}\label{Pmclosed}
Let $\A$ be a $C^*$-algebra on $X$ and $m$ a mean on $\A$. Let $(\iota_{\A},K_{\A},\mu)$ be the compactification
where $(\iota_\A,K_{\A})$ is the natural compactification for $\A$ and
$\mu$ is the complete and regular Borel measure on $K_{\A}$ which
satisfies $\int F\, d\mu = m(F \circ \iota_{\A})$ for all $F \in C(K_{\A})$.
Then $$\Riemanntop_{\mu}^*:=\{F\circ\iota_{\A}: F\in \Riemanntop_{\mu}\} \subseteq \overline{\A}^{(m)}$$ for the $m$-completion $\overline{\A}^{(m)}$ of $\A$. Furthermore, if $\A$ separates points of $X$, then $\Riemanntop_{\mu}^*=\overline{\A}^{(m)}$.
\end{pro}
\begin{proof}
Let $F\in \Riemanntop_{\mu}$, then it is straight-forward to check that
$F\circ\iota_{\A}$ is in the $m$-closure of $\A$, i.e.\ $\Riemanntop_{\mu}^*\subseteq \overline{\A}^{(m)}$.
Assume now that $\A$ separates points and take the real-valued function $f\in \overline{\A}^{(m)}$. By definition,
for every $\eps>0$ there are real-valued $F_1,F_2 \in C(K_{\A})$ such that $F_1\circ\iota_{\A}\le f \le F_2\circ\iota_{\A}$ with $\int (F_2-F_1)d\mu \le \eps$. Observe
\begin{equation*}\label{Eqsupinf}
 F_{\flat}:=\sup_{F_1\circ\iota_{\A}\le f} F_1 \le F\le  \inf_{f\le F_2\circ\iota_{\A}} F_2=:F^{\#}, \quad F_1,F_2 \in C(K_{\A}).
\end{equation*}
The fact that $f$ is in the $m$-closure of $\A$ implies
that every $F$ with $F_{\flat}\le F\le F^{\#}$ is $\mu$-Riemann integrable.
Since $\A$ separates points of $X$ the map $\iota_{\A}: X \to K_{\A}$ is one-one, cf.\ Proposition \ref{Palgsepa}.
Thus we can define a function $F^\circ: K_{\A} \to \R$ via
\[F^\circ(k)= \left\{ { f(x) \atop 0} { \mbox{ if } k = \iota_{\A} (x) \atop \mbox{ otherwise }}\right. .\] Then $F:=\max\{F^{\circ}, F^{\flat}\}$ is $\mu$-Riemann integrable and satisfies $F\circ\iota_{\A}=f$. Hence $\Riemanntop_{\mu}^*\supseteq \overline{\A}^{(m)}$.
\end{proof}

\begin{rem}\label{Rpointseparating}
In any case $\Riemanntop_{\mu}^*$ is a $C^*$-algebra. $\iota_{\A}^*: F\mapsto F\circ\iota_{\A}$ is a bounded $*$-homomorphism which maps $\Riemanntop_{\mu}$ into $B(X)$ and the image of every bounded $*$-homo\-morphism is again a $C^*$-algebra, cf.\ \cite[Theorem I.5.5]{Davi96}.
\end{rem}

In the general case of a non point-separating algebra $\A\subseteq B(X)$ we can do a general construction: Consider the equivalence relation on $X$ defined by $x\approx y$ if for every $f \in \A$ we have $f(x)=f(y)$. Then $\A$ induces an algebra $\A/_{\!\approx} \subseteq B(X/_{\!\approx})$ which is isomorphic to $\A$ but point separating.

\begin{exa}\label{E-closure}
 Let $X=\{a,b,c\}$, $\A:=\{f: X \to \C: f(a)=f(b)\}$ and consider the fapm $\delta_{c}$. Then $K_{\A}:=\{\alpha,\beta\}$ is a two element set and $\iota_{\A}(a)=\iota_{\A}(b)=:\alpha$. The inclusion $\A=\Riemanntop_{\delta_{\beta}}^*\subset\overline{\A}^{m_{\delta_c}}=B(X)$ is strict, showing that in the last statement of Proposition \ref{Pmclosed} the point-separation property can not be omitted.
\end{exa}

However, the identity $\Riemanntop_{\mu}^*= \overline{\A}^{(m)}$ may hold even for certain non point separating algebras $\A$, e.g.\ if $\A$ consists of constant functions.

\section{The set of all means}\label{Smeans}

For an infinite discrete set $X$ there is an abundance of means on
the algebra $B(X)$ of bounded $f:X \to \C$. For a better understanding
of the structure of the set of all means, compactifications turn out
very useful. We start by restricting to very special means, namely
multiplicative ones. (A mean $m$ defined on an algebra $\A$ of
functions is called multiplicative if $m(f_1\cdot f_2)=m(f_1)m(f_2)$ for
all $f_1,f_2 \in \A$.) As a standard reference we mention \cite{Gran72a}.

Given a multiplicative mean $m$ on $B(X)$,
let $p=p_m$ be the corresponding fapm, defined on the whole power set
$\Algebra= \P(X)$ of $X$.
For every $A \tm X$ multiplicativity of $m$ yields
$p_m(A) = m(\Eins_A) = m(\Eins_A \cdot \Eins_A)
= m(\Eins_A)m(\Eins_A) = p_m^2(A)$,
hence $p_m(A) \in \{0,1\}$.

Conversely, every fapm $p$
defined for all $A \tm X$ and taking only
the values $0$ and $1$ induces a multiplicative mean on $B(X)$:
First check that in all four possible cases for $p(A_1),p(A_2) \in \{0,1\}$
one gets $p(A_1 \cap A_2) = p(A_1)p(A_2)$. This implies
$m_p(f_1\cdot f_2) = m_p(f_1) m_p(f_2)$ whenever $f_i = \Eins_{A_i}$.
By multiplicativity and distributivity this transfers to
$f_i \in \Simple_{\Algebra}$. Finally observe that $B(X)$ is the uniform
closure of $\Simple_{\Algebra}$ to conclude by standard approximation
arguments that $m_p$ is indeed a multiplicative mean on $B(X)$.

Thus multiplicative means are in a one-one correspondence with
fapm's on the power set taking only the values 0 and 1. For an
arbitrary such $p$ the system $\F_p$ of all $A \tm X$ with $p(A)=1$
is closed under finite intersections, supersets and contains $X$, i.e.\
$\F_p$ is a filter. Since for each $A$ either $p(A)=1$ or $p(X
\setminus A)=1$, $\F_p$ is an ultrafilter. Obviously also this
argument is reversible: Every ultrafilter $\F$ on $X$ induces a fapm
$p_{\F}$ by $p_{\F}(A) = \Eins_{\F}(A)$ which takes only the values
0 and 1. Consider the Stone-{\v C}ech compactification $\beta X$ as
the space of ultrafilters on $X$, according to Proposition
\ref{Pstonecechsummary}. Then the means $m$ on $B(X)$
transfer to positive linear functionals on $C(\beta X)$ and thus, by
Riesz' Representation Theorem, to Borel probability measures on $\beta X$. The
functionals, taking only the values 0 and 1, are point evaluations
$F \mapsto F(y)$, $F \in C(\beta X)$, corresponding to Dirac
measures $\delta_y$ concentrated in the point $y \in \beta X$. As an
ultrafilter, $y$ contains exactly those $A \tm X$ with $p(A)=1$.
Note that, in the set of all sub-probability Borel measures,
normalized point measures are exactly the extreme ones, i.e.\ they
can be represented as a convex combination only in the trivial way.
In the weak-*-topology the set of all sub-probability measures is compact.
Thus, by the Krein-Milman Theorem (cf.\ for instance \cite{Rudi91}),
an arbitrary Borel measure on $\beta X$ is in the weak-*-closure of
the convex hull of certain point measures. Going back to $X$ and
means on $X$ we thus have:

\begin{pro}\label{Pstonecechmm}
The set of all means on $B(X)$, $X$ discrete, is given by the
convex hull of all multiplicative means on $X$. The multiplicative
means on $X$ are in a natural bijective correspondence with the
points of the Stone-{\v C}ech compactification $\beta X$.
\end{pro}

Indlekofer has systematically used the relation between
means and fapm's on $\N$ or $\Z$ with probability measures on
the Stone-{\v C}ech compactification in probabilistic number theory
(cf.\ for instance \cite{Indl02}).\\

In Section \ref{Scompcont} we have seen that for $f:X \to \C$, $X$ discrete, there is a smallest
continuous representation which we called the natural one and which is unique up to equivalence.
Looking for Riemann integrable representations, also a measure has
to be involved and thus the situation is more complicated. This has
the consequence that there is not one unique (up to equivalence)
smallest Riemann integrable representation. Nevertheless, at least
for discrete $X$, one can easily find many minimal compactifications:

\begin{exa}\label{Ecompacttopology}
Let $X$ be discrete. For given bounded $f: X \to \C$ equip $K:=f(X)$ with a compact topology and let $\iota: x \mapsto f(x)$. Then $(\iota,K)$ is a compactification and $F: K \to \C$, $k \mapsto k$
is the only representation of $f$ in $(\iota,K)$.
This representation
clearly is minimal, provided $K$ carries an appropriate Borel probability measure $\mu$. If $K$ is finite, the discrete topology is compact and does the job as well as any probability measure $\mu$ defined on $\P(K)$.
In the infinite case we define a compact topology on $K$
by fixing any $k_0 \in K$ and taking as open sets all subsets of $K$
not containing $k_0$ and all cofinite sets which contain $k_0$.
Note that all $k \in K \setminus \{k_0\}$ are isolated points,
hence every function is continuous in such $k$. The only possible
discontinuity point is $k_0$. Thus, provided $\mu(\{k_0\})=0$, we have $\Riemanntop_{\mu}=B(K)$. In particular $F$ is $\mu$-Riemann integrable.
\end{exa}

For many reasons this construction is not very satisfactory.
One of them is that there is no
canonical choice of $\mu$.
The most natural way to find canonical measures is by invariance
requirements. In the forthcoming chapters
we will be concerned with invariance mainly with
respect to group or semigroup operations, to some extent also with
respect to a single transformation in the sense of topological
and symbolic dynamics.

\chapter{Invariance under transformations and operations}\label{CHinvariance}
\section{Invariant means for a single transformation}\label{Ssingle}

At the end of the previous chapter we have seen that there is an
abundance of means on an infinite discrete set $X$. If $X$ carries
further structure one asks for means and measures with
certain interesting additional, namely invariance properties.

\begin{defi}\label{Dinvariant}
Let $X$ be any nonempty set and $T: X \to X$. Then
$U_T: B(X) \to B(X)$ is defined by $f \mapsto f \circ T$. A set $\A \tm B(X)$ is
called \textbf{$\mathbf{T}$-invariant} if $U_T(\A) \tm \A$. Assume that $\A$ is a
\textbf{$\mathbf{T}$-invariant} vector space and $m$ is a mean on $\A$. Then $m$ is
called \textbf{$\mathbf{T}$-invariant} if $U_T^*(m)= \circ U_T=m$, i.e.\ if
\[ m(f \circ T) = m(f) \]
for all $f \in \A$. By $M(\A)$ we denote the set of all means on $\A$,
$M(X):=M(B(X))$, and by $M(\A,T)$ the set of all $T$-invariant $m \in M(\A)$.
For bijective $T$ we call $\A$ resp.\ $m$ \textbf{two-sided $\mathbf{T}$-invariant} if
it is both, $T$- and $T^{-1}$-invariant.
\end{defi}

Note that in the two sided invariant case one has $M(\A,T)=M(\A,T^{-1})$.
It is easy to check that $M(\A)$ and $M(\A,T)$ are
weak-*-closed subsets of the unit ball in $B(X)^*$, the
dual space of the Banach space $B(X)$. Thus, since by the Banach-Alaoglu
Theorem the dual unit ball is compact in this topology,
$M(\A)$ and $M(\A,T)$ are compact as well.
More directly, compactness becomes clear from applying Tychonoff's Theorem to
\[
M(\A) \tm \prod_{f \in \A} \{z \in \C: |z|\le ||f||_{\infty}\}.
\]
As a consequence, any sequence $m_n \in M(\A)$ has at least
one accumulation point (accumulation measure) $m \in M(\A)$.
In particular the set $M_{T,(m_n)}$ of accumulation means of the
sequence $(\overline m_n)$ is nonempty:
\[
\overline m_n:= \frac 1n \sum_{k=0}^{n-1} m_n(f \circ T^k) =
   \frac 1n \sum_{k=0}^{n-1} m_n (U_T^k(f)) =
   \frac 1n \sum_{k=0}^{n-1} {U_T^*}^k(m_n)(f).
\]

\begin{pro}\label{Paim}
$M_{T,(m_n)} \tm M(\A,T)$. In particular there are $T$-invariant
means. We can take for instance the point evaluation
$m_n := m_{\delta_x}: f \mapsto f(x)$ for any $x \in X$.
\end{pro}
\begin{proof}
Let $m \in M_{T,(m_n)}$. As $m$ is an accumulation mean, for every $\eps>0$ and bounded $f: X\to \C$ there is a sequence $n_1<n_2<\ldots$ such that both
$|m(f)-\overline{m}_{n_k}(f)|\le \eps$ and
$|m(f\circ T)-\overline{m}_{n_k}(f\circ T)|\le \eps$ for all $k\in \N$.
From the defining properties of $\overline m_{n_k}$ we obtain
\begin{align*}
\left|\overline {m}_{n_k} (f \circ T) - \overline {m}_{n_k} (f)\right|
=& \frac 1{n_k} \left|\sum_{j=0}^{n_k-1} m_{n_k}(f \circ T^{j+1} - f \circ T^j)\right|\\
=& \frac 1{n_k} \left|m_{n_k}(f \circ T^{n_k}) - m_{n_k}(f)\right|\le \frac 2{n_k}||f||_{\infty},\\&\\
|m(f)-m(f\circ T)|\le&\hspace{3.2ex} |m(f)-\overline{m}_{n_k}(f)| + \left|\overline {m}_{n_k} (f)-\overline {m}_{n_k} (f \circ T)\right|&\\&+\,|\overline{m}_{n_k}(f\circ T)-m(f\circ T)|&\\
\le&\, 2\eps + \frac 2{n_k}||f||_{\infty}.
\end{align*}
As this is true for all $k\in N$ and $\eps>0$ we obtain $T$-invariance of $m$.
\end{proof}

We study now which values of $m(f)$ are possible for
$m \in M(T,\A)$ and $f \in \A$.

\begin{pro}\label{Pvaluesofmeans}
Let $T: X \to X$, $\A \tm B(X)$ a $T$-invariant vector space,
$f \in \A_{\R}$. Then the set $\{m(f):\ m \in M(\A,T)\}$ coincides
with the interval $[a,b]$ where, with the short hand
\begin{equation}\label{EqArithmeticmean}
  a = \lim_{n \to \infty} \inf_{x \in X} s_n(x) \quad \mbox{and} \quad
  b = \lim_{n \to \infty} \sup_{x \in X} s_n(x), \quad
  s_n = s_{n,T,f}:= \frac 1n \sum_{k=0}^{n-1} f \circ T^k.
\end{equation}
In particular this set does not depend on $\A$.
\end{pro}
\begin{proof}
Note first that $m(s_n)=m(f)$ for every $m \in M(\A,T)$. For the
proof it suffices to show that, for every $\alpha \in \R$, there is
an $m \in M(\A,T)$ with $m(f) = \alpha$ if and only if the following
condition is satisfied:

Condition(C): \emph{For all $\e>0$ and $n \in \N$ there are $x=x(\e,n), y=y(\e,n) \in X$ such that $s_n(x) > \alpha - \e$ and $s_n(y) < \alpha + \e$.}

Necessity of (C): Let $m(f) = \alpha$ with $m \in M(\A,T)$ and suppose, by
contradiction, that (C) fails. Then there is an $\e >0$
and an $n \in \N$ such that $s_n(x) \le \alpha - \e$ for all $x \in X$
(the case $s_n(x) \ge \alpha + \e$ can be treated similarly), hence
$m(f) = m(s_n) \le ||s_n||_{\infty} \le \alpha - \e$, contradiction.

Sufficiency of (C):
Assume that (C) holds and consider the point measures
$m_n:= \delta_{x(1/n,n)}$. With the notation of Proposition
\ref{Paim} this means $\overline m_n(f) > \alpha - \tfrac 1n$ for all $n$.
Use Proposition \ref{Paim} to find an $m' \in M_{T,(m_n)} \tm M(\A,T)$.
Then $m'(f) \ge \alpha$. Similarly one finds an $m'' \in M(\A,T)$
with $m''(f) \le \alpha$. It follows that there is a $\lambda \in [0,1]$
such that $$\lambda m'(f) + (1-\lambda)m''(f) = \alpha.$$ Since $M(\A,T)$ is convex $m:= \lambda m' + (1-\lambda)m''$
has the required properties.
\end{proof}


Of particular interest are the functions $f$ with a unique mean value.
\begin{defi}\label{Dalmostconv}
Let $T: X \to X$ and $\A \tm B(X)$ be a $T$-invariant linear space.
A function $f \in \A$ is called \textbf{$\A$-almost convergent} if $m(f)$ has
the same value for all $m \in M(\A,T)$. The set of all $\A$-almost convergent
$f \in \A$ is denoted by $AC(\A,T)$, for $\A=B(X)$ we also write
$AC(B(X),T) = AC(T)$. We write $m_{\A}$ for the restriction of
$m \in M(\A,T)$ to $AC(\A,T)$. If $AC(\A,T) = \A$
we call $T$ \textbf{uniquely ergodic} (with respect to $\A$).
\end{defi}

By definition, $m_{\A}$ does not depend on $m$.
It is clear that $AC(\A,T)$ is a $T$-invariant uniformly closed
linear space containing all constant functions.
Furthermore $AC(\A,T) = AC(T) \cap \A$.
Finally, $f \in AC(\A,T)$ with $m_{\A}(f) = \alpha$
if and only if $f \in \A$ and, for all $x_n \in X$,
\[
\lim_{n \to \infty} s_{n,T,f}(x_n) = \frac 1n \sum_{k=0}^{n-1} f(T^k(x_n))
   = \alpha.
\]
The obvious way to define $T$-invariance of a set algebra $\Algebra$
on $X$ or a finitely additive measure $p$ defined on $\Algebra$ is
to require that $\{\Eins_A:\ A \in \Algebra\} \tm B(X)$ resp.\ $m_p$
as defined in Section \ref{Sfam} is $T$-invariant. Since $\Eins_A
\circ T = \Eins_{T^{-1}[A]}$ this is the case if and only if
$T^{-1}[A] \in \Algebra$ resp.\ $p(T^{-1}[A])=p(A)$ for all $A \in
\Algebra$. From Proposition \ref{Pvaluesofmeans} we get:

\begin{pro}\label{Pfapmvalues}
The possible values $p(A)$ for $T$-invariant fapm $p$ are given by the interval $[a,b]$ where
\[
a = \lim_{n \to \infty} \inf_{x \in X} d_n(x) \quad\mbox{ and }\quad
b = \lim_{n \to \infty} \sup_{x \in X} d_n(x), \quad
d_n(x) = \frac 1n \left|\{k:\ T^k(x) \in A,\, 0\le k< n\}\right|.
\]
In particular, $p(A)$ takes the same value for all $T$-invariant $p$ if and only if $a=b$.
\end{pro}

\section{Applications}\label{Sapplications-means}

\subsection{Finite $X$}\label{SSfinite}
Let $X$ be finite, $T:X \to X$ and
$\A=B(X)=\C^{X}$. For every $x \in X$ there is a minimal $m \ge 0$
and a minimal $k > n$ such that $T^k(x)=T^m(x)$. We call $C_x:=
\{T^m(x), T^{m+1}(x), \ldots,T^{k-1}(x)\}$ the \textbf{cycle (cyclic attractor)}  induced by $x$ and $B_x = B(C_x) := \{y \in X:\ C_y=C_x\}$ the \textbf{basin} of $C_x$.
It is clear that $C_x \tm B_x$, the $B_x$ forming a partition.
$C_x=B_x$ if and only if the restriction of $T$ to  this set is
bijective. Furthermore the $s_n = s_{n,T,f}$ defined by
\[s_n(x) = \frac 1n \sum_1^{n-1} f(T^k(x))\]
converge to a function $\overline f$ which, on each $C=C_x$, takes
the constant value
\[m_C(f):= \frac 1{|C|} \sum_{y \in C} f(y).\]
It is clear that $m_C \in M(B(X),T)$ for each cycle $C$. The same holds
for all convex combinations. We claim that, conversely, every
$m \in M(B(X),T)$ is of this type, i.e.\ $m= \sum_C \lambda_C m_C$ with $0 \le \lambda_C \le$ for all $C$ and $\sum_C \lambda_C = 1$.
To see this, define $\lambda_C:= m(\Eins_{B(C)})$ and observe that
\[\overline f = \sum_C m_C(f) \Eins_{B(C)}.\] This implies
\[m(f) = m(s_n) = m(\overline f) = \sum_C m_C(f) m(\Eins_{B(C)}) =
   \sum_C \lambda_C m_C(f).\]
The uniqueness of the $\lambda_C$ follows since the $m_C$ are linearly independent. This gives an obvious description of almost convergent functions: $f \in AC(B(X),T)$ if and only if $m_C(f)$ takes the same value for all cycles $C$.

In terms of measures this means that every $T$-invariant $p$ is a convex combination of the ergodic measures $p_C$ defined by $p_C(A):= \tfrac{|A \cap C|}{|C|}$. This is the finite version of the ergodic decomposition given by Birkhoff's Ergodic Theorem (cf.\ for instance \cite{Walt82}). Infinite $X$ would have to be treated in this context, but we do not go further into this direction.

\subsection{$X = \Z, T: x \mapsto x+1$}
\label{SSintegersshift}
First note that whenever $\A\subseteq B(\Z)$
is two sided $T$-invariant then $T$-invariance of a mean $m$ or a
fam $p$ implies invariance with respect to all translations on the
additive group $\Z$. In Section \ref{Sinvgroups} we will focus on
this aspect. Here we want to apply our analysis from Section
\ref{Ssingle}. For this reason we have to consider the quantities
\[s_{N,f}(n):= \frac 1N \sum_{k=n}^{n+N-1} f(k)\]
and, for real valued $f$, the corresponding lower and upper limits
\[m_*(f) := \lim_{N \to \infty} \inf_{n \in \Z} s_{N,f}(n)
\quad \mbox{and} \quad
m^*(f) := \lim_{N \to \infty} \sup_{n \in \Z} s_{N,f}(n).\]
For $f= \Eins_A$ these values $m_*(A):= m_*(\Eins_A)$ and
$m^*(A):= m^*(\Eins_A)$ are known as \textbf{lower Banach density} resp.\  \textbf{upper Banach density} of $A$. The possible values of $T$-invariant measures are given by \[\{m(f):\ m \in M(T,\A)\} = [m_*(f),m^*(f)],\]
hence \[AC(\A,T) = \{f \in \A:\ m_*(f)=m^*(f)\}.\]
The restriction of $m_*$ and $m^*$ to $AC(\Z):= AC(B(X),T)$ is known as \textbf{Banach density} and denoted by $m_B$.

We have already mentioned that the set $AC(\Z)$ of almost convergent
$f$ on $\Z$ is a linear space and, furthermore, uniformly closed.
Having the results about compactifications and complex valued
functions in mind, we ask whether $AC(\Z)$ is an algebra as well.
But this is not the case as the following example shows.
\begin{exa}
Consider the sets $A=2\Z$ of even numbers and $B = (B_1 \cup
B_2)\cup(-B_1 \cup -B_2)$ with
\begin{eqnarray*}
B_1 &=& \bigcup_{n=1}^{\infty} ((2n-1)!,(2n)!] \cap 2\Z,\\
B_2 &=& \bigcup_{n=1}^{\infty} ((2n)!,(2n+1)!] \cap (2\Z+1).
\end{eqnarray*}
In $B$ one has very long blocks of even numbers alternating with
very long blocks of odd numbers. It is clear that both $A$ and $B$
have Banach density $\tfrac 12$ while $B_1 = A \cap B$ has
lower Banach density $0$ not coinciding with its upper Banach density $\tfrac 12$.
It follows that $\Eins_{B_1} = \Eins_A\cdot\Eins_B \notin AC(\Z)$
although $\Eins_A, \Eins_B \in AC(\Z)$. Thus $AC(\Z)$ is
not an algebra. In particular there is no compactification $(\iota,K)$ of $\Z$ such that $AC(\Z)$ is the set of all $f$ having a continuous (or Riemann integrable) representation in $(\iota,K)$.\end{exa}

\subsection{$X$ compact, $\A = C(X)$, $T$ continuous}\label{SScompact}
Let $X$ be a compact space, $\A=C(X)$ the algebra of complex valued
continuous functions on $X$ and $T: X \to X$ continuous.
Then $C(X)$ is $T$-invariant since $f \in C(X)$ implies
$T^*f = f \circ T \in C(X)$. (If $T$ is bijective then it is a homeomorphism, hence $C(X)$ is even both sided $T$-invariant.) This is the framework of classical topological dynamics.

Proposition \ref{Paim} guarantees that there are $T$-invariant means
on $C(X)$. By Riesz' Representation Theorem \ref{riesz} every $m \in M(T,C(X))$ induces a unique regular Borel probability measure $\mu_m$ with
\[ m(f) = \int_X f\, d\mu_m \quad \mbox{for all $f \in C(X)$}.\]

\begin{defi}\label{Dinvmeanmeasure}
Let $\mu$ be a measure defined on a $\sigma$-algebra $\Algebra$ on $X$ and $T: X \to X$ measurable. Then $\mu$ is called $T$-invariant if $\mu(T^{-1}[A]) = \mu(A)$ for all $A \in \Algebra$.
\end{defi}

In the context of compact $X$ we are particularly interested in the
case that $\Algebra$ contains all Borel sets and that $\mu$ is
complete. In the case of regular Borel measures the invariance of a mean $m$ is equivalent with invariance of the corresponding $\mu_m$:

\begin{pro}\label{Pinvmeanmeasure}
Let $X$ be a compact space, $T: X \to X$ continuous, $m$ a
mean on $C(X)$ and $\mu$ a Borel measure, i.e.\ defined on a $\sigma$-algebra $\Algebra$ containing all Borel sets, such that $m(f) = \int_X f\, d\mu$ for all $f \in C(X)$. Then:
\begin{enumerate}
\item If $\mu$ is $T$-invariant then $m$ is $T$-invariant.
\item Assume that $\mu$ is regular. If $m$ is $T$-invariant then $\mu$ is $T$-invariant.
\item The implication in the second statement does not hold if one
  drops the regularity assumption on $\mu$.
\end{enumerate}
\end{pro}
\begin{proof}[Sketch]
\begin{enumerate}
\item In order to show that $m(f \circ T)=m(f)$ for all $f \in C(X)$  one first considers $f = \Eins_A$
  with $A \in \A$, then linear combinations of such $f$ and then uses
  that any $f \in C(X)$ can be uniformly approximated by such linear
  combinations.
\item Let $A \in \Algebra$ and $\e>0$. By regularity of $\mu$ there
  are a closed set $C$ and an open set $O$ such that $C \tm A \tm O \tm X$
  and $\mu(O \setminus C) < \e$, and corresponding Urysohn functions, i.e.\
  continuous $f,g: X \to [0,1]$ with
  $\Eins_C \le f \le \Eins_A \le g \le \Eins_O$. By invariance of $m$
  we obtain
  \[
  \mu(T^{-1}[A]) \le m(g \circ T) = m(g) \le \mu(O) \le \mu(A) + \e
  \]
  and similarly $\mu(T^{-1}[A]) \ge \mu(A) - \e$, hence
  $\mu(T^{-1}[A]) = \mu(A)$.
\item We use the example $X = [0,\omega_1]$ from the end of
  Section \ref{Scompintegral} and the constant mapping
  $T: X \to X$, $x \mapsto \omega_1$. Then the point evaluation mean
  $m: f \mapsto f(\omega_1)$ is $T$-invariant (as well as the corresponding
  point measure $\delta_{\omega_1}$ concentrated in $\omega_1$).
  The measure $\mu$ from the end of Section \ref{Scompintegral}
  satisfies $m(f) = \int_X f(x)\, d\mu$ for all $f \in C(X)$. Nevertheless
  for $A=\{\omega_1\}$ we have $\mu(A) = 0 \neq 1 = \mu(X) = \mu(T^{-1}[A])$. \qedhere
\end{enumerate}
\end{proof}

If $M(C(X),T)$ consists of only one measure, we say that $T$ is \textbf{uniquely ergodic}.
\begin{cor} Let $M(C(X),T)=\{m\}$.
For every $f\in C(X)$ \[s_{n,T,f}(x) := \frac 1n \sum_{k=0}^{n-1} f( T^k(x)) \to m(f), \quad \mbox{ uniformly in } x\in X,\]
in particular, the uniform closure of the convex hull of the $T$-orbit contains the constant function $m(f)\cdot\Eins_X$.
\end{cor}
\begin{proof}
By Proposition \ref{Pvaluesofmeans} we know that
\[\lim_{n\to\infty} \left(\inf_{x\in X} s_{n,T,f}(x)\right)=m(f)
=\lim_{n\to\infty} \left(\sup_{x\in X} s_{n,T,f}(x)\right),\]
i.e.\ for every $\eps>0$ and large enough $n$ we have $\sup_{x\in X} |s_{n,T,f}(x)-m(f)|<\eps$.
\end{proof}

If $T$ is uniquely ergodic, $m(f)\Eins_X$ is the only constant function in the uniform closure of the convex hull of the $T$-orbit. For a arbitrary constant $\lambda\in \C$ in this closure we have
\[\lambda = m(\lambda\Eins_X)=\lim_{k\to \infty} m(f\circ T^{n_k})=m(f),\quad \mbox{whenever } f\circ T^{n_k}\to \lambda\Eins_X.\]

\subsection{Shift spaces and symbolic dynamics}
\label{SSshiftspaces}
We now consider a special case of the situation treated in
Section \ref{SScompact}. Let $A$ be a finite set. In this context $A$ is called
an alphabet and its members are called symbols. Furthermore let $X$
be a closed subset of the compact space $A^{\Z}$ which is
shift invariant, i.e.\ $\sigma(X) = X$ for the shift
$\sigma: (a_k)_{k \in \Z} \mapsto (a_{k+1})_{k \in \Z}$.
Such dynamical systems $(X,\sigma)$ are also called subshifts and
are the main objects of symbolic dynamics.

The importance of this apparently simple objects is due to the
abundance of $\sigma$-invariant closed subsets $X$ of $A^{\Z}$
by which a quite big class of dynamical systems can be represented
in a reasonable way. Assume that $Y$ is a compact space,
$T: Y \to Y$ a continuous transformation,
$Y = Y_0 \cup Y_1 \cup \ldots \cup Y_{s-1}$ (pairwise disjoint union),
$s \in \N$, a finite partition of $Y$ and $F: Y \to \{0,1,\ldots,s-1\}$,
$y \mapsto i$ if $y \in Y_i$.

Note that $F$ is in general
(in particular for connected $Y$) not continuous but,
if the $Y_i$ are continuity sets w.r.t.\ an appropriate measure on $Y$,
Riemann integrable. Therefore we must not expect that the induced
mapping $y \mapsto x(y):= (F(T^ky))_{k \in \Z}$ is continuous.
But we may consider the closure $X$ of its image which is shift invariant.
In many interesting cases of a metrizable $Y$ there is a continuous
surjection $\varphi: X \to Y$ such that $y = \varphi(x)$
whenever $x=x(y)$. One says that $(Y,T)$ is a \textbf{factor} of $(X,\sigma)$ and we have the following commutative diagram
\begin{diagram}
X            &\rTo^{\phi}&Y\\
\dTo>{\sigma}&         &\dTo>{T} \\
X            &\rTo^{\phi} &Y.
\end{diagram}

\begin{exa}[Sturmian sequence]\label{Esturmsequence}
Let $Y=\R/\Z$, $\alpha \in Y$ and $T$ be the
homeomorphism $y \mapsto y + \alpha$ (rotation). Consider the partition of $Y$ into two segments  $Y_0 = [0,\alpha) + \Z$ and
$Y_1 = [\alpha,1) + \Z$. Hence $A=\{0,1\}$.
For the $T$-orbit $x(0)$ of $y=0$, i.e.\
$x(0)=(a_k)_{k \in \Z}$ with $a_k=1$ if $k \alpha \in Y_1$ and $a_k=0$
if $k \alpha \in Y_0$ it turns out that $X$ is the closure of the
shift orbit $\{\sigma^k(x(0)):\ k \in \Z\}$ in $\{0,1\}^{\Z}$ of $x(0)$
in $\{0,1\}^{\Z}$. The two sided sequence $x(0)$ is an example of a
so called Sturmian sequence, a class providing some of the most simple
but typical examples of Hartman functions, the main topic of
the forthcoming chapters.
\end{exa}

In order to understand invariant means or, equivalently by Riesz' Representation Theorem \ref{riesz} and Proposition \ref{Pinvmeanmeasure}, invariant measures
on shift spaces $X \tm A^{\Z}$, note that every invariant
Borel probability measure on $X$ is uniquely determined by its values
on cylinder sets $S=[b_0,\ldots,b_{n-1}] = \{(a_k)_{k \in \Z}:\
a_k=b_k,\ k=0,\ldots,n-1\}$. Thus it suffices to consider functions
$f=\Eins_S$ for such $S$ and apply the arguments from Section \ref{Ssingle}.
For instance the numbers $a,b$ in Proposition \ref{Pvaluesofmeans}
can be described in terms of relative frequencies of blocks
$(b_0,\ldots,b_{n-1})$ in symbolic sequences $x \in X$.

\subsection{The free group $F(x,y)$}
\label{SSfreegroup}
Let $X=F(x,y)$ denote the group with two free generators $x$ and $y$.
As usual we assume
each member $w \in F(x,y)$ to be a reduced word built up from the four allowed
symbols $x,y,x^{-1},y^{-1}$, including the empty word $w=\emptyset$.
Denote by $W_x$ all reduced words ending with the symbol $x$; $W_{x^{-1}},
W_{y}$ and $W_{y^{-1}}$ are defined analogously.
Consider the transformation $T_1: w \mapsto wx^{-1}$. In particular we have $T_1(W_x) = W_x \cup W_y \cup W_{y^{-1}} \cup
\{\emptyset\}$.
Assume $m_1 \in M(T_1,B(X))$. Then, for the associated fapm $p_1=m_{p_1}$ we compute
\begin{equation}0 = p_1(T_1(W_x))-p_1(W_x) = p_1(W_y \cup W_{y^{-1}} \cup \{\emptyset\}).\end{equation}
By symmetry,
for $T_2: w \mapsto wy^{-1}$ we have $T_2(W_y) = W_y \cup W_x \cup W_{x^{-1}} \cup \{\emptyset\}$ and every $m_2 \in M(T_2,B(X))$
satisfies
\begin{equation}0 = p_2(T_2(W_y))-p_2(W_y) = p_2(W_x \cup W_{x^{-1}} \cup \{\emptyset\}).\end{equation}
In particular we have $T_1(W_x)\setminus W_x \cup T_2(W_y)\setminus W_y =X$.
Assume $m \in M(T_1,B(X)) \cap M(T_2,B(X))$. The associated fapm $p=m_p$ satisfies
\begin{eqnarray}
1 = p(X)&=&p\left(T_1(W_x)\setminus W_x \cup T_2(W_y)\setminus W_y\right)\\
 &\le& p(T_1(W_x))-p(W_x)\, +\,  p(T_2(W_y))-p(W_y) = 0,
 \end{eqnarray}
contradiction. Thus there is no mean on $B(X)$ which is both $T_1$- and $T_2$-invariant.

Since $T_1$ and $T_2$ are group translations this shows that
$X=F(x,y)$ is not an amenable group, see also Section \ref{SSamenable}. Together with the fact that $F(x,y)$ can be realized as a group of orthogonal transformations of $\R^3$ this is the core of the celebrated Banach-Tarski paradox, see \cite{Wago85}. We will focus on (semi-)group structures systematically in Section \ref{Sinvgroups}.

\section{Compactifications for transformations and actions}\label{Strafocomp}

We have seen the importance of compactifications for means and measures
already in Chapter \ref{CHprel}. The r\^ole
of transformations for identifying interesting measures in terms
of invariance properties was pointed out in Section \ref{Ssingle}. We now combine both points of view
by investigating the following setting.

\begin{defi}\label{Dcompcompatilble}
Let $X$ be a topological space,
$T: X \to X$ continuous and $(\iota,K)$ a compactification of $X$. Then $(\iota,K)$ is called
a compactification \textbf{compatible} with $T$ and, vice versa,
$T$ a transformation \textbf{compatible} with $(\iota,K)$ if there is a
continuous $T_K: K \to K$ such that $\iota \circ T = T_K \circ \iota$,
i.e.\ making the diagram
\begin{diagram}
X            &\rTo^{\iota} &K\\
\dTo>{T}&         &\dTo>{T_K} \\
X            &\rTo^{\iota} &K.
\end{diagram}
commutative. $T_K$ is called a \textbf{continuous extension} of $T$
in $(\iota,K)$.
\end{defi}

Continuous extensions are unique and compatible with composition:

\begin{pro}\label{Pcontextension}
Let $T_K$ and $T'_K$ be continuous extensions of $T:X \to X$ in the
compactification $(\iota,K)$ of $X$. Then $T_K=T'_K$. Furthermore,
if $S_K$ is a continuous extension of $S: X \to X$ in $(\iota,K)$
then $(S \circ T)_K = S_K \circ T_K: K \to K$, is the continuous
extension of $S \circ T: X \to X$.
\end{pro}
\begin{proof}
The first statement follows since $T_K$ and $T'_K$ are continuous
and coincide on the dense set $\iota(X) \tm K$. For the second statement
observe
\begin{eqnarray*}
(S_K \circ T_K) \circ \iota = S_K \circ (T_K \circ \iota) &=&
S_K \circ (\iota \circ T)\\
& =& (S_K \circ \iota) \circ T =
(\iota \circ S) \circ T = \iota \circ (S \circ T)
\end{eqnarray*}
and use the uniqueness of continuous extensions to obtain
$S_K \circ T_K = (S \circ T)_K$.
\end{proof}

Note that from a certain point of view
for continuous extensions of transformations $T: X \to X$ the situation
is more complicated than for representations of complex valued
$f: X \to \C$ in the sense of Section \ref{Sopcomp}.
This is due to the fact that there is no obvious analogue of
Proposition \ref{Palgebracontinuous} which yields a natural
compactification for $f$ or even for a unital $C^*$-algebra $\A$
which is minimal. Proposition \ref{Palgebracontinuous} was based on Proposition \ref{Pcomprepr}
yielding that for each compactification allowing a continuous
representation every bigger compactification has the same property.
The following example shows that this is not true for $T: X \to X$.\\

\begin{exa}\label{Eodometer}
Let $X=\Z$, $T: k \mapsto k+1$, $\alpha \in \R \setminus \Q$,
and consider the compactifications $(\iota_i,K_i)$, $i=1,2$, given by
$K_1:= \T = \R/\Z$, $\iota_1: k \mapsto k\alpha + \Z$, and
$K_2:=[0,1]$, $\iota_2: k \mapsto \{k\alpha\}:=k\alpha-
\max \{m \in \Z:\ m \le k\alpha\}$. Obviously $(\iota_1,K_1) \le (\iota_2,K_2)$
via $\pi: K_2=[0,1] \to \R/\Z = K_2$, $x \mapsto x+\Z$. The compactification $(\iota_1,K_1)$ is compatible with $T$ since $T_1: x \mapsto x+\alpha$
is a continuous extension of $T$ in $(\iota_1,K_1)$.
There is no continuous extension $T_2$ of $T$ in $(\iota_2,K_2)$.
Let us prove this by contradiction, i.e.\ by assuming that such a $T_2: [0,1] \to [0,1]$ exists.
Since $\iota_2(\Z)$ is dense in $[0,1]$ we can find a sequence
$x_n = \{k_n \alpha\} = \iota_2(k_n)$ such that $x_n \to 1 - \alpha$. Furthermore we can arrange $0\le x_n < 1- \alpha$ for all $n\in \N$.
Using $T_2 \circ \iota_2 = \iota_2 \circ T$ we get
\begin{eqnarray*}
T_2(1-\alpha) &=& \lim_{n \to \infty} T_2(x_n) =
      \lim_{n \to \infty} T_2 \circ \iota_2(k_n) =
      \lim_{n \to \infty} \iota_2 \circ T (k_n) = \\
   &=& \lim_{n \to \infty} \iota_2(k_n+1) =
      \lim_{n \to \infty} \{k_n\alpha + \alpha\} =
      \lim_{n \to \infty} x_n + \alpha = 1.
\end{eqnarray*}
We now pick another sequence
$y_n = \{l_n \alpha\} = \iota_2(l_n)$ such that $y_n \to 1 - \alpha$ but
now with the requirement $1- \alpha<y_n\le 1$ for all $n\in \N$.
Similarly, we get
\begin{eqnarray*}
T_2(1-\alpha) &=& \lim_{n \to \infty} T_2 (y_n) =
      \lim_{n \to \infty} T_2 \circ \iota_2(l_n) =
      \lim_{n \to \infty} \iota_2 \circ T (l_n) = \\
   &=& \lim_{n \to \infty} \iota_2(l_n+1) =
      \lim_{n \to \infty} \{l_n\alpha + \alpha\} =
      \lim_{n \to \infty} y_n + \alpha -1 = 0,
\end{eqnarray*}
contradiction.
\end{exa}
We see that taking bigger compactifications does not guarantee that we find
continuous extensions. Nevertheless, for the Stone-{\v C}ech
compactification $(\iota_{\beta},\beta X)$, cf.\ Section \ref{Sstonecech},
everything works. In particular the first statement of
Proposition \ref{Pstonecechsummary} applies to
$K=\beta X$: For every continuous $T: X \to X$ the map $\iota_{\beta} \circ T: X \to \beta X$ is continuous. Therefore there is a continuous $T_{\beta}: \beta X \to \beta X$ such that $T_{\beta} \circ \iota_{\beta} = \iota_{\beta} \circ T$.

\begin{defi}\label{Dbeta}
Let $X$ be a completely regular space and $T:X \to X$ be continuous.
Then $T_{\beta}: \beta X \to \beta X$ denotes the (unique) continuous
extension of $T$ in the Stone-{\v C}ech compactification
$(\iota_{\beta},\beta X)$ of $X$.
\end{defi}

There is no obstacle to consider families of transformations
instead of a single $T$. In order to proceed into this direction
recall the notion of (semi)group actions.

\begin{defi}\label{Dsemigroupaction}
Let $S$ be a semigroup and $X$ a set. A mapping $\alpha: S \times X \to X$,$(s,x) \mapsto \alpha(s,x)$ is called a \textbf{semigroup action} of $S$ on $X$ if $\alpha(s_1,s_2(x)) = \alpha (s_1s_2,x)$ for all $s_1,s_2 \in S$ and all $x \in X$.
\end{defi}

This construction carries over in the obvious way to groups:
\begin{defi}\label{Dgroupaction}
Let $G$ be a semigroup and $X$ a set. A mapping $\alpha: G \times X \to X$,$(g,x) \mapsto \alpha(g,x)$ is called a \textbf{group action} of $G$ on $X$ if it is a semigroup action of $G$ considered as a semigroup and $\alpha(e_G,x)=x$, for the unit element $e_G\in G$ and all $x\in X$.
\end{defi}

For a semigroup action the maps $s^{\alpha}: x\mapsto \alpha(s,x)$ are self-maps of $X$. If we impose a semigroup structure on $X^X$, the set of all maps $f: X \to X$, by using the composition of maps as semigroup operation, a semigroup action  of $S$ on $X$ is nothing else than a homomorphism $\alpha: S \to X^X$, $s\mapsto s^\alpha$. Similarly we can impose a group structure on $\mbox{Sym}(X)=\{f\in X^X: f \mbox{ bijective}\}$. In the group case we have
\[ s^{\alpha}(s^{-1})^{\alpha} = (ss^{-1})^{\alpha} = \mbox{Id}_X =
(s^{-1})^{\alpha}s^{\alpha},\]
thus a group action is a homomorphism $G\to \mbox{Sym}(X)$. So far the set $X$ on which $\alpha$ acts upon carries no structure itself.

\begin{defi} Let $X$ be a topological space and $\alpha: S\times X\to X$ a semigroup action. If $s^\alpha$ is continuous for every $s\in S$ we say that $S$ \textbf{acts by continuous maps} on $X$.
\end{defi}

Let $\alpha$ be a semigroup action of $S$ on $X$ by continuous maps.
Suppose that all $s^{\alpha}: X \to X$ have continuous extensions
$s_K^{\alpha}: K \to K$ in the compactification $(\iota,K)$ of $X$.
As a consequence of Proposition \ref{Pstonecechsummary} we have
$(st)_K^{\alpha} = s_K^{\alpha} \circ t_K^{\alpha}$ for all $s,t \in S$, hence $\alpha_K:S\times K \to K$, $(s,c) \mapsto s_K^{\alpha}(c)$, defines a semigroup action of $S$ on $K$ by continuous maps.

\begin{defi}
The action $\alpha_K$ defined as above is called the
\textbf{extension} of the action $\alpha$ in the compactification
$(\iota,K)$ of $X$. For $(\iota,K) = (\iota_{\beta},\beta X)$
all $s^{\alpha}$ have continuous extensions, denoted by $\alpha_{\beta}$.
\end{defi}

Assume now that $S$ is equipped with a topology for which the semigroup operation $S\times S \to S: (s_1,s_2) \mapsto s_1 s_2$ is
jointly continuous on $S \times S$. Then
$S$ is called a topological semigroup, see also Section \ref{Sopcomp}.

\begin{defi}\label{Dtoposemigroup} Let $S$ be a topological semigroup which acts by contionuous maps on $X$. The semigroup action $\alpha: S\times X\to X$ is called a \textbf{jointly continuous semigroup action} if
$\alpha$ is jointly continuous on $S\times X$.
\end{defi}

In the next Section we will analyze when the extension $\alpha_\beta$ of a jointly continuous semigroup action $\alpha$ is again a jointly continuous semigroup action.

\section{Separate and joint continuity of operations}\label{Sjoint}

Let us now focus on the case of a discrete semigroup $S$.
Then the semigroup operation $\alpha: S \times S \to S$, $(s,t) \mapsto st$
is an action of $S$ on itself which has a continuous extension
$\alpha_{\beta}: S \times \beta S \to \beta S$ to its Stone-{\v C}ech
compactification $(\iota_{\beta}, \beta S)$. All $r_a: S \to \beta S$,$r_a(s):= \alpha_{\beta}(s,a)$, $a \in \beta S,$ have continuous extensions
$\rho_a: \beta S \to \beta S$ fulfilling
$\rho_a(\iota_{\beta}(s)) = r_a(s) = \alpha_{\beta}(s,a)$. Consider now
the operation $*: \beta S \times \beta S \to \beta S$,
$(a,b) \mapsto a*b:= \rho_b(a)$. This operation is described in the
following statement.

\begin{pro}\label{Pstonecechsemigroup}
Let $S$ be a discrete semigroup. Then there is a unique semigroup operation
$*: \beta S \times \beta S \to \beta S$ on the Stone-{\v C}ech
compactification $(\iota_{\beta},\beta S)$ of $S$ such that:
\begin{enumerate}
\item
$*$ extends the semigroup operation on $S$, i.e.\
$\iota_{\beta}(s) * \iota_{\beta}(t) = \iota_{\beta}(st)$
for all $s,t \in S$.
\item
The right translations $\rho_{a}: \beta S \to \beta S$, $x \mapsto x*a$,
are continuous for all $a \in \beta S$.
\item
The left translations $\lambda_{a}: \beta S \to \beta S$, $x \mapsto a*x$,
are continuous for all $a \in \iota_{\beta}(S)$.
\end{enumerate}
\end{pro}
\begin{proof}
It suffices to prove that the operation $*$ defined before the proposition
is associative. For all $s,t,u \in \iota_{\beta}(S)$ we have
\[
\lambda_{st}(u) = (st)u = s(tu) = \lambda_s \circ \lambda_t(u).
\]
Since $\lambda_{st}$ and $\lambda_s \circ \lambda_t$ are continuous
and $\iota_{\beta}(S)$ is dense in $S$ this
equation extends to $\lambda_{st}(z) = \lambda_s \circ \lambda_t(z)$ for
all $z \in \beta S$. Hence
\[
\rho_z \circ \lambda_s(t) = (s*t)*z = \lambda_{st}(z) =
   \lambda_s \circ \lambda_t(z) = s*(t*z) = \lambda_s \circ \rho_z(t).
\]
Since $\rho_z \circ \lambda_s$ and $\lambda_s \circ \rho_z$ are continuous
this equation similarly
extends to $\rho_z \circ \lambda_s(y) = \lambda_s \circ \rho_z(y)$ for
all $y \in \beta S$. Hence
\[
\rho_z \circ \rho_y(s) = (s*y)*z = \rho_z \circ \lambda_s(y) =
   \lambda_s \circ \rho_z(y) = s*(y*z) = \rho_{y*z}(s).
\]
Once again, since $\rho_z \circ \rho_y$ and $\rho_{y*z}$ are continuous
this equation extends to $\rho_z \circ \rho_y (x) = \rho_{y*z}(x)$ for
all $x \in \beta S$, hence
\[
(x*y)*z = \rho_z \circ \rho_y (x) = \rho_{y*z}(x) = x*(y*z)
\]
for all $x,y,z \in \beta S$.
\end{proof}

For much more information about the algebraic structure of $\beta S$  we refer to \cite{HiSt98}. There one can also find information about related constructions as the enveloping semigroup of a semigroup of continuous transformations etc.
We are now going to show that $*$ is jointly continuous only in
very special cases (which are not strikingly interesting).
Not for maximizing generality but in order to identify the natural
context we use terminology from General Algebra.

\begin{defi}\label{Dnaryoperation}
For any set $X$, a function $\omega: X^n \to X$, $n \in \N$, is called an
\textbf{$n$-ary operation on $X$}. If $(\iota,K)$ is a compactification
of $X$, $\omega_K: K^n \to K$ is called an \textbf{extension} of $\omega$ if
$\iota \circ \omega = \omega_K \circ \iota^n$ with
$\iota^n: X^n \to X^n$, $(x_1,\ldots,x_n) \mapsto
\left(\iota(x_1),\ldots,\iota(x_n)\right)$. $(\iota,K)$ is called \textbf{compatible}
with $\omega$ and vice versa if a continuous extension $\omega_K$ of
$\omega$ exists.
\end{defi}

\begin{rem}\label{Rclone}
If $\omega_0$ is an $m$-ary operation on $X$ and
$\omega_1,\ldots,\omega_m$ are $n$-ary operations on $X$ then
$\omega(x_1,\ldots,x_n) := \omega_0(\omega_1(x_1,\ldots,x_n), \ldots,
\omega_m(x_1,\ldots,x_n))$ defines an $n$-ary operation $\omega$ on $X$,
called the composition of $\omega_0$ and the $\omega_i$, $i=1,\ldots,m$.
If all involved operations are continuous then so is $\omega$.
Other (trivial) examples of continuous $n$-ary operations are the projections
$\pi_i^n: (x_1,\ldots,x_n) \mapsto x_i$, $1 \le i \le n$. A set $\Omega$ of
operations on $X$ which contains all projections and is closed under
composition is called a \textbf{clone} on $X$. Thus, for every family of continuous
$n_i$-ary operations $\omega_i$ on $X$, $i \in I$, all operations in
the clone generated by the $\omega_i$ are continuous as well.
A standard reference on the clone of continuous functions is \cite{Tayl86}.
From this point of view the traditional approach in General Algebra,
namely to define a universal algebra as an object of the type
$(X,(\omega_i)_{i \in I})$, is intimately connected with the investigation
of clones. In particular for infinite $X$ there is indeed much
current research on clones on $X$, cf. \cite{GoPi08}. But here we proceed in a different direction.
\end{rem}

Recall that, by definition, a topological space is 0-dimensional if there exists a topological basis of clopen sets.
\begin{lem}\label{Lclopen}
Let $X$ be a 0-dimensional compact Hausdorff space and $R\tm X^n$ a clopen subset.
Then $R = \bigcup_{i=1}^k R_i$ is a finite union of generalized rectangles $R_i = A_{i,1} \times \ldots\times A_{i,n}$ with clopen $A_{i,j} \tm X$,$i=1,\ldots,k$, $j=1,\ldots,n$. The $R_i$ can be taken pairwise disjoint and such that all the $A_{i,j}$ are from a fixed finite partition $P=\{A_1,\ldots,A_k\}$.
\end{lem}
\begin{proof}
Pick $x = (x_1,\ldots,x_n) \in R$.
Since $R$ is open and $X$ has a clopen basis there are clopen
neighborhoods $A_{x,i}$ of $x_i$ such that $A_x:= A_{x,1} \times \ldots\times A_{x,n} \tm R$. Hence $R = \bigcup_{x \in R} A_x$.
This covering is open. Since $R$, being a closed subset of $X$,
is compact, finitely many $R_i:= A_{x_i}$, $i=1,\ldots,k$, form
a covering as well. It is clear that by finite refinements, the $R_i$ can be taken pairwise disjoint and all the resulting $A_{i,j}$ from one finite partition.
\end{proof}

This lemma yields a characterization of operations having a
continuous extension in the Stone-{\v C}ech compactification.

\begin{thr}\label{Tstonecech}
Let $X$ be discrete and $\omega: X^n \to X$ an $n$-ary operation on $X$. Then $\omega$ has a continuous extension in the Stone-{\v C}ech
compactification $(\iota_{\beta},\beta X)$ if and only if for every
$S \tm X$ the preimage is a finite union of rectangles, i.e.\
\[\omega^{-1}[S] = \bigcup_{i=1}^k R_i \quad \mbox{with} \quad
     R_i = A_{i,1} \times \ldots \times A_{i,n}. \]
\end{thr}
\begin{proof}
Necessity of the condition: Assume that $\omega_{\beta}$
is the continuous extension of $\omega$ in $(\iota_{\beta}, \beta X)$.
$S^* = \overline{\iota_{\beta}(S)}$ (notation as in Proposition
\ref{Pstonecechsummary}) is clopen, hence, by continuity of $\omega_{\beta}$,
$\omega_{\beta}^{-1}[S^*]$ is clopen as well. So Lemma \ref{Lclopen}
applies, yielding that this set is a finite union of rectangles.
This immediately translates to the same property of $\omega^{-1}[S]$.

Sufficiency of the condition: Assume that for $\omega: X^n \to X$
all preimages $\omega^{-1}[S]$, $S \tm X$, are finite unions of rectangles. We have to construct a continuous extension $\omega_{\beta}$ of $\omega$ in $(\iota_{\beta},\beta X)$. We use the ultrafilter description from Proposition \ref{Pstonecechsummary}. So let $p_1,\ldots,p_n$ be ultrafilters on $X$.
We define an ultrafilter $p:= \omega_{\beta}(p_1,\ldots,p_n)$ on $X$
by letting $F \tm X$ be a member of $p$ if and only if
$\omega(F_1 \times \ldots \times F_n) \tm F$ for some sets $F_i \in p_i$.
It is straight forward to check that $\emptyset \notin p$,
that $F \in p$ and $F \tm F' \tm X$ implies $F' \in p$ and that
$F,F' \in p$ implies $F \cap F' \in p$. But $p$ is even maximal:
For arbitrary $F \tm X$ our assumption yields that $R:=\omega^{-1}[F]$
can be taken as stated in Lemma \ref{Lclopen}. Since for each $j=1,\ldots,n$, $p_j$ is an ultrafilter on $X$ there is exactly one $A_{k_j}\in p$ such that $A_{k_j}\in p_j$. For the rectangle $R':=A_{k_1}\times\ldots\times A_{k_n}$ we either have $R'\tm R$ or $R'\tm X\setminus R$. In the first case this implies $F\in p$, in the second case $X\setminus F \in p$, showing that $p$ is an ultrafilter.

Finally we have to prove that $\omega_{\beta}$ is continuous on $(\beta X)^n$. We use Proposition \ref{Pstonecechsummary} several times. Take arbitrary ultrafilters $p_1,\ldots,p_n \in \beta X$
and any neighborhood $U$ of $p := \omega_{\beta}(p_1,\ldots,p_n)$.
By the definition of the topology on $\beta X$ there is a set $F \tm X$ such that $F \in p_1$ and $U$ contains all ultrafilters $p$ with $F \in p$. By the definition of $\omega_{\beta}$ there are $F_i \in p_i$ such that $\omega (F_1 \times \ldots \times F_n) \tm F$. Each $F_i$ defines a neighborhood $U_i$ of $p_i$ consisting of all ultrafilters which contain $F_i$. It is clear that $\omega_{\beta}(U_1 \times \ldots\times U_n) \tm U$, showing that $\omega_{\beta}$ is continuous.
\end{proof}

\begin{cor}\label{Cnogroupextension}
Let $S$ be an infinite discrete group. Then there is no continuous
extension of the group operation on $S$ in $(\iota_{\beta},\beta S)$.
\end{cor}
\begin{proof}
Preimages of singletons are infinite but contain only singleton rectangles, thus can not be finite unions of rectangles.
\end{proof}

Similar arguments apply for many semigroups as $\N$ with addition or with multiplication or infinite totally ordered sets with $\min$ or $\max$ as semigroup operation.

\begin{defi}\label{Dessunary}
An $n$-ary operation $\omega: X^n \to X$ is called \textbf{essentially unary} (depending on the $i$-th component) if there is an $f: X \to X$ such that $\omega(x_1,x_2,\ldots,x_n) = f(x_i)$ for all $x_1,\ldots,x_n \in X^n$. $\omega$ is called \textbf{locally essentially unary} if there is a finite partition of $X$ into sets $A_i$, $i=1,2,\ldots,k$,  such that the restriction of $f$ to each rectangle $R = A_{i_1} \times\ldots\times A_{i_n}$, $i_j \in \{1,2,\ldots,k\}$ is essentially unary.
\end{defi}

\begin{pro}\label{Tlocallyunary}
Let $\omega: X^n \to X$ be locally essentially unary.
Then there is a continuous extension $\omega_{\beta}$ of $\omega$
to $(\iota_{\beta},\beta X)$.
\end{pro}
\begin{proof}
As the reader checks easily, every locally unary operation
$\omega$ satisfies the condition of Theorem \ref{Tstonecech}.
\end{proof}

Continuing work of van Douwen \cite{Douw81}, Farah was able to show in \cite{Fara01, Fara02}
that the converse of Proposition \ref{Tlocallyunary} also holds true

\begin{pro}[Farah]
Let $X$ be an infinite discrete set and assume that $\omega: X^n \to X$ has a continuous extension in $(\iota_{\beta}, \beta X)$. Then $\omega$ is locally essentially unary.
\end{pro}

\section{Compactifications for operations}\label{Sopcomp}

We have seen in the previous section that many interesting binary
operations can not be extended continuously to the Stone-{\v C}ech
compactification. Nevertheless some ideas presented in Section
\ref{Scompcont} can be adapted. In order to be more flexible it is
useful to consider the following setting.

\begin{defi}\label{Dsemitopological}
Let $I$ be an index set, $n_i \in \N$ and $\gamma_i\subseteq \P(\{1,2,\ldots,n_i\})$ for all $i \in I$.
A \textbf{semitopological} (general) \textbf{algebra} of type
$\tau=((n_i)_{i \in I},(\gamma_i)_{i \in I})$ is a topological space
$X$ together with a family of $n_i$-ary operations $\omega_i: X^{n_i} \to X$
for which $(x_{j_1},\ldots,x_{j_s}) \mapsto \omega_i(x_1,\ldots,x_{n_i})$
is continuous for all $\{j_1,\ldots,j_s\} \in \gamma_i$ and all fixed $x_i\in X$, $i\notin\{j_1,\ldots,j_s\}$.
This semitopological algebra is called a \textbf{topological algebra}
if furthermore $\{1,\ldots,n_i\} \in \gamma_i$ for all $i \in I$.
In this case one might omit the information contained in the
$\gamma_i$ and consider $\tau$ to be given by the $\tau=(n_i)_{i\in I}$.
\end{defi}

\begin{exa}[Semitopological algebras]\label{Etopoalgebras}
\hfill\par\vspace{-2ex}
\begin{itemize}
\item
Topological groups are groups
considered as topological algebras of type $\tau = (2,1)$, requiring joint continuity of the binary operation as well as continuity of the operation $x\mapsto x^{-1}$.
 \item Topological groups can also be seen as topological algebras of type $\tau = (2,1,0)$ if one prefers to emphasize that the neutral element may be considered as a $0$-ary operation.
\item Topological semigroups are semigroups
which are topological algebras of type $\tau = (2)$.
\item Semitopological semigroups are semigroups considered as semitopological algebras
of type $(2,\{\{1\},\{2\}\})$, i.e.\ the semigroup operation
is continuous in each component but not necessarily jointly continuous.
Similarly the type of left and right topological semigroups
is $\tau = (2,\gamma)$ with $\gamma = \{\{i\}\}$ with $i=1$ resp.\ $i=2$.
\end{itemize}
\end{exa}

The value of the rather technical concept of a semitopological algebra
gets clear by considering compactifications of general algebras.

\begin{defi}\label{Dcompalgebra}
Let $X$ and $K$ be (semi)topological algebras of type $\tau$. If
$(\iota,K)$ is a compactification of the set $X$ such that each operation
on $K$ extends the corresponding operation on $X$ we call $(\iota,K)$
a \textbf{$\tau$-compactification of $X$}. In the case of topological groups, (semi)topological semigroups
etc.\ these compactifications are also called group, (semi)topological
(semi)group etc.\ compactifications in the obvious way.
\end{defi}

Later we will discuss the special cases of group, semigroup
and semitopological semigroup compactifications in more detail.
In the general context the following observations hold.

\begin{pro}\label{Psemitopologicalalgebrafacts}\hfill\par\vspace{-2ex}
\begin{enumerate}
\item The direct product of a family of (semi){-}topological algebras
of type $\tau$ is again a (semi)\-topological algebra of type $\tau$.
\item Every semitopological algebra of type $\tau$ has a maximal
$\tau$-compactification.
\end{enumerate}
\end{pro}
\begin{proof}
The first statement is obvious. For the second statement
the product compactification, cf.\ Definition \ref{Dprodcomp},
of all (semi)topological compactifications of type $\tau$
has the required properties. To justify this construction it suffices to show that there is a set ${\cal S}$ of compactifications of $X$
such that for every $\tau$-compactification $(\iota,K)$ of $X$
there is an equivalent compactification in ${\cal S}$. Since
$|K| \le |\beta X|$ one can take for ${\cal S}$ the set of all
compactifications $(\iota,K)$ of $X$ with $K \tm \beta X$
(as a set, not necessarily as a topological subspace or subalgebra).
\end{proof}

\begin{exa}\label{EBohr}
\hfill\par\vspace{-2ex}
\begin{enumerate}
\item
For a topological group $G$ the maximal group compactification
is called the almost periodic or \textbf{Bohr compactification}
and denoted by $(\iota_b,bG)$.
\item
For a semitopological semigroup $S$ the maximal semitopological
semigroup compactification is called the \textbf{weak almost periodic compactification} and denoted by $(\iota_w,wS)$, see also Section \ref{Sweakhartman}. For the realization of $wS$ as space of filters in the spirit of Proposition \ref{Pstonecechsummary} we refer to \cite{BeHi84}.
\end{enumerate}
\end{exa}

\section{Invariance on groups and semigroups}\label{Sinvgroups}
\subsection{The action of a semigroup by translations}
With every (semi)group $S$ comes a natural action, namely the action of $S$ by right translations.

\begin{defi} Let $S$ be a semitopological semigroup. Then $S$ acts on $B(S)$ by right translations in the following way (notation as in Proposition \ref{Pstonecechsemigroup}):
\[R: S\times B(S) \to B(S),\quad R_s(f)(t):=f(\rho_s(t))= f(ts).\]
\end{defi}
For every $s\in S$ the map $R_s$ is a bounded linear operator.
As for the left translations
\[L: S\times B(S) \to B(S),\quad L_s(f)(t):=f(\lambda_s(t))= f(st)\]
we have $L_s L_t=L_{ts}$ the map $(s,f)\mapsto L_f $ is not a semigroup action, but merely an "anti"-action of $S$ on $B(S)$. However, in the group case we can define an action $\Lambda: G\times B(G) \to B(G)$ via $\Lambda_s(f)(t)=f(\lambda_{s^{-1}}t)$ after all.

In the sequel we use compatibility with respect to these translations to single out a unique measure or mean on certain algebras $\A\subseteq B(S)$. We will mainly focus on the group case. As in contrast to the previous sections this section will be less self-contained. As standard references  (which also extensively treat the semigroup case) we mention \cite{BeJM89,Burc70,Gree69,Pate88,Rudi90,Rudi91,Rupp84}.

\subsection{Means}
\label{SSmeans}
\begin{defi}\label{Dinvmean}
Let $\A\subseteq B(S)$ be a $*$-algebra which is invariant under translations. A mean $m\in M(\A)$ is left (right) invariant, if $m(f)=m(L_sf)$ $(m(f)=m(R_sf))$ for all $s\in S$ and $f\in \A$. A mean which is both left and right invariant is called bi-invariant, or, simply invariant.
\end{defi}

It is a nontrivial task to find conditions on $\A$ which ensure the existence of an invariant mean. It turns out that the closure of the convex hull of the orbits with respect to various topologies play an important r\^ole.

\begin{pro}\label{Pinvmean}
Let $S$ be a semitopological semigroup and $\A\subseteq B(S)$ a $C^*$-algebra such that
\begin{enumerate}
\item
for each $f\in \A$, $\overline{\mbox{co}}(\{L_sf: s\in S\})$ contains a constant,
\item
for each $f\in \A$, $\overline{\mbox{co}}(\{R_sf: s\in S\})$ contains a constant,
\end{enumerate}
where $\overline{\mbox{co}}$ indicates the closure of the convex hull w.r.t\ uniform convergence. Then there exists a mean $m\in M(\A)$ which is bi-invariant. Furthermore $m$ is unique.
\end{pro}
The proof of this assertion can be found for example in \cite{Burc70}.
It is a general principle in the theory of function spaces on semigroups that constants in convex closures are intimately linked to invariant means, see \cite[Chapter 2]{BeJM89}. Proposition \ref{Pinvmean} states that if we can find constants in the uniform closures of the convex hull, there exists already a corresponding bi-invariant mean which is unique.

Having Proposition \ref{Dinvmean} at hand, we can establish the existence and uniqueness of an invariant mean for (weak) almost periodic functions. Recall the notion of almost periodicity:

\begin{defi}\label{Dweakap}
A bounded function $f: S \to \C$ on a semitopological semigroup $S$ is called \textbf{(weak) almost periodic}, if the set of left translations $\{L_s f: s\in S\}$ is relatively compact in the norm (weak) topology. Let us denote the algebra of almost periodic functions by $AP(S)$ and the algebra of weak almost periodic functions by $\W(S)$.
\end{defi}
Evidently $AP(S)\subseteq \W(S)$. Weakly almost periodic functions may be characterized using the following double limit criterion.

\begin{pro}[Grothendieck]\label{Tdoublelimit}
Let $S$ be a semitopological semigroup. A bounded function $f:S\to \C$ is weakly almost periodic if and only if
\[\lim_{n\to\infty} \lim_{m\to\infty} f(t_ns_m)=\lim_{m\to\infty}\lim_{n\to\infty}f(t_ns_m)\]
whenever $(t_n)_{n=1}^{\infty}, (s_m)_{m=1}^{\infty} \subseteq S$ are sequences such that the involved limits exist.
\end{pro}
Using weak compactness of the translation orbit it is easy to check that weakly almost periodic functions satisfy the double limit condition. The complete proof can be found for instance in \cite{Krem86}.

\begin{rem}[Weak almost periodicity]
\hfill\par\vspace{-2ex}
\begin{itemize}
\item
This definition of (weak) almost periodicity does not depend on the given topology on $S$ since the norm (weak) topology on $C_b(S)$ coincides with the relative topology inherited from the norm (weak) topology on $B(S)=C_b(S_{dis})$. The statement is obvious for the norm topology; for the weak topology it follows from the Hahn-Banach Theorem.

\item For any (weak) almost periodic function the right orbit
$\{R_s f: s\in S\}$ is relatively (weak) compact as well. However, in the weak case left and right orbit closures will in general not coincide, while for almost periodic functions this always is the case.

\item  One can show that the set of (weak) almost periodic functions is a $C^*$-algebra. Its structure space (see also Section \ref{Scompcont}) is a topological group (semitopological semigroup). The structure space coincides with the \textbf{Bohr compactification} in the almost periodic case resp.\ with the\text{weak almost periodic compactification} in the weak case.
\end{itemize}
\end{rem}
 Before we go on, we quote the celebrated fixed point Theorem of C.\ Ryll-Nardzewski which is vital to the theory of almost periodicity
\begin{pro}[Ryll-Nardzewski] Let $X$ be a Banach space and $K\subseteq X$ a weakly compact convex set. Let $S$ be a semigroup which acts on $K$ by affine mappings, i.e.\ for each $s\in S$ there
is a linear operator $A_s: X\to X$ and an element $x_s \in X$ such that $T_s:=x_s+A_x$ and such that $T_tT_s=T_{ts}$ for all $s,t\in S$.
If, furthermore, $\inf_{s,t \in S} \|T_sx-T_tx\|>0$ for every $x\in K$, then there exists a common fixed-point.
\end{pro}

\begin{proof} See the original paper \cite{Ryll67} for a probabilistic or \cite{Gree69} for a geometric proof.
\end{proof}

\begin{pro}\label{Pwapmean}
Let $S$ be a semitopological semigroup. There exists a unique invariant mean on $\W(S)$ and hence also on $AP(S)$.
\end{pro}
We only give a sketch of the argument for the case where $S=G$ is a group. We will employ the Ryll-Nardzewski fixed-point Theorem to show that $\W(S)$ meets the requirements of Proposition \ref{Pinvmean}. Observe that a function $f: S\to S$ satisfying $L_sf=f$ for all $s\in S$ has to be constant. The weak closure of the convex hull of a weakly compact set is again weakly compact (Krein-\v{S}mullyan Theorem). So for $f\in \W(S)$ the set  $K:=\overline{\mbox{co}}^{(w)}(\{L_sf: s\in S\})$ is again weakly compact. As $K$ is convex, the norm closure and the weak closure coincide (Hahn-Banach Theorem). The action of $S$ by the translations $L_s$ leaves the set $K$ invariant, in fact $S$ acts by linear isometries. Thus we can use the Ryll-Nardsewski Theorem to conclude that there exists a common fixed-point, i.e.\ a constant.

\begin{exa}
\hfill\par\vspace{-2ex}
\begin{itemize}
\item Let $\chi: S \to \C$ be a (semi-)character, i.e.\ a continuous  (semi-)group homomorphism which satisfies $|\chi(s)|= 1$. Then $\chi$ is almost periodic. For the invariant mean on $AP(S)$ we have
    \[m(\chi)=m(L_s\chi)=\chi(a)m(\chi).\]
    Thus, if $\chi$ is not the constant character $\Eins_S$, then  $m(\chi)=0$.
\item Let $S$ be locally compact. Every $f\in C_0(S)$ is weakly almost periodic and $m(f)=0$ for the unique mean $m$ on $\W(S)$.
\item Let $S=G$ be a locally compact abelian (LCA) group. Then the Fourier-transform $\hat{\mu}: \hat{G} \to \C$ of a Borel measure $\mu$ on $G$
    \[\hat{\mu}(\chi):=\int_{G} \chi d\mu,\]
    is weakly almost periodic and $m(\hat{\mu})=\mu(\{e_G\})$ for the unique mean $m$ on $\W(\hat G)$ and $e_G$ the neutral element of $G$. By Bochner's Theorem every positive definite function is weakly almost periodic. Recall that a function $f: \hat G \to \C$ is positive definite if for all $\chi_1,\ldots, \chi_n \in \hat G$ the matrix $\left(f(\chi_i \overline{\chi_j})\right)_{i,j=1}^n\in \C^{n\times n}$ is positive definite.
\end{itemize}
\end{exa}

\subsection{Measures}
\label{SSmeasures}
Proposition \ref{Pwapmean} takes a particularly nice form if $S=G$ is a compact topological group.

\begin{pro}[Haar measure]\label{haar}
Let $G$ be a compact topological group. Then there exists a regular Borel probability measure $\mu$ on $G$ which is invariant under left- and right-translations, i.e.\ $\mu(A)=\mu(gA)=\mu(Ag)$ for every Borel set $A\subseteq G$ and $g\in G$. This measure is unique and called the \textbf{Haar measure}.
\end{pro}
\begin{proof} For compact $G$ and continuous $f:G\to \C$ the map $L^f: G \to B(G)$, $g\mapsto L_g f$ is continuous in the norm (weak) topology on $G$. Thus the norm (weak) closure of $\mbox{co}\left(\{L_g f: g\in G\}\right)$ is compact, i.e.\ $f$ is (weak) almost periodic;
$AP(G)=\W(G)=C_b(G)$. Proposition \ref{Pwapmean} together with Riesz' Representation Theorem \ref{riesz} yields the existence of a unique left-invariant measure $\mu_m$. Note that $\mu_m(G)=m(\Eins_G)=1$, so $\mu$ is a probability measure.

Inversion $g\mapsto g^{-1}$ turns any left invariant measure (mean) on $G$ into a right invariant measure (mean). Since the unique mean on $AP(G)$ is bi-invariant the left invariant Haar measure on a compact group is also right invariant.\end{proof}

\begin{rem}
If $G$ is only \emph{locally compact} we can still construct a left invariant measure. In this setting however, uniqueness holds only up to a multiplicative constant and left invariance does in general not imply right invariance. For a rigorous treatment of the Haar measure we refer to
\cite{HeRo79,Stro06}
\end{rem}

We can use the Haar measure to give an alternative approach to the unique invariant mean on the (weak) almost periodic functions. As the Bohr compactification $(\iota_b, bG)$ of a topological group $G$ is compact, there exists the Haar measure $\mu_b$ on $bG$. $bG$ is the structure space of the $C^*$-Algebra $AP(G)$, thus there is an isomorphism $C(bG)\cong AP(G)$ given by $F\mapsto F\circ \iota$, see Section \ref{Scompcont}. Consequently $m(F\circ\iota):=\int_{bG} F d\mu_b$ defines an invariant mean on $AP(G)$, which, by uniqueness of the Haar measure, is the unique invariant mean on $AP(G)$.

\begin{cor}\label{Capnull}
Let $G$ be a topological group and $f\in AP(G)$. If $m(|f|)=0$ for the unique invariant mean $m$ on $AP(G)$, then $f=0$. \end{cor}
\begin{proof} As the Haar measure gives positive measure to open sets in $bG$ the equality $$0=m(|F\circ\iota|)=\int_{bG} |F| d\mu_b, \quad F\in C(bG)$$ can only hold if $F=0$.\end{proof}

\begin{cor}
Let $G$ be a topological group and let $(C_1,\iota_1)$, $(C_2,\iota_2)$ be group compactifications of $G$. If $(C_1,\iota_1)\le (C_2,\iota_2)$ via $\pi: C_2 \to C_1$ then $\mu_1 = \pi \circ \mu_2$, where $\mu_1$ resp.\ $\mu_2$ is the Haar measure on $C_1$ resp.\ $C_2$. In particular Proposition \ref{Pcompmeasure} applies.
\end{cor}
\begin{proof} Using continuity of $\pi$ and density of $\iota_i(C) \in C_i$ it is straight-forward to check that $\pi \circ \mu_2$ is an invariant Borel measure on $C_1$, hence it must be the Haar measure.
\end{proof}

\subsection{Amenability}
\label{SSamenable}
Finally we will drop any uniqueness assumptions and focus on (semi)\-groups such that there exists at least one invariant mean.
\begin{defi}
A discrete (semi)group $S$ is called \textbf{amenable} if there exists a bi-invariant mean on $B(S)$.
\end{defi}
In the group case one can use the inversion $g\mapsto g^{-1}$ to show that existence of a one-sided invariant mean is equivalent to the existence of a bi-invariant mean. Indeed, let $m_l, m_r: B(X)\to \C$ be left resp.\ right-invariant. For bounded $f: G\to G$ the function $M_lf(g):=m_l(R_g f)$ is again bounded. Then
$m(f):= m_r(M_lf)$ defines a bi-invariant mean. Furthermore, if $G$ carries a locally compact topology then the existence of an invariant mean on $UCB(G)$, the algebra of uniformly continuous functions, implies the existence of an invariant mean on $B(G)$ and vice versa. For details we refer to \cite{Gree69,Pate88}.

As we have seen in the previous section every group which admits a compact group topology is amenable as the Haar measure defines an invariant mean on $C_b(G)$. Among many other more or less remarkable properties amenable groups have interesting dynamical behavior.
\begin{pro}[Markov-Kakutani]\label{PMKfixpoint}
Let $S$ be amenable. If $S$ acts on a compact convex subset $K$ of a topological vector space $X$ by affine mappings, then there exists a common fixed point.
\end{pro}
The proof of this statement can be found in \cite{Gree69}.
Note that the classical Markov-Kakutani fixed-point theorem is stated for the case $S=\Z$, or, slightly more general, for abelian $S$.

The class $\textbf{Am}$ of amenable \emph{groups} is closed under elementary group theoretical constructions, i.e.\ if $G$ is amenable, then so is every subgroup and homomorphic image of $G$. Similarly, let $G_1$ and $G_2$ be amenable, then $G_1\times G_2$ is amenable, also the extension of an amenable group by an amenable group is amenable. Finally, also directed unions of amenable groups are amenable.

It is trivial that every finite group is amenable and it is well-known that every abelian group is amenable. The class $\textbf{ElAm}$ of elementary amenable groups is defined as follows: all finite groups and all abelian groups belong to $\textbf{ElAm}$ and $\textbf{ElAm}$ is closed under taking subgroups, homomorphic images, finite direct products, extensions and direct unions. One might ask whether $\textbf{ElAm} = \textbf{Am}$. The answer to this question is negative. For details see \cite{Pate88}.

In Section \ref{SSfreegroup} we have seen that the free group $F_2:=F(x,y)$ is not amenable. Consequently no group containing $F_2$ as a subgroup (such as $SO(n)$ for $n\ge 3$) can be amenable. Define the class $\textbf{NFree}$ of groups which do not contain $F_2$ as subgroup. The so called \emph{von Neumann conjecture} states that $\textbf{AM} = \textbf{NFree}$. However, this long standing conjecture was proven to be wrong by Ol'shanskii, for details see again \cite{Pate88}.

We conclude this detour to amenable groups with the notion of extreme amenability. A group is called \textbf{extremely amenable} if there exists a \emph{multiplicative} invariant mean. Extremely amenable groups arise as transformation groups of infinite dimensional Hilbert spaces. They are intimately linked to concentration of measure phenomena; compact groups which are extremely amenable must be trivial (uniqueness of the Haar measure), but also locally compact groups can never be extremely amenable, see \cite{GiMi00a, Uspe02}.

\chapter{Hartman measurability}\label{CHhartman}
\section{Definition of Hartman functions}
\label{Shartman}

The following definition fixes the main objects for the rest of the paper.
\begin{defi}\label{Dhartman}
Let $G$ be a topological group. We call
a bounded function $f: G\to \C$ \textbf{Hartman measurable} or a \textbf{Hartman function}
if $f$ can be extended to a Riemann integrable
function on some group compactification. By $\H(G)$ we denote the set of Hartman functions, by $\Hartset(G)$ the system of \textbf{Hartman sets}
\[ \{\iota_b^{*\,-1}[A]: A\in \mathfrak{C}_{\mu_b}(bG)\}.\]
\end{defi}

According to Proposition \ref{Pcompmeasure} such a compactification
can always be taken to be the maximal one, i.e.\ the Bohr compactification $(\iota_b, bG)$.  The Haar measure on $bG$ is denoted by $\mu_b$. Let $\mathfrak{C}_{\mu_b}(bG)$ denote the $\mu_b$-continuity sets on the Bohr compactification, see Definition \ref{DJordan}. Furthermore, it is easy to verify that $\Hartset(G)$ is a set algebra on $G$. We define a fapm $p$ on $\Hartset(G)$ via $p(\iota_b^{*\,-1}[A]):=\mu_b(A)$. $p$ is well-defined by Proposition \ref{Priemannmean}.

\begin{pro}\label{PHartman}The following assertions are equivalent:
\begin{enumerate}
\item $f\in\H(G)$, i.e.\ by definition $f=F\circ \iota$ with  $F\in \Riemanntop_{\mu_K}(K)$, $\Riemanntop_{\mu_K}(K)$ denoting the set of all $F: K\to\C$ which are Riemann-integrable w.r.t.\ the Haar measure $\mu_K$ on $K$, for some group compactification $(\iota, K)$ of $G$.
\item $f=F\circ \iota_b$ with $F\in \Riemanntop_{\mu_b}(bG)$.
\item $f\in B(\Hartset(G)$.
\end{enumerate}
Furthermore, if $\iota_b:G \to bG$ is one-one (1), (2) and (3) are equivalent
to
\begin{enumerate}
\item[4.] $f\in \overline{AP(G)}^{(m)}$, the $m$-completion of the almost periodic functions with respect to the unique invariant mean $m$.
\end{enumerate}
\end{pro}
\begin{proof}
(1)$\Leftrightarrow$(2): Apply Proposition \ref{Pcompmeasure}.

(2)$\Leftrightarrow$(3): Consider the map $\iota_b^*: F\mapsto F\circ\iota_b$ which sends a function defined on the  Bohr compactification $bG$ to a function defined on the group $G$. $\iota_b^*$ maps $\Riemanntop_{\mu_b}(bG)$, the set of Riemann integrable functions on $bG$ (Definition \ref{Driemannmu}) onto $\H(G)$. Thus $\H(G)= \iota_b^* \Riemanntop_{\mu_b}(bG)$.
The map $\iota_b^*: \Riemanntop_{\mu_b}(bG) \to B(G)$ is a bounded homomorphism of $*$-algebras as the reader may quickly verify. Consequently its image, $\H(G)$, is a $C^*$-algebra, see \cite[Theorem I.5.5]{Davi96}. In particular $\H(G)$ is closed.

Recall that $\mathfrak{C}_{\mu_b}(bG)$ denotes the set algebra of $\mu_b$-continuity sets on the Bohr compactification, cf. Definition \ref{DJordan}.  We then have the inclusions $\Simple_{\tiny\Hartset}\subseteq \H(G)\subseteq B(\Hartset(G))$.
The first inclusion is valid since due to linearity of $\iota_b^*$ every $f\in \Simple_{\tiny\Hartset}$ is of the form $F\circ\iota_b$ for some $\mathfrak{C}_{\mu_b}$-simple function $F$. The second inclusion is true by the following argument: $f\in \H(G)$ if there are
$\mathfrak{C}_{\mu_b}$-simple functions $F_n$ such that
$\limn \| F_n\circ \iota_b^*-f\|_{\infty}=0$.  Every function $F_n\circ \iota_b^*$ is $\Hartset$-simple, thus $f$ is in the uniform closure $B(\Hartset(G))$.
Since $\H(G)$ is closed we have
$\H(G)= \overline{\Simple_{\tiny\Hartset}}=B(\Hartset(G))$ using the notation of Section \ref{Salgebras}.

(3)$\Leftrightarrow$(4): Apply Proposition \ref{Pmclosed}.
\end{proof}

In \cite{Hart61} Hartman has used the $m$-closure of almost periodic functions  to define a class of functions called "R-fast\-periodisch" ("R-almost periodic"). According to \cite{Hart61} this nomenclature was suggested by C.\ Ryll-Nardzewski. In our terminology the R-almost periodic functions
coincide with $\H(\R)$, the Hartman functions on the reals.

The equivalence of (1) and (4) in Proposition \ref{PHartman} for $G=\Z^n, \R^n$ has independently been obtained by J.-L.\ Mauclaire (oral communication). In \cite{Mauc86,Mauc88} J.-L.\ Mauclaire used extensions of arithmetic functions to (semi)group compactifications to prove number-theoretic results.

While the inclusion $\H(G)\subseteq \overline{AP(G)}^{(m)}$ is always valid, the converse does not hold true. The crucial property is injectivity of the map $\iota_b: G\to bG$. Topological groups where $\iota_b$ is one-one are called maximally almost periodic.

\begin{exa}\label{Emap}
 Let $G_1$ be a topological group such that $bG_1=\{e\}$ is trivial (such groups exist in abundance, cf. \cite{GlNe01,Remu88,Roth80}, they are called \emph{minimally almost periodic}) and $G_2=\T$, the torus. Denote by $\mu$ the Haar measure on $\T$. Consider $G:=G_1\times G_2$. Then $\iota_b: G\to bG$ is the projection onto the second factor $\iota_b: (x,y)\mapsto y$.    Consider $\A=AP(G)$ and $m$ the unique invariant mean on $AP(G)$. Then $\Riemanntop_{\mu_b}(G)=C(´\{e\}\times \T)$, i.e.\ a function $f:G_1\times G_2 \to \C$ belongs to $\H(G)$ if $f(x,y)=F(y)$ for a function $F\in \H(\T)$.
    Let $F: \T \to \R$ be Riemann integrable such that
    \[F_{\flat}:=\sup_{F_1 \le F \atop F_1\in C(\T)} F_1\neq \inf_{F_2 \ge F \atop F_2\in C(\T)} F_2=: F^{\sharp}.\]
    We can take for instance $F=\Eins_A$, where $A$ is the Cantor middle-third set (in this case $F_{\flat}=0$ and $F^{\sharp}=F$).
    Pick any $y_0\in \T$ such
    that $\alpha:=F_{\flat}(y_0)<F^{\sharp}(y_0)=:\beta$ and pick any non-constant function $F_0: G_1 \to [\alpha,\beta]$. Define
    \[ f(x,y)=\left\{ \begin{array}{ccl} F(y)& \mbox{for} & y \neq y_0 \\
    F_0(x)& \mbox{for} & y = y_0.\end{array}\right.\]
    Then $F\in \overline{AP(G)}^{\,m}$ (since $F_{\flat}(y)\le f(x,y)\le F^{\sharp}(y)$), but $F\notin \H(G)$.
\end{exa}

\begin{pro}\label{PHartmanfacts}
$\H(G)$ is a translation invariant  $C^*$-subalgebra of $B(G)$ and there exists a unique invariant mean on $\H(G)$.
\end{pro}
\begin{proof}

Translation invariance is a consequence of the fact that $\iota_b$ is a group homomorphism. In Proposition \ref{PHartman} we have already seen that $\H(G)$ is a $C^*$-algebra.

Every mean $m$ on $\H(G)$ lifts to a mean $m_b$ on $\Riemanntop_{\mu_b}(bG)$ via the definition $m_b(F):=m(F\circ\iota_b^*)$ for $F\in \Riemanntop_{\mu_b}(bG)$. For invariant $m$ one has
$m_b(F) = \int_{bG} F d\mu_b $ for all \emph{continuous} $F: bG \to bG$ (Riesz' Representation Theorem \ref{riesz} and uniqueness of the Haar measure). Since $\Riemanntop_{\mu_b}(bG)$ is the $\mu_b$-closure of $C(bG)$, $m_b$ is not only unique on $C(bG)$ but also on $\Riemanntop_{\mu_b}(bG)$. This settles the uniqueness of $m$. On the other hand
\[m(F\circ \iota_b):=\int_{bG} F d\mu_b, \quad F \in \Riemanntop_{\mu_b}(bG)\] defines such an invariant mean on $\H(G)$.
\end{proof}

In light of Proposition \ref{PHartman} the fapm $p$ on $\Hartset(G)$ resp. the mean $m$ on $\H(G)$ has a nice completeness property.
\begin{cor} \label{Chartmancomplete}
Let $G$ be a topological group such that $\iota_b: G \to bG$ is one-one
\begin{enumerate}
\item Let $A\in \Hartset(G)$ be a null-set, i.e.\ $p(A)=0$. If $B\subseteq A$ then $B\in \Hartset(G)$.
\item Let $f\in \H(G)$ be a function with zero absolute mean-value, i.e.\ $m(|f|)=0$. If $f: G\to \C$ is such that $|g|\le |f|$ then $g\in \H(G)$.
\end{enumerate}
\end{cor}

\section{Definition of weak Hartman functions}
\label{Sweakhartman}
We need some results concerning the weak almost periodic
compactification $(\iota_q, wS)$ of a semitopological semigroup $S$.
Recall from Section \ref{Sopcomp} that a semitopological semigroup $S$ is a semigroup where all left translations $\lambda_s: S\to S$ and all right translations
$\rho_s: S\to S$ are continuous.

\begin{defi}\label{Dwapcomp}
Let $S$ be a semitopological semigroup. By Proposition \ref{Psemitopologicalalgebrafacts} there exists a maximal compactification $(\iota_w, wS)$ which is a semitopological semigroup. $(\iota_w, wS)$ is called the \textbf{weak almost periodic compactification} of $S$.
\end{defi}

$wS$ is a compact semitopological semigroup which contains $S$
as dense sub-semigroup.
\begin{cor}
Let $S$ be an abelian semitopological semigroup. Then $wS$, the weak almost periodic compactification of $S$, is also abelian.
\end{cor}
\begin{proof}
For every $s\in S$ the continuous maps $\lambda_{\iota_w(s)},\rho_{\iota_w(s)}: wS\to wS$ coincide on the
dense set $\iota_w(S)$. Therefore $\lambda_{\iota_w(s)}=\rho_{\iota_w(s)}$. For arbitrary $x\in wS$ and $s\in S$ we have \[\lambda_x(\iota_w(s))=\rho_{\iota_w(s)}(x)\stackrel{!}{=}
\lambda_{\iota_w(s)}(x)=\rho_x(\iota_w(s)),\]
thus also $\lambda_x$ and $\rho_x$ coincide on a dense set and therefore are equal.
\end{proof}

\begin{defi}
Let $S$ be a semigroup. A subset $I\subseteq S$ of is called a (two-sided) \textbf{ideal} if $\lambda_s(I)\subseteq I$ and $\rho_s(I) \subseteq I$ for every $s\in S$. The by $K(S)$ we denote the \textbf{kernel} of $S$, i.e.\ the intersection of all ideals in $S$.
\end{defi}

From now on we will stick to the special case that $S=G$ is an algebraically an abelian group. Here the kernel $K(G)$ has particularly nice properties.

\begin{pro} Let $G$ be a semitopological abelian group. Then the kernel $K(G)$ is a compact topological group.
\end{pro}
The proof of this assertion can be found in \cite{BeJM89,Rupp84}.

Let $e\in G$ denote the neutral element of the group $K(wG)$. Then $K(wG) = e+wG$ and the mapping $\rho: wG \to K(wG)$ defined via $x\mapsto e+x$ is a continuous retraction, i.e.\ $\rho(x)=x$ for all $x\in K(wG)$.

\begin{pro}\label{Pwapbohr}
Let $G$ be an abelian topological group and $(\iota_w,wG)$ the weakly almost periodic compactification of $G$. Then the compactification $(\rho\circ\iota_w, K(wG))$ is equivalent to the Bohr compactification of $G$.
\end{pro}
\begin{proof}
Note that $(\rho\circ\iota_w, K(wG))$ is a group compactification. We show that it has the universal property of the Bohr compactification.
Each almost periodic function $f$ on $G$
may be extended to a continuous function $F$ on $wG$. Consider the function $F-F\circ \rho$. Since $F\circ\rho$ may be regarded as a continuous function on the group compactification $(\rho\circ\iota_w, K(wG))$ the function $|F-F\circ\rho|$
induces a nonnegative almost periodic function on $G$. Since
this function vanishes on $K(wG)$ the induced almost periodic
function has zero mean-value (note that the mean-value is given by
integration over $K(wG)$ with respect to $\mu_b$). By Corollary \ref{Chartmancomplete} and continuity this implies $F=F\circ\rho$.
Thus we have $f=F\circ\iota_w=(F\circ\rho)\circ\iota_w$. So $F\circ
\rho$ is a continuous extension of $f$ on $(\rho\circ\iota_w, K(wG))$.
\end{proof}

Similarly one proves that for an arbitrary semitopological semigroup compactification $(\iota, C)$ of $G$
the kernel $K(C)$ is a compact topological group and coincides with $\rho(C)=e+C$, where $e$ is the neutral element of $K(C)$ and $\rho$ the retraction defined as above. In this setting, $(\rho\circ\iota, K(C))$ constitutes a group compactification of $G$.

\begin{pro}\label{Pwapmeasure}
Let $G$ be an abelian topological group and $(\iota, C)$ a semitopological semigroup compactification. Then there exists a unique translation invariant Borel measure on $C$.
\end{pro}
\begin{proof}
Suppose $\mu$ is an invariant measure on $C$, then $\mu(C)=\mu(e+C)=\mu(K(C))$ implies
that $\mu$ is supported on the compact
group $K(C)$. Consequently $\mu_{|K(C)}$ coincides with the Haar measure
on $K(C)$. Thus
\begin{equation}\mu(A)=\mu_b(A\cap K(C)), \label{EqExtension}\end{equation}
where $\mu_b$ denotes the Haar measure
on $K(C)$. This settles the uniqueness of $\mu$. Since equation
(\ref{EqExtension}) indeed defines a translation invariant Borel measure also the existence is guaranteed.
\end{proof}

Let us use the framework of semigroup compactifications to
define weak Hartman functions:

\begin{defi}\label{Dweakhartman}
Let $G$ be an abelian topological group and $f: G \to \C$ a bounded function.
Then $f$ is called \textbf{weak Hartman measurable} or a \textbf{weak Hartman function}
if there exists a semitopological semigroup compactification $(\iota,C)$ of $G$
and an $F \in \Riemanntop_{\mu_C}(C)$ such that
$f = F \circ \iota$ for the unique translation invariant measure $\mu_C$ on $C$. The set of all weak Hartman functions on $G$ is denoted by $\H^w(G)$.
\end{defi}

It is almost, but not quite, entirely analog to the strong case to check that
$\H^w(G)$ is a translation invariant $C^*$-subalgebra of $B(G)$ on which a unique invariant mean exists. Furthermore the universal property
of the weakly almost periodic compactification $(\iota_w, wG)$
implies that a bounded function $f$ is weak Hartman
if there exists a $\mu_w$-Riemann integrable function
$F \in \Riemanntop(wG) = \Riemanntop_{\mu_w}(wG)$, $\mu_w$ denoting the unique translation invariant measure
on $wG$, such that $f = F \circ \iota$.
From Definition \ref{Dweakhartman} it is also obvious that $\H^w(G)\supseteq\H(G)\cup\W(G)$.

\begin{defi}
Let $G$ be an abelian topological group.
By $\H^w_0(G)$ we denote the set of all weak Hartman functions
$f$, such that $|f|$ has zero mean value.\end{defi}
$\H^w_0(G)$ is a closed ideal of $\H(G)$. We will now identify the corresponding quotient space. Given any $\mu_w$-Riemann integrable function $F: wG\to \C$ we can write
\[F=(\underbrace{F-F\circ \rho}_{=:F_0})+\underbrace{F\circ \rho}_{=:F_h}.\]
Note that the retraction $\rho:wG \to wG$ is measure-preserving:
\begin{eqnarray*}
\rho\circ\mu_w (A) &=& \mu_w(\rho^{-1}[A])\\
                   &=& \mu_w(\{x\in wG: e+x\in A\})\\
                   &=& \mu_b(\{x\in wG: x=e+x \mbox { and } e+x \in A)\\
                   &=& \mu_b(A\cap K(wG))=\mu_w(A)
\end{eqnarray*}
Thus both $F_0$ and $F_h$ are $\mu_w$-Riemann integrable.
The induced weak Hartman function $f:=F\circ \iota_w$ can be written
as the sum $f=f_0+f_h$ where $f_0:=F_0\circ\iota_w$ is a weak Hartman
function with zero mean-value and $f_h:=F_h\circ\iota_w$ is an ordinary Hartman function. This decomposition
is unique. So we have proved:

\begin{thr}\label{Tweakhartmandeco}
Let $G$ be an abelian topological group and denote by $\H^w(G)$ and $\H^w_0(G)$
the space of weak Hartman functions resp. the space of
weak Hartman functions with zero mean value. Then
$\H^w(G)= \H^w_0(G)\oplus \H(G)$.
\end{thr}

\section{Compactifications of LCA groups}
\label{Slcacomp}
In the following sections we will deal with \emph{locally compact abelian (LCA)} groups. If $H$ is a subgroup of the topological group $G$, we will denote this by $H\le G$. By $G_d$ we mean the group $G$ equipped with the discrete topology. We will use standard notation such as $\hat{G}$ for the Pontryagin dual, $\chi$ for characters, $H^{\perp}$ for the annihilator of a subgroup and $\phi^*$ for the adjoint of a homomorphism without further ado and refer the reader instead to standard textbooks on this topic such as \cite{Arma81,HeRo79,Stro06}.

Similarly to Proposition \ref{Palgebracontinuous} where we used a function algebra $\A$ to construct a compactification we will now use a group of characters. Let $H\le \hat{G}_{d}$ be an (algebraic) subgroup of the dual of $G$. $H$ induces a group compactification $(\iota_H,K_H)$ of $G$ in the following way:
\[\iota_H: g \mapsto (\chi(g))_{\chi\in H}, \quad K_H:=\overline{\iota_H(G)}\le \T^H\]
and for every such $H\le \hat{G}_d$ the kernel of $\iota_H$ coincides with the annihilator $H^{\perp}\le G$. Remarkably, also the converse is true.

\begin{pro}\label{Plcacomp}
Let $G$ be an LCA group and let $(\iota, C)$ be a group compactification of $G$. Then there exists a unique subgroup $H\le \hat{G}_d$ such that $(\iota, C)$ and $(\iota_H,K_H)$ are equivalent, namely $H=\iota^*(\hat C)$.
\end{pro}
\begin{proof} As $\iota: G\to C$ has dense image, the adjoint homomorphism $\iota^*: \hat C \to \hat G$ is one-one. Let $H:=\iota^*(\hat C)\le \hat G$ and
consider the group compactification $(\iota_H, K_H)$.
Note that for $g\in G$ we have, due to injectivity of $\iota^*$,
$$(\chi(\iota (g)))_{\chi \in \hat C} = (\iota^*(\chi)(g))_{\chi \in \hat C} = (\eta (g))_{\eta \in H} = \iota_H (g).$$
Define $\pi: C \to \T^{\hat C}$ via $c\mapsto (\chi (c))_{\chi\in \hat C}$. Then $\pi$ is a continuous homomorphism and
maps the dense subgroup $\iota(G)\le C$ onto the dense subgroup $\iota_H(G)\le K_H$. As $C$ is compact $\pi (C)$ is closed and thus contains $K_H$.
On the other hand $\pi^{-1}[K_H]$ is closed since $\pi$ is continuous and so contains $C$. Thus $\pi$ maps $C$ onto $K_H$. Since $\iota_H=\pi\circ\iota$ this implies $(\iota, C)\ge (\iota_H, K_H)$.
If $\pi(c)=0$ then $\chi(c)=0$ for all $\chi\in \hat C$. Thus $c=0$ and $\pi$ must be one-one. So $(\iota, C)\cong (\iota_H, K_H)$ via $\pi$, cf.\ Definition \ref{Dcompvergleich}.
\end{proof}

\section{Realizability on LCA Groups}
\label{Srealize}
Let us turn now towards the \emph{realizability} of Hartman functions as Riemann integrable functions, cf.\ Definition \ref{Dcomp}. It follows from Theorem 4 in \cite{Wink02}, that every $f \in \H(\Z)$ which is a characteristic function can be realized in a metrizable compactification. We are going to generalize this result. As a corollary we prove that metric realizability of every $f \in \H(G)$ is possible precisely for LCA groups with separable dual. First we have to establish some useful concepts.

\subsection{Preparation}
\begin{defi}\label{Dtopweight}
Let $G$ be an LCA-group. The \textbf{topological weight} $\kappa(G)$ is defined as the cardinal number
\[ \kappa(G)= \min  \{|I|: (O_i)_{i\in I} \mbox { is an open basis of } G\}.\]
\end{defi}

Note that this minimum exists (and is not merely an infimum) since cardinal numbers are well-ordered. The topological weight behaves well to products, i.e.\ for infinite $G$ we have $\kappa(G\times H)=\max\{\kappa(G),\kappa(H)\}$ and $\kappa(G^I)=\kappa(G)\cdot|I|$. For $H\subseteq G$ we clearly have $\kappa(H)\le \kappa(G)$.

From the theory of LCA-groups it is known that $\kappa(G)=\kappa(\hat{G})$, see \cite[\S 24.14]{HeRo79}. The topological weight of the group compactification $K_H$ can thus be computed very easily via $\kappa(K_H)=\kappa(\hat{K}_H)= \kappa(H_d)=|H|$.

\begin{defi}\label{Dcoweight}
Let $G$ be an LCA-group. The \textbf{co-weight} $c(G)$ is defined as
$\min \{|H|: H \le \hat G \mbox { and } \overline H = \hat G\}$.
\end{defi}

We collect some facts concerning the co-weight:
\begin{enumerate}
\item $c(G)\le \aleph_0$ $\Leftrightarrow$ $G$ has separable dual.
\item $c(G)<\infty$  $\Leftrightarrow$ $G$ is finite.
\item If $H$ is a closed subgroup of $G$, then $c(H)\le c(G)$.
\end{enumerate}
The statements (1) and (2) are obvious. For the sake of completeness we give the argument for (3): By duality $\hat{H}$
and $\hat{G}/H^{\perp}$ are isomorphic LCA groups. Let $G_0\le \hat{G}$ be a subgroup of $\hat{G}$ with $|G_0|=c(G)$. Then $H_0:=G_0+H^{\perp}$ is a subgroup of $\hat{G}/H^{\perp}$ and $|H_0|\le |G_0|$. As the canonical projection $\pi_{H^{\perp}}: \hat{G} \to \hat{G}/H^{\perp}$ is continuous and onto, dense sets are mapped onto dense sets. So pick $G_0\le \hat{G}$ which is dense with $|G_0|=c(G)$ to conclude $c(H)\le|\pi_{H^{\perp}}(G_0)|\le|G_0|= c(G)$.

\begin{defi}\label{Dhartmanweight}
Let $G$ be an LCA-group and $f\in\H(G)$ a Hartman function. The \textbf{weight} of $f$, $\kappa(f)$, is defined by
\[\kappa(f):= \min \{\kappa(K): f \mbox{ can be realized on } (\iota, K)\}.\]\end{defi}

By virtue of Proposition \ref{Plcacomp} we can compute $\kappa(f)$ as the minimum of all $|H|$ such that $f$ can be realized on $(\iota_H, K_H)$.
We want to prove the following:
\begin{thr}\label{Trealize}
Let $G$ be an LCA group, then $$\max\{\kappa(f):f\in \H(G)\}=c(G),$$ i.e.\ every Hartman measurable function on $G$ can be realized on a compactification whose topological weight is at most $c(G)$, and this is best possible.
\end{thr}

\begin{cor}\label{Crealizemetric}
Let $G$ be an LCA group, then the following are equivalent
\begin{enumerate}
\item $\hat G$ is separable,
\item Every $f\in\H(G)$ can be realized on a metrizable compactification
\end{enumerate}
\end{cor}
\begin{proof}
Separability of $\hat G$ is equivalent to $c(G)\le\aleph_0$.
\end{proof}
The rest of this Section is devoted to the proof of Theorem \ref{Trealize}.
\subsection{Estimate from above}
%
\begin{lem}\label{Linjective}
Let $G$ be an LCA group.
Then there exists an injective group compactification of $G$, i.e.\ a group compactification $(\iota,C)$ such that $\iota: G \to C$ is one-one. Furthermore
\[\min\{ \kappa(C): (\iota,C) \mbox{ is injective} \} = c(G).\]
\end{lem}

\begin{proof}
The result follows from the fact that any group compactification is equivalent to some
$(\iota_H, K_H)$ with $H\le \hat G$ and that $\ker \iota_H = H^{\perp}$. Thus $\iota_H$ is injective
if $H^{\perp}=\{0\}$ and this is equivalent to $H$ being dense in $\hat G$. Thus $c(G)\ge |H| = \kappa(K_H)=\kappa(C)$ and equality is obtained if $|H|=c(G)$ and $(\iota,C)=(\iota_H,K_H)$.
\end{proof}

In the following let us call a group compactification $(\iota, C)$ of an LCA group a \textbf{finite dimensional compactification} if $C\le \T^s$ for some $s\in \N$. (If $(\iota, C) \cong (\iota_H, K_H)$ this is equivalent to $|H|\le s$.)

\begin{lem} \label{Lfinitedim}
Let $G$ be an LCA group and $T\subseteq G$ a Hartman set.
For every $\eps>0$ there are Hartman sets $T_{\eps}$ and $T^{\eps}$,
realized on a finite dimensional compactification $(\iota,C)$ such that
$T_{\eps}\subseteq T \subseteq T^{\eps}$,
$m(T^{\eps}\setminus T_{\eps})<\eps$.
\end{lem}
\begin{proof}
We proceed similar to \cite[Theorem 2]{Wink02}.
Let $M\subseteq bG$ be a $\mu_b$-continuity-set
realizing $T$, i.e. $T=\iota_b^{-1}[M]$.
Use the inner regularity of the Haar measure on $bG$ to find a compact
inner approximation $K\subset M^o$ with
$\mu_b(M\setminus K) = \mu_b(M^o\setminus K) < \tfrac{\eps}2$.

Recall that one can construct the Bohr compactification
as $(\iota_{\hat G {}_{d}},K_{\hat G {}_{d}})$, see Proposition \ref{Plcacomp}.
As $bG=K_{\hat G{}_{d}}\tm \prod_{\chi\in\hat{G}} \overline{\chi(G)}$ one can obtain a basis of open sets $(B_i)_{i\in I}$ in $bG$ by restricting the standard basis of the product space to the subspace $bG$.
The sets $B_i$ can be chosen to be finite intersections of sets of the form
\[
D_{\chi_0;a,b}:=\{ (\alpha_{\chi})_{\chi\in\hat{G}} \in bG:
   \alpha_{\chi_0} \in (a,b)\}
\]
where $(a,b)$ denotes an open segment in $\T$ and such that the basis $(B_i)_{i\in I}$ consists of
$\mu_b$-continuity sets.

We can cover $K$ by finitely many sets of the form
$O_j = B_{i_j} \cap M^o,\; j=1,\ldots n$
with $i_j \in I$. Each $O_j$ is a $\mu_b$-continuity set and
induces a Hartman set $T_j = \iota_b^{-1}[O_j]$
on $G$ that may be realized on a finite dimensional
group compactification $(\iota_j, C_j)$, i.e.\ $C_j \le \T^{s_j}$.
Let $(\iota_0,C_0)$ denote the supremum of all $(\iota_j,C_j)$,
$j=1,\ldots,n$. It is easy to check that $C_0 \le \T^{s_0}$
with $s_0 = \sum_{j=1}^n s_j$ and that
$T_{\eps} = \iota_b^{-1} \left(\bigcup_{j=1}^n O_j \right)$ is
a Hartman set which can be realized in $(\iota_0,C_0)$, see Definition \ref{Dprodcomp}.

In a similar way one finds an outer approximation $T^{\eps}$ which
can be realized in some compactification $(\iota^0,C^0)$ with
$C^0 \le \T^{s^0}$. Then we can take the supremum $(\iota,C)$
of $(\iota_0,C_0)$ and $(\iota^0,C^0)$ and $s = s_0+s^0$.
\end{proof}

\begin{lem}\label{Lmetricrealset}
Every Hartman set $T$ on an infinite LCA group $G$
can be realized on a group compactification with topological weight $c(G)$.
\end{lem}

\begin{proof}
We follow the lines of \cite[Theorem 4]{Wink02}. Let $T$ be a
Hartman set and $(T_{1/n})_{n=1}^{\infty}, (T^{1/n})_{n=1}^{\infty}$
sequences of Hartman sets as in Lemma \ref{Lfinitedim},
approximating $T$ from inside
resp.\ outside. Let $(\iota, C)$ be the supremum of all involved
at most countably many finite dimensional compactifications. As $\kappa(\T^s)=\aleph_0$ for every $s\in \N$,
the topological weight of $C$ can not exceed $\aleph_0\cdot\aleph_0=\aleph_0$.
By Lemma \ref{Linjective} we can find an injective group compactification, covering
$(\iota,C)$ and having topological weight $\max\{c(G),\aleph_0\}=c(G)$.
For notational convenience we call this compactification again $(\iota, C)$.

Denote by $M_n$ resp. $M^n$ the
$\mu_C$ continuity-sets in $C$ that realize the Hartman
sets $T_{1/n}$ resp. $T^{1/n}$.
Thus $M_{\infty}:= \bigcup_{n=1}^{\infty} M_n^{\circ}$ is open,
$M^{\infty}:= \bigcap_{n=1}^{\infty} \overline{M^n}$ is closed and
$\iota^{-1}[M_{\infty}]\subseteq T \subseteq \iota^{-1}[M^{\infty}]$.
Let $M:=M_{\infty}\cup \iota(T)$. Since $\iota$ is one-one the preimage
of $M$ under $\iota$ coincides with the given Hartman set $T$.
Furthermore
\[\mu_C(\partial M)\le \mu(M^{\infty}\setminus M_{\infty})= \lim_{n\to\infty}
 \mu (M^n\setminus M_n)=0\]
shows that $M$ is a $\mu_C$-continuity set.
\end{proof}

\begin{cor}\label{Crealsimple}
Let $G$ be an infinite LCA group and $f \in \H(G)$
with $f(G)$ finite. Then $f$ can be realized in
a group compactification $(\iota, C)$ with topological weight $c(G)$
by a simple $\mu_C$-continuity function.
\end{cor}
\begin{proof}
By assumption $f = \sum_{i=1}^n \alpha_i \Eins_{ T_i}$. It is clear
that the $T_i$ can be taken to be Hartman sets. By Lemma \ref{Lmetricrealset}
$T_i$, $i=1,\ldots,n$ can be realized on a compactification
$(\iota_i,C_i)$ with $\kappa(C_i)=c(G)$. The supremum $(\iota,C)$ of the $(\iota_i,C_i)$
has again topological weight $c(G)$  and, as a consequence of Proposition \ref{Pcomprepr},
each $T_i$ and hence $f$ can be realized in $(\iota,C)$.
\end{proof}

Now we can prove the first part of Theorem \ref{Trealize}, namely $\kappa(f)\le c(G)$ for every
$f\in \H(G)$:

\begin{proof}
First consider the finite (compact) case: If $G$ is compact then
every $f\in \H(G)$ can be realized on $(\mbox{id}_G,G)$. Hence
$\kappa(f)\le \kappa(G)=\kappa(\hat{G})=|\hat{G}|=c(G)$.

Now we show that on an infinite LCA group $G$ every $f\in \H(G)$ can be realized on a group compactification with topological weight not exceeding $c(G)$.
W.l.o.g.\ we may assume that $f$ is real-valued. By Proposition \ref{PHartman} $f$ can be realized
in the maximal compactification $(\iota_b,bG)$ by some $F^b \in \Riemanntop_{\mu_b}(bG)$,
i.e.\ $f = F^b \circ \iota_b$. By Lemma \ref{Priemannchar} there is a
sequence of simple $\mu_b$-continuity functions $F_n^b$ on $bG$ converging
to $F^b$ uniformly. Consider the Hartman functions
$f_n = F_n^b \circ \iota_b$ and
note that each $f_n$ takes only finitely many values.
Corollary \ref{Crealsimple} guarantees
that each $f_n$ can be realized on a group compactification
$(\iota_n,C_n)$ with $\kappa(C_n)=c(G)$ by simple $\mu_{C_n}$-continuity functions $F_n^0$, i.e.\
$f_n = F_n^0 \circ \iota_n$. The supremum of countably many
group compactifications of topological weight $c(G)$ has a topological weight not exceeding
$\aleph_0\cdot c(G)=c(G)$. By technical convenience we use
Lemma \ref{Linjective} to get an injective group
compactification $(\iota,C)$ with $\kappa(C)=c(G)$ covering all $(\iota_n,C_n)$. For each $n$
let $\pi_{n}: C \to C_n$ denote the canonical projection, i.e.
$\iota_n = \pi_{n} \circ \iota$. Consider the functions
$F_n = F_n^0 \circ \pi_{n}$ which are in $\Riemanntop_{\mu_C}(C)$ by Proposition \ref{Pcomprepr}
and in fact simple $\mu_C$-continuity functions.
\begin{diagram}[textflow]
        &               &        C   &            &\\
        &\ruTo^{\iota}  &\dTo>{\pi_n}&\rdTo^{F_n} &\\
G       & \rTo^{\iota_n}&        C_n &\rTo^{F_n^0}&\C.
\end{diagram}

In order to realize $f$ in
$(\iota,C)$ by $F \in \Riemanntop_{\mu_C}(C)$ we have to define $F(x)=f(g)$ whenever
$x = \iota(g)$ for some $g \in G$. Since $\iota$ is one-one
$F$ is well-defined on $\iota(G)$. For $x \in C \setminus \iota(G)$ we define
\[
F(x) = \limsup_{\iota(g) \to x} F(\iota(g)).
\]
It remains to show that $F$ is
$\mu_C$-Riemann integrable. For each $n \in \N$ let
$F_n = \sum_{i=1}^{k_n} \alpha_i \Eins_{A_{n,i}}$ be a representation
of $F_n$ with pairwise disjoint continuity sets $A_{n,i}$, $i=1,\ldots, k_n$.
The open sets $U_n = \bigcup_{i=1}^{k_n} A_{n,i}^o$ have full $\mu_C$-measure.
Thus the dense $G_{\delta}$-set $U = \bigcap_{n \in \N} U_n$ has
full $\mu_C$-measure as well. If we can prove the following
claim, we are done.

\begin{claim} Each $x \in U$ is a point of continuity for $F$.
\end{claim}

Fix $x \in U$ and $\eps > 0$. We are looking for an open neighborhood
$V\in\mathfrak{U}(x)$ ($\mathfrak{U}(x)$ denoting the filter of neighborhoods of $x$) such that $\iota(g_1), \iota(g_2) \in V$ implies
$|f(g_1)-f(g_2)| < \eps$. This suffices to guarantee
$|F(x_1) - F(x_2)| \le \eps$ for all $x_1,x_2 \in V$, in particular
$|F(x_1)-F(x)| \le \eps$, yielding continuity of $F$ in $x$.
To find such a $V$ note that, by construction, the $F_n$ converge
uniformly to $F$ on the dense set
$\iota(G)$. Choose $n \in \N$ in such a way that
$|F_n(\iota(g)) - F(\iota(g))| < \tfrac {\eps}2$ for all $g \in G$.
There is a unique $i \in \{1, \ldots,k_n\}$ such that $x \in A_{n,i}^{\circ}$.
The set $V := A_{n,i}^{\circ}$ has the desired property:
For $\iota(g_1), \iota(g_2) \in V$ we have $F_n(\iota(g_1)) = F_n(\iota(g_2))$
and $$|f(g_1) - f(g_2)| = |F(\iota(g_1)) - F_n(\iota(g_1))| +
|F_n(\iota(g_2)) - F(\iota(g_2))| < \eps.\qedhere$$
\end{proof}

\subsection{Estimate from below}
In this section we are concerned with the construction of  a Hartman function with $\kappa(f)=c(G)$ for a given infinite group $G$.

\begin{lem}\label{Linjectiveset}
Let $G$ be an uncountable LCA-group.
Then there exists a subset $A\subseteq G$ such that $\Eins_A$ can be realized only in injective group compactifications of $G$.
\end{lem}
\begin{proof}
Let $\alpha\to g_{\alpha}$ be a bijection between the set of all ordinals
$\alpha< |G|$ and  $G\setminus\{0\}$.
By choosing for each $\alpha$ elements of the co-sets of the subgroup $\langle g_{\alpha} \rangle$ generated by $g_{\alpha}$ we find (for each $\alpha$) $x_{\alpha}^{(i)}\in G$ with $|I|\le |G|$ such that
$$ G= \bigcup_{i\in I}´\left(\,\langle g_{\alpha}\rangle+ x_{\alpha}^{(i)}\right), \quad \alpha < |G|,$$ is a disjoint union. As $G$ is uncountable and $\langle g_{\alpha} \rangle$ is at most countable we must have $|I|=|G|$. We start with the construction of $A$: Assume, by transfinite induction, that for given $\alpha_0<|G|$ we have already constructed elements $x_{\alpha},y_{\alpha}, \alpha<\alpha_0$  such that $y_{\alpha}\in \langle g_{\alpha} \rangle + x_{\alpha}$ with $x_{\alpha}=x^{(i)}_{\alpha}$ for some $i=i(\alpha)\in I$ and such that all the $x_{\alpha},y_{\alpha}$ with $\alpha<\alpha_0$ are pairwise distinct.
To find $y_{\alpha_0}$ we first observe that
$$N_1:=\{y_{\alpha},x_{\alpha}:\alpha<\alpha_0\}\subseteq \bigcup_{\alpha<\alpha_0} \langle g_{\alpha} \rangle + x_{\alpha},$$
hence we have $|N_1| \le |\alpha_0|\cdot\aleph_0<|G|$. Therefore, by the cardinality of
\[G=\bigcup_{i\in I}\left( \langle g_{\alpha_0} \rangle + x_{\alpha_0}^{(i)}\right),\]
there are $|G|$ many indices $i\in I$ with $\langle g_{\alpha_0} \rangle + x_{\alpha_0}^{(i)}$ disjoint from $\{y_{\alpha},x_{\alpha}: \alpha<\alpha_0\}$.
Pick such an $i=i(\alpha_0)$ and $x_{\alpha_0}=x_{\alpha_0}^{(i)}, y_{\alpha_0}=y_{\alpha_0}^{(i)} \in \langle g_{\alpha_0} \rangle + x_{\alpha_0}^{(i)}$ with $x_{\alpha_0}\neq y_{\alpha_0}$. Let $A=\{x_{\alpha}:\alpha<|G|\}$.

Note that by its very definition $y_{\alpha}\notin A$ for all $\alpha<|G|$. Suppose $A$ can be realized in some non-injective group compactification $(\iota, C)$, i.e. there exists a set $M\subseteq C$ such that $\iota^{-1}[M]=A$. As $\iota$ is not injective, there exists some $g_{\alpha}\in \ker \iota$, $g_{\alpha}\neq 0$.
Then $\langle g_{\alpha}  \rangle \subseteq \ker \iota$. As
$y_{\alpha}\in \langle g_{\alpha} \rangle + x_{\alpha}$ we have $\iota (y_{\alpha})=\iota(x_{\alpha})\in M$ , i.e.\ $y_{\alpha}\in \iota^{-1}[M]=A$. Contradiction.
\end{proof}

\begin{rem} For countable $G$ we could use a similar and even simpler construction. However, we will use a different approach for this case.\end{rem}

\begin{lem}\label{Lhmsubgroup}
Let $G$ be an infinite LCA group. Then there exists a closed subgroup $G_0\le G$ such that $c(G_0)=c(G)$ and $G_0$ is a Hartman null-set.
\end{lem}
\begin{proof}
First we distinguish two cases and employ in each of them the existence of a nontrivial closed subgroup with zero Hartman measure.
\begin{enumerate}\item
There exists a $\chi\in\hat G$ such that $\chi(G)$ is infinite. The annihilator $G_0:=\{\chi\}^{\perp}$ is a Hartman null-set. To see this consider the preimage of the singleton $\{0\}$ in the group compactification $(\chi, \T)$.
\item
All characters are torsion elements, i.e.\ $\chi(G)$ is finite for every $\chi\in\hat{G}$. Pick a sequence of pairwise (algebraically) independent characters $\chi_n$ and consider $G_0:=\bigcap_{i=1}^{\infty}\{\chi_i\}^{\perp}$ and the group compactification
$(\iota,C):=\bigvee_{i=1}^{\infty}(\chi_i, \Z_{m_i})$ with $m_i=|\chi_i(G)|$. As $C$ is infinite, every singleton in $C$ has zero $\mu_C$-measure. Therefore $G_0=\iota^{-1}[\{0\}]$ is a Hartman null-set.
\end{enumerate}

As $G_0\le G$ is a closed subgroup $c(G_0)\le c(G)$.
Let $H$ be a dense subgroup of $\hat{G_0}$ with $|H|=c(G_0)$. Any $\eta\in H$ can be extended from a character on $G_0$ to a character $\tilde\eta$ on the whole of $G$. Let us denote the set of these extended characters by $\tilde H$. Note that as $H$ is infinite
the subgroup
\[H_0:=\left\{ \begin{array}{cl}
\tilde H & \mbox{ for } $G$ \mbox{ as in case } (1)\\
\langle\tilde H\cup\{\chi_i: i\in \N\}\rangle & \mbox{ for } $G$ \mbox{ as in case } (2)
\end{array}\right. \]
has the same cardinality as $H$. We show that $H_0$ is dense in $\hat{G}$ by computing its annihilator. Pick any $g\in G$ such that $g\in H_0^{\perp}$. If $g\in G_0$ we have in particular $\chi_{|G_0}(g)=0$ for every $\chi \in \tilde{H}$ and thus $\eta(g)=0$ for every $\eta\in H$. As $H$ is dense in $\hat{G_0}$ we have $g=0$. Thus $H_0^{\perp}\cap G_0=\{0\}$. On the other hand, by its very definition, $H_0^{\perp}\le G_0$. So $H_0^{\perp}=\{0\}$ and thus $\overline{H_0}=G$. This gives the reverse inequality $c(G)\le |H_0|=|H|=c(G_0)$.
\end{proof}

By proving the next statement we conclude the remaining part of Theorem \ref{Trealize}.

\begin{pro}\label{Pexistshmfunction}
Let $G$ be an LCA group. Then there exists a Hartman function $f\in\H(G)$ such that $\kappa(f)= c(G)$.
\end{pro}
\begin{proof}
\hfill\par\vspace{-2ex}
\begin{enumerate}
\item Finite $G$: Let $\hat{G}=\{\chi_1,\ldots \chi_n\}$ and take $f=\chi_1+\ldots+\chi_n$. Then $\kappa(f) = n=c(G)$.
\item
Countable $G$:  $|G|=\aleph_0$ implies $c(G)=\aleph_0$. Note that every countable LCA group is discrete.
Thus $\hat{G}$ is isomorphic to a closed subgroup of $\T^{\aleph_0}$. In particular $\hat G$ is uncountable.

Consider any sequence $(\chi_i)_{i=1}^{\infty}$ of pairwise (algebraically) independent characters on $G$ and define the almost periodic function $$f=\sum_{i=1}^{\infty} \frac{\chi_i}{2^n}$$ on $G$. Denoting by $m_G$ the unique invariant mean on $AP(G)$, we see that
\[\Gamma(f):=\langle \{ \chi \in \hat G: m_G(f\cdot\chi)\neq 0\}\rangle = \langle\{\chi_i: i \in \N\}\rangle.\]
Then, by \cite[Corollary 13]{Mare05}, we have $\Gamma(f)\subseteq H$ for any compactification $(\iota_H, K_H)$ on which $f$ can be realized. As $f$ can only be realized on compactifications $(\iota_H,K_H)$ with infinite $H$ we conclude $\kappa(f)\ge |H| = \aleph_0=c(G)$.

\item Uncountable $G$: Take $G_0\leq G$ as in Lemma \ref{Lhmsubgroup} and let $f:=\Eins_A$ for the set $A\subseteq G_0$ from Lemma \ref{Linjectiveset}. As $G_0$ has zero Hartman measure Corollary \ref{Chartmancomplete} implies $f\in\H(G)$. By construction for every compactification $(\iota, C)$ where $f$ can be realized $\iota$ must be injective on $G_0$, hence $\kappa(f)\ge c(G_0)=c(G)$.\qedhere
\end{enumerate}
\end{proof} 
\chapter{Classes of Hartman functions}\label{CHwap}
\section{Generalized jump discontinuities}\label{Sgjd}
The concept of generalized jump discontinuities is useful for
comparing Hartman functions and weakly almost periodic functions. In the present section we do not need the group setting.

\begin{defi}\label{dgjd}
Let $X,Y$ be topological spaces.
A function $f: X\to Y$
has a \textbf{generalized jump discontinuity (g.j.d.)} at $x\in
X$ if there are (disjoint) open sets $O_1$ and $O_2$, such that
$x\in \overline{O_1}\cap\overline{O_2}$ but
$\overline{f(O_1)}\cap\overline{f(O_2)}=\emptyset$.\end{defi}

\begin{exa}\label{egjd}
\hfill\par\vspace{-3ex}
\begin{enumerate}
\item The function $f_1(x)=\Eins_{\left[0,\tfrac{1}{2}\right)}(x)$ on $X=[0,1]$
has a g.j.d.\ at $x=\tfrac{1}{2}$. The function
$f_2(x)=\Eins_{ \left\{\tfrac{1}{2}\right\} }(x)$ has no g.j.d.\
\item The function $f: [0,1]\to \R$
$$f(x)=\left\{\begin{array}{cl} \sin\left(\tfrac{1}{x}\right)&x\neq 0\\
0&x=0\end{array}\right.$$ has a g.j.d.\ at $0$. To see this, consider the open set
$O_1:=f^{-1}[ (\tfrac 12,1]]$ and the open set $O_2:=f^{-1}[[-1,-\tfrac 12)]$.
\item Generalizing the first example, let $X$ be compact and $\mu$ be a finite complete regular Borel measure on $X$ and $A$ a $\mu$-continuity set.  Then the characteristic function $\Eins_A$ has g.j.d.s on $\partial(A^{\circ})\cap\partial(A^{\circ c})$, the common boundary of $A^{\circ}$ and its complement.
\item Let $X$ be compact and $\mu$ a finite complete regular Borel measure with $\mbox{supp}(\mu)=X$. If $f:X\to \R$ is constant $\mu$-a.e.\ then $f$ has no
    g.j.d.
\end{enumerate}
\end{exa}

\begin{pro} \label{Pgjdclosed}
Let $X$ be a topological space. Denote by $J(X)$ all bounded functions
$f: X\to \R$ having a g.j.d.\ Then $J(X)\tm B(X)$ is open in the topology of uniform convergence. \end{pro}
\begin{proof}
Let $f\in J(X)$. Then there exist disjoint open sets $O_1$ and $O_2$
with $\partial O_1 \cap \partial O_2\neq\emptyset$ but
$\overline{f(O_1)} \cap \overline{f(O_2)}=\emptyset$. Let $\eps:=
d(f(O_1),f(O_2))>0$ and suppose $\|f-g\|_{\infty}<\tfrac{\eps}{8}$.
Then $d(f(Y),g(Y))<\tfrac{\eps}{4}$ for any set $Y\subseteq X$, hence
\[|d(g(O_1),g(O_2))-d(f(O_1),f(O_2))|< \tfrac{\eps}{2}.\] In particular
$d(g(O_1),g(O_2))>0$, i.e.\ $g$ has a g.j.d.
\end{proof}

\begin{lem} \label{Lriemann}
Let $X$ be compact and $\mu$ a finite complete regular Borel measure with $\mbox{supp}(\mu)=X$.
Let $f,g \in \Riemanntop_{\mu}(X)$ be Riemann integrable functions.
\begin{enumerate}
\item
If $f$ and $g$ coincide on a dense set,
then they coincide on a co-meager set of full $\mu$-measure.
\item
If $f$ and $g$ coincide on a dense set and $f$ has a g.j.d.\ at $x\in X$,
then also $g$ has a g.j.d.\ at $x\in X$.
\end{enumerate}
\end{lem}

\begin{proof}
\hfill\par\vspace{-2ex}
\begin{enumerate}
\item
By Proposition \ref{Pfsigma} it suffices to show
$[f = g]:=\{x\in X: f(x)=g(x)\} \supseteq X\setminus\left(\mbox{disc}(f)\cup\mbox{disc}(g)\right)$.
Let $x\in X$ be a point of continuity both for $f$ and $g$,
$U\subseteq X$ a neighborhood of $x$ such that $y\in U$ implies
$|f(y)-f(x)|<\tfrac{\eps}{2}$ and $|g(y)-g(x)|<\tfrac{\eps}{2}$.
As $[f=g]$ is dense in $X$ we can pick $y_{\eps} \in U\cap [f=g]$. Thus
\[|f(x)-g(x)|\le |f(x)-f(y_{\eps})|+|g(y_{\eps})-g(x)|<\eps.\] Since $\eps>0$
was arbitrary this implies $f(x)=g(x)$.
\item
Choose $O_1$ and $O_2$ according to the definition of a
g.j.d.\ of $f$ at $x$. By 1. $f$ and $g$ coincide
on a dense set of common continuity points.
Thus for every $x\in X$ and $U\in \U(x)$ we can pick
$x_i^{U}\in U\cap O_i$, $i=1,2$, which are both points of continuity
for $f$ and $g$ and such that $f(x_i^{U})=g(x_i^{U})$, $i=1,2$.
Pick open neighborhoods $O_i^{U}$ of $x_i^U$ such that
$O_i^{U}\subseteq U\cap O_i$ and
\[\mbox{diam}(g(O_i^U))<\tfrac{1}{3}\mbox{ dist}(f(O_1),f(O_2)),\quad i=1,2.\]
Consider the open sets $U_i:=\bigcup_{U \in \mathfrak{U}(x)} O_i^U$, $i=1,2$.
Then $g(U_1)$ is separated from $g(U_2)$ and $x_i^U\in U_i$ for
all $U \in \mathfrak{U}(x)$. This implies $x\in \overline{U_i}$, $i=1,2$, proving that $x$ is a g.j.d.\ for $g$. \qedhere
\end{enumerate}
\end{proof}

\begin{cor} \label{Criemann}
Let $X$ be compact and $\mu$ a finite complete regular Borel measure, $\mbox{supp}(\mu)=X$. Let $f,g$ be simple $\mathfrak{C}_{\mu}$-functions (see Definition \ref{DJordan}). If $f$ and $g$ coincide on a dense set, then they coincide on an open set of full $\mu$-measure.
\end{cor} \begin{proof}
Lemma \ref{Lriemann} implies that $[f=g]$ has full
$\mu$-measure. Since $f$ and $g$ are simple $\mathfrak{C}_{\mu}$-functions $[f=g]$ is a $\mu$-continuity set. Thus
$[f=g]^{\circ}$ and $[f=g]$ have the same
$\mu$-measure $\mu(X)$.\end{proof}

\begin{pro}\label{Pregularization}
Let $X$ be compact and $\mu$ a finite complete regular Borel measure with $\mbox{supp}(\mu)=X$. Let $f\in \Riemanntop_{\mu}(X)\setminus J(X)$, i.e.\ $f$ is Riemann integrable without a g.j.d. Then there exists a \emph{unique} continuous function $f_r\in C(X)$, the \emph{regularization} of $f$, such that $f$ and $f_r$ coincide on $X\setminus\mbox{disc}(f)$.
\end{pro}

\begin{proof}
For Riemann-integrable $f$ the set $X \setminus \mbox{disc}(f)$ is dense in $X$ by Proposition \ref{Pfsigma}. Hence
there is at most one continuous $f_r$ with $f_r(x)=f(x)$ for $x\notin \mbox{disc}(f)$ .

Let $x\in \mbox{disc}(f)$. For each $U \in \mathfrak{U}(x)$, the neighborhood system $\mathfrak{U}$ of $x$,
$y \in U \setminus \mbox{disc}(f)$ and $\eps>0$
pick an open neighborhood $O = O(U,\eps,y) \in \mathfrak{U}(y)$ such that
$O \subseteq U$ and $\mbox{diam}(\overline{f(O)}) < \eps$. Let
$$O(U,\eps):= \bigcup_{y \in U \setminus \mbox{\tiny disc}(f)} O(U,\eps,y).$$
\begin{claim}The set $\Lambda(x)$ consists of exactly one point $\lambda_x$, where
$$\Lambda(x):=\bigcap_{\eps>0}\bigcap_{U\in\U(x)} \overline{f(O(U,\eps))}.$$
\end{claim} $\Lambda(x)\neq\emptyset$ by the finite-intersection
property of the compact sets $\overline{f(O(U,\eps))}$,
$U \in \mathfrak{U}(x)$, $\eps > 0$.
Suppose by contradiction that $\lambda_1,\lambda_2 \in \Lambda(x)$ and
$\lambda_1\neq\lambda_2$. Consider the open sets
$O_i := \bigcup O(U,\eps,y)$, $i=1,2$, where the union is taken
over all triples $(U,\eps,y)$ with $U \in \mathfrak{U}(x)$, $\eps < \tfrac{|\lambda_1 - \lambda_2|}4$ and
$y \in U \setminus \mbox{disc}(f)$ such that
$|f(y) - \lambda_i|<\eps $.
By construction we have
$x\in \overline{O_1}\cap\overline{O_2}$ and $\overline{f(O_2)}\cap
\overline{f(O_2)}=\emptyset$. Hence $x$ is a g.j.d. of $f$. Contradiction.

\begin{claim}
$f_r: X\to \R$,  $f_r(x)=\left\{\begin{array}{cl}
f(x)&\mbox{ for } x\notin\mbox{disc}(f)\\
\lambda_x& \mbox{ for } x\in\mbox{disc}(f)
\end{array}\right.$ is continuous.\end{claim}
It is immediate to check that
$\mbox{disc}(f_r) \subseteq \mbox{disc}(f)$.
Suppose by contradiction that there exists $x \in \mbox{disc}(f_r)\tm\mbox{disc}(f)$.
Then an inspection of the argument above shows that $x$ would be a g.j.d.\ for $f$.
\end{proof}

Note that for $f: X \to \R$ meeting the requirements of Proposition \ref{Pregularization} we have
\begin{align*}
\|f_r\|_{\infty}&=\sup_{x\in X} |f_r(x)|=\!\sup_{x\in X\setminus\mbox{\tiny disc}(f)}\! |f_r(x)|=\!\sup_{x\in X\setminus\mbox{\tiny disc}(f)}\!|f(x)|\\
                &\le \sup_{x\in X} |f(x)|=\|f\|_{\infty}.
\end{align*}
Thus the mapping $f\mapsto f_r$ is continuous w.r.t.\ the topology of uniform convergence on its domain of definition, i.e.\ on
$\Riemanntop_{\mu}(X)\setminus J(X)\subseteq B(X)$.

\begin{cor}\label{Cgjdchar}
Let $X$ be compact and $\mu$ a finite complete regular Borel measure, $\mbox{supp}(\mu)=X$. For $f\in \Riemanntop_{\mu}(X)$ the following assertions are equivalent:
\begin{enumerate}
\item There exists $g\in C(X)$ such that $f$ and $g$ coincide on a
co-meager set of full $\mu$-measure.
\item $f$ has no g.j.d.
\end{enumerate}
\end{cor}
\begin{proof}
$(1)\Rightarrow (2)$: Suppose $f$ has a g.j.d. at $x\in X$. Pick open sets
$O_1$ and $O_2$ according to the definition of a g.j.d.\ at $x$. Next, pick nets
$\{x_{\nu}^{(i)}\}_{\nu \in \mathcal{N}_i}$, where $(\mathcal{N}_i,\le)$ are directed sets, such that
\[x_{\nu}^{(i)} \in O_i \cap [f=g], \quad \lim_{\nu \in \mathcal{N}_i}=x,\quad i=1,2.\]
This gives the desired contradiction
\begin{align*}
\overline{f(O_1)}&\ni
\lim_{\nu \in \mathcal{N}_1}f(x_{\nu}^{(1)})=
  \lim_{\nu \in \mathcal{N}_1}g(x_{\nu}^{(1)})\\
  &=\lim_{\nu \in \mathcal{N}_2}g(x_{\nu}^{(2)})=
  \lim_{\nu \in \mathcal{N}_2}f(x_{\nu}^{(2)})
\in\overline{f(O_2)}.\end{align*}
$(2)\Rightarrow (1)$: The statement follows from Proposition \ref{Pregularization} and
Lemma \ref{Lriemann}. \end{proof}

\section{Hartman functions that are weakly almost periodic}\label{Shartmanwap}
Recall the notion of weak almost periodicity from Sections \ref{Smeans} and
\ref{Sweakhartman}.

\begin{thr} \label{Tjump}
Let $G$ be a topological group and
$f\in\mathcal{H}(G) \cap \W(G)$ a weakly almost periodic
Hartman function. Let $(\iota, C)$ be a group compactification on which
$f$ can be realized by $F \in \Riemanntop_{\mu_C}(C)$. Then $F: C\to \C$ has no g.j.d.
\end{thr}

\begin{proof}
Assume, by contradiction, that $f\in \H(G) \cap \W(G)$ can be realized on the group compactification
$(\iota,C)$ by $F\in\Riemanntop_{\mu_C}(C)$, where $F$ has a g.j.d.\ at $x_0\in C$. Pick $O_1,O_2\subseteq C$ as in the definition of a g.j.d.\ at $x_0$.

Pick a net $(g_{\nu})_{\nu \in \mathcal{N}}$
in $G$, where$(\mathcal{N},\le)$ is a directed set,
in such a way that $\iota(g_{\nu})\in O_1$
and $\lim_{\nu\in \mathcal{N}} \iota(g_{\nu})=x_0$. W.l.o.g.\ we take $(g_{\nu})_{\nu \in \mathcal{N}}$ to be a universal net, i.e.\ for every $A\subseteq G$  $(g_{\nu})_{\nu \in \mathcal{N}}$ stays eventually in $A$ or $G\setminus A$.
Furthermore we define
\[\phi_{\mathcal{N}}: \tilde{f} \mapsto  \lim_{\nu\in \mathcal{N}}
   \tilde{f}(g_{\nu})\quad(=\lim_{\nu\in \mathcal{N}}
   \delta_{g_{\nu}}(\tilde{f})),\]
where $\delta_{g_{\nu}}$ denotes the evaluation functional at the point $g_{\nu}$. By universality of $(g_{\nu})_{\nu \in \mathcal{N}}$, $\phi_{\mathcal{N}}$ is well-defined and a bounded linear functional on $B(G)$.
Since $x_0 \in \overline{O_2}$ for every neighborhood $V\subseteq C$ of $e$, the neutral element of the group $C$, we can find a neighborhood $U\subseteq C$ of $x_0\in C$ and a $g=g_{U,V} \in G$ such that
$\iota(g_{U,V}) \in V$ and $\iota(g_{\mu})\cdot U \subseteq O_2$.
All such pairs $(U,V)$ form a directed set $\mathcal{M}'$ equipped with the order
$(U_1,V_1)\le(U_2,V_2):\Leftrightarrow U_1\supseteq U_2$ and $V_1 \supseteq V_2$.
The net  $(g_{\mu'})_{\mu \in \mathcal{M}'}$ has the property that for every $\mu'\in\mathcal{M}'$ the net $(\iota(g_{\nu} g_{\mu'}))_{\nu\in\mathcal{N}}$ stays eventually in $O_2$.

Pick a directed subset $(\mathcal{M},\preceq)$ such that
$(g_{\mu})_{\mu\in\mathcal{M}}$ is a universal refinement of
net $(g_{\mu'})_{\mu\in\mathcal{M}'}$. Then $\lim_{\mu\in \mathcal{M}} \iota(g_{\mu})=e \in C$. As $f$ is weakly almost periodic, the closure of the left-translation orbit $$O_L(f)=(L_g f: g\in G)$$ is weakly compact in $B(G)$. This implies that there exists $f_0 \in B(G)$ in the weak closure of $O_L(f)$ such that
weak-$\lim_{\mu \in \mathcal{M}} L_{g_{\mu}} f= f_0$. Consider the evaluation functionals $\delta_{g_{\nu}}\in B(G)^*$:
\[\delta_{g_{\nu}}(f_0)=
f_0(g_{\nu})=\lim_{\mu\in \mathcal{M}} L_{g_{\mu}}f(g_{\nu})
= \lim_{\mu\in \mathcal{M}} f(g_{\nu}g_{\mu})=\lim_{\mu\in \mathcal{M}}
F(\iota(g_{\nu}g_{\mu})).\]
As $\lim_{\mu\in\mathcal{M}} \iota(g_{\mu})=e\in C$ for fixed $\nu\in \mathcal{N}$ the net $(\iota (g_{\nu} g_{\mu}))_{\mu\in\mathcal{M}}$ stays in $O_1$ eventually. Hence $\delta_{g_{\nu}}(f_0) \in \overline{F(O_1)}$ and
thus $\phi_{\mathcal{N}} (f_0)= \lim_{\nu\in \mathcal{N}} f_0(g_{\nu}) \in \overline{F(O_1)}$. Let us now compute the value of the functional $\phi_{\mathcal{N}}$
at $f_0$ directly:
\begin{align*}
\phi_{\mathcal{N}}(f^0)&= \lim_{\mu\in \mathcal{M}} \phi_{\mathcal{N}} (L_{g_{\mu}}f)\;\;\,= \lim_{\mu\in \mathcal{M}}\lim_{\nu\in {\mathcal{N}}}L_{g_{\mu}}f(g_{\nu})&\\&= \lim_{\mu\in \mathcal{M}}\lim_{\nu\in {\mathcal{N}}} f(g_{\nu} g_{\mu})= \lim_{\mu\in \mathcal{M}}\lim_{\nu\in {\mathcal{N}}}
   F(\iota(g_{\nu} g_{\mu})).&
\end{align*}
Thus $\phi_{\mathcal{N}}\in \overline{F(O_1)}\cap\overline{F(O_2)}$. Contradiction.
\end{proof}

\begin{rem} The proof of Theorem \ref{Tjump} employs the same argument (but regarding nets instead of sequences), that may be used to establish the easy direction of Grothendieck's Double Limit Theorem \ref{Tdoublelimit}. \end{rem}

\begin{cor}\label{Chnichtw}
Let $G$ be an infinite LCA group.
Then there exists a Hartman function which is not weakly almost periodic. In particular $\mathcal{H}(G)\neq\mathcal{W}(G)$.
\end{cor}
\begin{proof}
Let $(\iota,C)$ be any infinite metrizable compactification of $G$.
(This can be obtained by taking pairwise distinct characters
$\chi_n$, $n \in \N$, and $\iota: g \mapsto (\chi_n(g))_{n \in \N}$,
$C:= \overline{\iota(G)} \le \T^{\N}$.)
It suffices to find two
disjoint open $\mu_C$-continuity sets $O_1$ and $O_2$ in $C$ with
a common boundary point $x\in \partial O_1\cap \partial O_2$.
Then $x$ is a g.j.d.\ for $F=\Eins_{O_1}$ and, by Theorem \ref{Tjump},
$f = F \circ \iota \in\mathcal{H}(G)\setminus \mathcal{W}(G)$.

Let $d: C\times C \to [0,1)$ be a bounded metric which generates the topology of $C$. We employ the fact that for every $x\in C$ there are open balls $B(r,x):=\{y\in C: d(x,y)< r\}$
with center $x$ and arbitrarily small radius $r>0$ which are
$\mu_C$-continuity sets, see \cite[Example 1.3]{KuNi74}, or an argument similar to the proof of our Proposition \ref{Lfinitedim}.

Construction of $O_1, O_2$: Pick any $x\in C$.
We define two sequences of disjoint open $\mu_C$-continuity sets
$\left(O_j^{(1)}\right)_{j=0}^{\infty}$ and
$\left(O_j^{(2)}\right)_{j=0}^{\infty}$. Let $O_0^{(1)}$ and $O_0^{(2)}$ be any two disjoint open balls which are $\mu_C$-continuity sets, separated from $x$ and have $\mu_C$ measure smaller than $\tfrac 12$. We proceed by induction:
Suppose we have already defined
$O_0^{(1)},\ldots,O_n^{(1)}$ and $O_0^{(2)},\ldots,O_n^{(2)}$
such that
\[\mu\left(\bigcup_{j=0}^n O_j^{(i)}\right)<\tfrac{1}{2}(1-\tfrac{1}{2^n})
\quad\mbox{and}\quad
0<\mbox{dist}\left(\bigcup_{j=0}^n O_j^{(i)},x\right)<\tfrac{1}{2^n}, \quad i=1,2.\]
Let
\[r<\min \left\{\mbox{dist}\left(O_j^{(i)},\right): j=0,\ldots,n \mbox{ and } i=1,2\right\}\]
and pick distinct $x_1,x_2\in B(r,x)$ and $\rho<\min\left\{\tfrac r2,\tfrac{1}{2^{n+1}}\right\}$
such that $O_{n+1}^{(i)}:=B(\rho,x_i)$, $i=1,2$,
are $\mu_C$-continuity sets of $\mu_C$-measure less than $\tfrac{1}{2^{n+1}}$.
Choosing $O_i:=\bigcup_{j=0}^{\infty} O_j^{(i)}$, $i=1,2$,
we obtain two disjoint open sets $O_1,O_2\subseteq C$ with the required properties.
\end{proof}


The converse problem, namely to find weakly almost periodic
functions that are not Hartman measurable appears to be harder. We
content ourselves with the special case $G=\mathbb{Z}$. The key ingredient
for our example are ergodic sequences. These sequences
were extensively studied by Rosenblatt and Wierdl in their paper \cite{RoWi95}. Also in the context of Hartman measurability ergodic sequences
were already mentioned in \cite{ScSW00}.

\begin{exa}[Ergodic sequences]\label{Ewnichth}
A sequence $n_k$ of nonnegative integers is called ergodic if for every measure preserving system $(X,T,\mu)$ with ergodic transformation $T: X\to X$ and every $\mu$-integrable $f$ \[
\lim_{N\to \infty} \frac 1N \sum_{k=0}^{N-1} f\circ T^{n_k}(x) = \int_{X}f d\mu
\] for $\mu$-almost every $x\in X$ holds true. Birkhoff's ergodic theorem (cf. \cite{Walt82}) states that $n_k=k$ is an ergodic sequence.
It is known (\cite[Theorem 11]{ScSW00} and the examples therein)
that there are other ergodic sequences, such as $(k \log k)_{k\in \mathbb{N}}$, which can not be Hartman measurable. On the other hand 0-1 sequences with the property that the length between consecutive 1s tends to infinity while
the length of consecutive 0s stays bounded are weakly almost periodic,
\cite[Theorem 4.2]{BeJM89}. Thus
$\mathcal{E}\subseteq \mathcal{W}(\mathbb{Z})\setminus \mathcal{H}(\mathbb{Z})$, where $\mathcal{E}$ is the set of all ergodic sequences on $\Z$.
\end{exa}

\begin{problem}
Construct $f\in \W(G)\setminus \H(G)$ on more general LCA, or even arbitrary topological groups.
\end{problem} 
\section{Hartman functions without generalized jumps}\label{Shartmannogjd}

Theorem \ref{Tjump} motivates us to further investigate Hartman functions having no g.j.d. First we show that the property of having a g.j.d.\ does not depend on the particular compactification.

\begin{pro} \label{Pgjd}
Let $G$ be a topological group and $f\in \H(G)$ a Hartman function. Let
$F_1\in \Riemanntop_{\mu_1}(C_1)$ and $F_2\in\Riemanntop_{\mu_2}(C_2)$ be realizations of $f$ on the group compactifications $(\iota_1,C_1)$ resp.\ $(\iota_2,C_2)$. If $F_1$ has a g.j.d., then $F_2$ also has a g.j.d.
\end{pro}
\begin{proof} Let $x\in G$ be a g.j.d.\ for $F_1$.
Suppose $(\iota_1,C_1)\ge(\iota_2,C_2)$, i.e.\ that there is a
continuous surjection $\pi: C_1 \to C_2$ with
$\iota_2 = \pi \circ \iota_1$ and $f=F_1 \circ \iota_1 = F_2 \circ \iota_2$. Thus $F_1$ and $F_2 \circ \pi$ coincide on $\iota_1(G)$.
\begin{diagram}[textflow]
        &               &       C_1   &           &\\
        &\ruTo^{\iota_1}&\dTo>{\pi}   &\rdTo^{F_1}&\\
G       & \rTo^{\iota_2}&        C_2  &\rTo^{F_2}&\C.
\end{diagram}
(Note that the right triangle in the diagram does not necessarily commute on the whole set $C_1$.)
Hence Lemma \ref{Lriemann}.1 implies that $F_1=F_2\circ\pi$ $\mu_1$-a.e. and Lemma \ref{Lriemann}.2 implies that
$F_2\circ \pi$ has a g.j.d.\ at $x\in C_1$ whenever $F_1$ has a g.j.d.\ at $x\in C_1$. Pick disjoint open sets $O_1,O_2 \subseteq C_1$ according to the definition of a g.j.d.\ for $F_2\circ\pi$ at $x\in C_2$, i.e.\
$x\in \overline{O_1}\cap \overline{O_2}$ but
$\overline{F_2\circ\pi(O_1)}\cap \overline{F_2\circ\pi(O_2)}=\emptyset$.
Thus $\pi(O_1)$ and $\pi(O_2)$ are disjoint. Since $\pi$ is an open mapping, $\pi(O_1)$ and $\pi(O_2)$ are open sets and such that $\pi(x)\in \overline{\pi(O_1)}\cap\overline{\pi(O_2)}$. Thus $\pi(x)$ is a g.j.d.\ for $F_2$.

In the general case let $\pi$ be the canonical projection $bG\to C_1$
and define $F^b:=F_1\circ\pi$. It is easy to check that if $F_1$ has a
g.j.d.\ at $x\in C_1$, then $F^b$ has a g.j.d.\ at every point of
$\pi^{-1}[\{x\}]$. Moreover $F^b,F_1$ and $F_2$ induce the same
Hartman function $f$ on $G$. Now apply the first part of this proof to the two functions $F^b$ and $F_2$.
\end{proof}

This result shows that being realized by a function with a g.j.d.\ is an intrinsic property of a Hartman function and does not depend on the particular realization. In virtue
of this result we can consider the set of all Hartman
functions such that one (and hence all) realizations lack a g.j.d.

\begin{defi}\label{defhc} Let $G$ be a topological group. Let
\begin{eqnarray*}
\H_c(G)&:=&\{f\in \H(G): \forall (\iota,C)
\quad f=F\circ\iota \mbox{ with }
F\in R_{\mu_C}(C)\\
&&\mbox{ implies that } F \mbox{ has no g.j.d.} \}\\
&=&\{f\in \H(G): \exists (\iota,C)  \quad f=F\circ\iota
\mbox{ with some } F\in R_{\mu_C}(C)\\
&&\mbox{ without any g.j.d.} \}.
\end{eqnarray*}
\end{defi}
In the next section we will see that $\H_c(G)$ enjoys nice algebraic and topological properties.

\section{Hartman functions with small support}\label{Shartmansmall}
Similar to the situation of g.j.d.s for different realizations of a Hartman function also the property of vanishing outside a meager null set does not depend on the special choice of the realization.

\begin{pro}\label{PH0}
Let $G$ be a topological group and $f\in \H(G)$ a Hartman function. Let
$F_1\in \Riemanntop_{\mu_1}(C_1)$ and $F_2\in\Riemanntop_{\mu_2}(C_2)$ be realizations of $f$ on the group compactifications $(\iota_1,C_1)$ resp.\ $(\iota_2,C_2)$. If $[F_1\neq 0]$ is a meager $\mu_1$-null set, then $[F_2\neq 0]$ is a meager $\mu_2$-null set.
\end{pro}
\begin{proof}
\hfill\par\vspace{-2ex}
\begin{enumerate}
\item
First consider the case where $(\iota_1,C_1)\le (\iota_2,C_2)$ via $\pi: C_2\to C_1$. By assumption $[F_1\neq 0]$ is a meager $\mu_1$-null set.
Use Lemma \ref{Lriemann}.1 to see that $[F_2=F_1\circ\pi]$ is a co-meager set of full $\mu_2$-measure.
Thus $\pi^{-1}[[F_1\neq 0]]\triangle[F_2\neq 0]$
is a meager $\mu_2$-null set. This implies
\[ \mu_2([F_2\neq 0])=\mu_2(\pi^{-1}[[F_1\neq 0]])
=\mu_1([F_1\neq 0])=0.\]

Next we show that $[F_2\neq 0])$ is meager. Indeed, $\pi: C_2\to C_1$ is open, closed, continuous and surjective. Thus one easily verifies that preimages of meager sets are meager, in particular if $[F_1\neq 0]$ is meager
in $C_1$, then $\pi^{-1}[[F_1\neq 0]]$ is meager in $C_2$.
Since $\pi^{-1}[[F_1\neq 0]]$ and $[F_2\neq 0]$ differ at most on
a meager $\mu_2$-null set $[F_2\neq 0]$ is meager, proving the claim.

\item Suppose $(\iota_1,C_1) \ge (\iota_2,C_2)$.  We use again that
$\pi: C_1\to C_2$ is an open and continuous surjection of compact spaces to conclude that $\pi[ [F_1\neq 0] ]$ is meager in $C_2$ whenever $[F_1\neq 0]$ is meager in $C_1$. The rest of the proof is analog to the first case.
\end{enumerate}
In the general case the property of vanishing outside a meager null-set transfers first by 1. from $(\iota_1,C_1)$ to
$(\iota_b,bG)$ and then by 2. from $(\iota_b,bC)$ to $(\iota_2,C_2)$.
\end{proof}

We define the set of those Hartman functions all realizations of which
vanish outside a meager null set.

\begin{defi} \label{defh0} Let $G$ be a topological group. Let
\begin{eqnarray*}
\H_0(G)&:=&\{f\in \H(G): \forall (\iota,C)\quad f=F\circ\iota \mbox{ with }
F\in R_{\mu_C}(C)\\
&&\mbox{ implies that } [F\neq 0] \mbox{ is a meager $\mu_C$-null set }\}\\
&=&\{f\in \H(G): \exists (\iota,C)\quad f=F\circ\iota \mbox{ with some }F\in R_{\mu_C}(C)\\
&&\mbox{ such that } [F\neq 0] \mbox{ is a meager $\mu_C$-null set }\}.
\end{eqnarray*}
\end{defi}

\begin{pro}\label{Ph0hcalgebras}
Let $G$ be a topological group.
$\H_0 (G)$ and $\H_c(G)$ are translation invariant $C^*$-subalgebras of $B(G)$. Furthermore, $\H_c(G)$ contains all constant functions.
\end{pro}

\begin{proof}
By their definition  it is clear that $\mathcal{H}_0(G)$
and $\mathcal{H}_c(G)$ are subalgebras of $B(G)$, invariant under translations and that $\H_c(G)$ contains all constants. It remains to prove that $\H_0(G)$ and $\H_c(G)$ are closed in the topology of uniform convergence.
\begin{enumerate}
\item
Let $\Riemanntop_0(bG):=\{f\in \Riemanntop_{\mu_b}(bG): [f\neq 0] \mbox{ is a meager } \mu_b\mbox{-null set}\}$.
Note that $\Riemanntop_0(bG)$ is a closed subalgebra of $\Riemanntop_{\mu}(bG)$ (due to
the fact that a countable union of meager null sets is again
a meager null set). Since $\iota_b^*: \Riemanntop_{\mu_b}(bG)\to B(G)$ is a continuous homomorphism of $C^*$-algebras and $\iota_b^*(\Riemanntop_0(bG)) = \H_0(G)$
(Definition \ref{defh0}), $\H_0(G)$ is closed (\cite[Theorem I.5.5]{Davi96}).

\item $J(bG)$, the set of all bounded functions
on $bG$ having a g.j.d.\ , is open in the topology of uniform
convergence (Proposition \ref{Pgjdclosed}).
Thus $C(bG)\oplus \Riemanntop_0(bG)$, the set of all bounded functions on $bG$ without a g.j.d.\ (Corollary \ref{Cgjdchar}), is closed (Proposition \ref{Pgjdclosed}).
$\H_c(G)=\iota_b^* (C(bG)\oplus R_0(bG))$ by \cite[Theorem) I.5.5]{Davi96}, thus $\H_c(G)$ is closed.\qedhere
\end{enumerate}
\end{proof}
The last part of this section is devoted to the relations of the algebras $\H_0$, $\H_c$ and $AP$. Note that $AP(G)\cap \H_0(G)=\{0\}$. This is due to the fact that $f\in \H_0(G)$ implies $m(|f|)=0$,
which is impossible for a nonzero almost periodic function.

\begin{lem}\label{LH0tmHc}
Let $G$ be a topological group. Then $\mathcal{H}_0 (G)\subseteq\mathcal{H}_c(G)$.\end{lem}
\begin{proof}
It suffices to show that for every $F\in \Riemanntop_{\mu_C}(C)$ on a group compactification $(\iota, C)$ such that $F\circ\iota \in \mathcal{H}_0 (G)$ there are no two distinct open sets $O_1,O_2 \subseteq  C$ with
$\overline{F(O_1)}\cap\overline{F(O_2)}=\emptyset$.
As $[F\neq 0]$ is a $\mu_C$-null set (Proposition \ref{PH0})
the set $[F=0]$ is dense in $C$, i.e.\ $0\in \overline{F(O_1)}\cap\overline{F(O_2)}$.
\end{proof}

\begin{pro}\label{directsum}
Let $G$ be a topological group. For every $f\in\mathcal{H}_c(G)$ there exists a unique almost periodic function
$f_a\in AP(G)$ and a unique function $f_0\in \H_0(G)$ such that $f:=f_a+f_0$. Furthermore if $f\ge 0$, then $f_a\ge 0$.
\end{pro}
\begin{proof}
Let $F$ be a realization of $f$ on a group
compactification $(\iota,C)$. Using Proposition \ref{Pregularization} we can decompose $F=F^r+(F-F^r)$, the first summand being continuous and
the second one having support on a meager $\mu_C$-null set.

Existence: Let $f_a:=F\circ\iota\in AP(G)$ and $f_0:=(F-F^r)\circ\iota \in H_0(G)$. By construction $f=f_a+f_0$.

Uniqueness: Suppose $f=f_a^{(1)}+f_0^{(1)}=f_a^{(2)}+f_0^{(2)}$ with
$f_a^{(1)},f_a^{(2)}\in AP(G)$ and $f_0^{(1)},f_0^{(2)}\in \H_0(G)$.
This implies $f_a^{(1)}-f_a^{(2)}=f_0^{(1)}-f_0^{(2)}\in AP(G)\cap\H_0(G)= \{0\}$, i.e.\ $f_a^{(1)}=f_a^{(2)}$ and $f_0^{(1)}=f_0^{(2)}$.

Positivity: Let $f\ge 0$. We claim that $F\ge 0$ outside $\mbox{disc}(F)$. To see this, let $x\in C$ be a point of continuity for $F$ and suppose by contradiction that $F(x)<0$. Pick an open neighborhood $V$ of $x$ such that
$F(y)<0$ for any $y\in V$. As $\iota(G)$ is dense in $C$ there exists
an element $\iota(g)\in V$ with $F(\iota(g))=f(g)\ge 0$. Contradiction.
$F$ and $F^r$ coincide on the dense set $C\setminus \mbox{disc}(F)$. By continuity of $F^r$ we have $F^r\ge 0$, implying $f_a\ge 0$.
\end{proof}

An immediate consequence is:
\begin{thr}\label{Tdecom}
Let $G$ be a topological group. Then $\H_c(G) =  AP(G) \oplus \H_0(G)$.
Furthermore the mapping $P: \H_c \to AP(G)$
defined via $f\mapsto f_a$, where
$f=f_a+f_0$ is the decomposition from Proposition \ref{directsum},
is a bounded positive projection with $\|P\|=1$ and
$m(Pf)=m(f)$ for the unique invariant mean $m$ on $\H(G)$.
\end{thr}

Recall from Example \ref{Emap} that a topological group $G$ is called \textbf{minimally almost periodic} (map) if $AP(G)$ consists only of the constant functions. $G$ is called \textbf{maximally
almost periodic} (MAP) if $AP(G)$ separates the points of $G$, cf. Proposition \ref{PHartman}. Every LCA group is maximally almost periodic.

\begin{cor}\label{Cmap}
Let $G$ be a topological group. The following assertions are equivalent:
\begin{enumerate}
\item $\H_c(G)=\H_0(G)$,
\item $G$ is minimally almost periodic.
\end{enumerate}
\end{cor}

\begin{proof} For map $G$ the Bohr compactification $bG=\{0\}$ consists of only one element.\end{proof}

\begin{problem}
 For which topological groups is the inclusion
 $\H_c(G)\supseteq \H(G)\cap \W(G)$ strict? Construct
 $f\in \H_c(G)\setminus( \H(G)\cap \W(G) )$. \end{problem}

\begin{lem}\label{Lnullset}
Let $G$ be a non compact
topological group and let $(\iota_b, bG)$ be the Bohr
compactification of $G$.
\begin{enumerate}
\item If $G$ is MAP then
$\mu_b(\iota_b(K))=0$ for every $\sigma$-compact $K\subseteq G$.
\item If $G$ is an LCA group and $\iota_b(G)$ is $\mu_C$-measurable
then $\mu_b(\iota_b(G))=0$.
\end{enumerate}
\end{lem}
\begin{proof}
\hfill\par\vspace{-2ex}
\begin{enumerate}
\item First suppose that $K$ is compact.
We inductively construct a sequence $(g_i)_{i=1}^{\infty}\subseteq G$ such that
$g_iK\cap g_jK=\emptyset$ for $i\not=j$: Suppose that $(g_iK)_{i=1}^n$
is a family of pairwise disjoint sets; we prove that
there exists $g_{n+1}\in G$ such that
$(g_iK)_{i=1}^{n+1}$ is also a family of pairwise disjoint sets.
Suppose by contradiction that
for every $g\in G$ there is a $j$ such that $g_jK\cap gK\not=\emptyset$. Then $g\in g_j KK^{-1}$. So $G=\bigcup_{j=1}^ng_j KK^{-1}$ would be compact. Contradiction.

Since $G$ is MAP, $\iota_b$ is one-one.
The sets $(\iota_b(g_iK))_{i=1}^{\infty}$
form an infinite sequence of pairwise disjoint translates of the
compact (and thus measurable) set $\iota_b(K)\subseteq bG$.
If $\mu_b(\iota_b(K))>0$ then
\[ 1=\mu_b(bG)\ge \sum_{i=1}^{\infty} \mu_b(\iota_b(g_i K)) =
\sum_{i=1}^{\infty} \mu_b(\iota_b(K))=\infty. \] Contradiction. Consequently $\mu_b(\iota_b(K))=0$. If $K$ is $\sigma$-compact the assertion follows from the
$\sigma$-additivity of $\mu_b$.

\item Follows from the fact that $\iota_b(G)$ has zero outer $\mu_b$-measure, see \cite{Glic62,Varo65}. \qedhere
\end{enumerate}
\end{proof}
If we replace in Lemma \ref{Lnullset} the Bohr compactification by an arbitrary compactification $(\iota,C)$ the measurability condition on the
set $\iota(G)$ becomes crucial.

\begin{exa}
Consider the compact group $\T = \R/\Z$ and any
fixed irrational $\alpha \in \T$. By Zorn's Lemma there is a maximal subgroup $G$ of $\T$ with $\alpha \notin G$.
$G$ equipped with the discrete topology is an LCA group. Let $\iota: G \to \T$ be the inclusion mapping and $C=\T$, then $(\iota,C)$ is a group compactification of $G$, distinct from the Bohr compactification.
Let $\mu_C$ be the Haar measure on $C$.
Assume, by contradiction, that $G$ is a $\mu_C$-measurable null set in $\T$. Consider the measure preserving mappings $\phi_k: \T \to \T$,
$x \mapsto kx$, $k \in \Z$. Then all sets
$\phi_k^{-1}[\alpha + G]$ are measurable $\mu_C$-null sets. Pick any $x \in \T$. If $x \notin G$ then, by the maximality property of $G$, $\alpha = kx+g$ for some $k \in \Z \setminus \{0\}$, $g \in G$.
This implies $\phi_k(x) \in \alpha + G$, i.e.\ $x \in \phi_k^{-1}[\alpha + G]$. We conclude that
\[\T = G \cup \bigcup_{k \in \Z \setminus \{0\}} \phi_k^{-1}[\alpha +G]\]
is the countable union of $\mu_C$-null sets, hence $1 = \mu_C(C) = 0$, contradiction.
\end{exa}

Let us by $\F_0(G)$ denote the set of all bounded (not necessarily continuous or even measurable) complex valued functions $f: G \to \C$ vanishing at infinity, i.e.\ $f \in \F_0(G)$ if for every $\eps > 0$
there is a compact set $K \subseteq G$ with $|f(x)| < \eps$
for all $x \in G \setminus K$. As usual $C_0(G)$ denotes the set of all \emph{continuous} $f \in \F_0(G)$.\label{f0andc0}

\begin{thr}\label{TC0}
Let $G$ be a MAP group.
Then $C_0(G)\subseteq \mathcal{H}(G)$. If $G$ is not compact then even $\F_0(G)\subseteq \mathcal{H}_0(G)$.
\end{thr}
\begin{proof}
In the first step we show $\F_0(G)\subseteq \mathcal{H}(G)$.
If $G$ is compact there is nothing to prove. Suppose $G$ is not compact. Let $f\in \F_0(G)$ and define $F: bG \to \C$ by
\[
F(x):=\left\{ \begin{array}{cl} f(g)& \mbox{ if } x=\iota_b(g),\;g\in G \\
0& \mbox{ else.} \end{array} \right.
\]
Then $f=F\circ\iota_b$. It suffices to consider $f$ such that $0\le f\le 1$.
For every $\eps>0$ there exists a
compact set $K_{\eps}\subseteq G$ such that $f(x) < \eps$ for $x\in
G\setminus K_{\eps}$. By Lemma \ref{Lnullset}, we have $\mu_b(A) = 0$,
where $A = \iota_b(K_{\eps})$. Regularity of the Haar measure implies that we can find an open set $O\supset A$
such that $\mu_b(O)<\eps$. Let $h$ be an Urysohn function for $A$ and $bG\setminus O$, i.e.\ $h: bG \to [0,1]$ is continuous with $h = 1$ on $A$ and $h = 0$ on $bG \setminus O$. Consider the continuous function
$g_{\eps}:=h+\eps \Eins_{bG}$. Since $0\le F\le g_{\eps}$ and
\[\int_{bG}g_{\eps}d\mu_b\le \mu_b( [h>0])+\eps\le 2\eps,\]
Proposition \ref{Priemannchar} implies $F\in \Riemanntop_{\mu}(bG)$.
Hence $f$ is a Hartman function.
It remains to show that $[F\neq 0]$ is a meager $\mu_b$-null set.
For each $n\in\N$ the set $[f \ge 1/n]$ is compact. Hence
$\iota_b( [f \ge 1/n]) = [F \ge 1/n]$ is
a compact $\mu_b$-null set and therefore
nowhere dense. Thus $[F\neq 0]=\bigcup_{n=1}^{\infty}[|F|\ge 1/n]$
is a meager $\mu_b$-null set and $f\in\H_0(G)$.
\end{proof}

\begin{cor}\label{Cnotmeasurable}
Hartman functions $f \in \H(G)$ need not be measurable
with respect to the completion of the Haar measure on $G$.
\end{cor}
\begin{proof}
As a counterexample take $G=\R$ with the Lebesgue measure and
any set $A \subset [0,1]$ which is not Lebesgue measurable.
Then $f = \Eins_A$ is a Hartman function by Theorem \ref{TC0}
but not Lebesgue measurable.
\end{proof}

As a further consequence of Theorem \ref{TC0} we get the following supplement to Corollary \ref{Chnichtw}.

\begin{cor}\label{CH0nichtW}
Let $G$ be a non-discrete MAP group. Then $\H_0(G) \setminus \W(G) \neq \emptyset$. In particular the inclusion $C_0(G)\subset \H_0(G)$ is strict.
\end{cor}
\begin{proof}
Let $f := \Eins_{\{0\}}$. Then $f \in \H_0(G)$ (trivially for
compact $G$, otherwise by Theorem \ref{TC0}). Since $f$ is not continuous $f\notin\W(G)$. (Recall that every weakly almost periodic $f$ has a representation $f = F \circ \iota_w$ with $F: wG \to \C$ continuous on the weakly almost periodic compactification
$(\iota_w,wG)$ and thus is continuous.)
\end{proof}

The following example shows that also for the integers the space $C_0(\Z)$ of functions vanishing at infinity is a proper subspace of $\H_0(\Z)$.

\begin{exa}\label{Eh0abernichtc0}
Let $T=\{t_n: n\in \N\}$
be a \textbf{lacunary} set of positive integers, i.e.\
$t_1<t_2<t_3<\ldots$ with $\limsup_{n\to\infty} \tfrac{t_n}{t_{n+1}}=\eps < 1$.
Then {\em $\Eins_T$} $\in \H_0(\Z)\setminus C_0(\Z)$.
\begin{proof}
By \cite[Theorem 9]{ScSW00} for
each $n\in\N$ there exists an $n$-dimensional
compactification $(\iota_n, C_n)$ and a compact $\mu_{C_n}$-continuity set $K_n \subseteq C_n$ with $\mu_n(K_n)\le 4n\eps^n$ such
that $\iota_n^{-1}[K_n]\supseteq T$. Furthermore we can arrange
$(\iota_n,C_n)\le (\iota_{n+1},C_{n+1})$ and
$\pi_{n+1,n}^{-1}[K_n]\supseteq K_{n+1}$,
where $\pi_{n+1,n}: C_{n+1}\to C_n$ is the canonical projection,
i.e.\ $\iota_n = \pi_{n+1,n} \circ \iota_{n+1}$.
Let $(\iota,C):=\bigvee_{n=1}^{\infty} (\iota_n,C_n)$
and let $\pi_n: C\to C_n$ be the canonical projection onto $C_n$. Thus
$K:=\bigcap_{n=1}^{\infty} \pi_n^{-1}[K_n]$ is a compact $\mu_C$-null
set (hence a $\mu_C$-continuity set) with
$\iota^{-1}[K]\supseteq T$. This shows $\Eins_T\in \H_0(\Z)$.
Since $T\subseteq \Z$ is infinite, we have $\Eins_T \notin C_0(\Z)$. \end{proof} \end{exa} 
\section{Hartman functions on $\Z$}\label{ShartmanZ}
For locally compact groups $G$ it is very easy to see that $C_0(G) \subseteq \W(G)$. A much harder problem is finding functions
in $\W(G)\setminus(AP(G)\oplus C_0(G))$, see for instance \cite{Rupp91}.

Topological groups with the property
$\mathcal{W}(G)=\mathcal{A}(G)\oplus C_0(G)$ are
called \textbf{minimally w.a.p.} A famous example, due to M.\ Megrelishvili, is $H_+[0,1]$, the group of all orientation-preserving self-homeomorphisms of the closed unit interval $[0,1]$ endowed with the compact-open topology, see \cite{Megr01}. For minimally w.a.p.\ groups our Theorem \ref{TC0} implies $\mathcal{W}(G)\subseteq\mathcal{H}(G)$. However, it is known that non-compact LCA groups are never minimally w.a.p., see \cite{Chou80}.

\begin{problem} Find a non-trivial topological group $G$ (necessarily not minimally w.a.p.) such that $\mathcal{W}(G)\subseteq\mathcal{H}(W)$.
\end{problem}

Throughout the rest of this section all results are stated for the case $G=\Z$. A quick way to obtain $f\in (\mathcal{W} \cap \H_0) \setminus (AP \oplus C_0)$ is implicated by the following result.

\begin{pro}\label{Placunary}
Let $(t_n)_{n=1}^{\infty}\subseteq \Z$ be a lacunary sequence of positive integers, i.e.
$$\limsup_{n\to\infty} \frac{t_n}{t_{n+1}}=\eps<1.$$
Let $T:=\{t_n : n\in \N\}\subseteq \Z$, then
$\Eins_T\in (\W\cap \H_0)\setminus (AP\oplus C_0)$.
\end{pro}
\begin{proof}
According to our example \ref{Eh0abernichtc0}
$f=\Eins_T$ is a member of $\H_0$. Since $T$ is a lacunary set \cite[Theorem 4.2]{BeJM89}
implies $f\in \W$. Furthermore, $\liminf_{k\to \pm \infty}=1$ implies $f\notin C_0$ and $f\notin AP$. Suppose by contradiction that $f=f_a+f_0$ where $f_a\in AP$ and $f_0\in C_0$. Then
\[0=\mbox{dens}(T)=m(f)=m(f_a)+m(f_0)\]
implies $m(f_a)=0$. As $f_a(k)\ge 0$ for all but finitely many $k$ this implies $f_a=0$. Contradiction.
\end{proof}

The main objective of this Section is now to illustrate a further method to construct functions $f\in(\mathcal{W} \cap \H_0) \setminus (AP \oplus C_0)$.

\subsection{Fourier-Stieltjes transformation}\label{SSfourier}
Let us recall some facts about the Fourier transformation of
measures on LCA groups. Let $G$ be an LCA group. By $\M(G)$ we denote the set of all finite complex Borel measures on $G$. Recall that $\M(G)$ can be regarded as the dual  $C_0(G)^*$ of the Banach space $C_0(G)$. The canonical pairing $C_0(G)\times C_0(G)^* \to \C$ is given by

\[\langle f,\mu\rangle := \int_G f(x) d \mu (x).\]
Also recall that we convolute two measures $\mu,\nu \in \M(G)$ according to the formula
\[\langle f, \mu*\nu\rangle = \int_{G\times G} f(x+y) d(\mu\otimes\nu)(x,y).\]

The Fourier-Stieltjes transform $\mu\mapsto \hat{\mu}$ assigns to a measure $\mu\in \M(G)$ the uniformly continuous function
\[ \hat{\mu}(\chi):=\int_G \chi(x)d\mu (x)\]
defined on the dual group $\hat G$. The map $\mu\mapsto \hat{\mu}$ is a continuous homomorphism of the convolution algebra $(\M(G),*)$ into the function algebra $(UCB(\hat{G}), \cdot)$ of uniformly continuous functions on $\hat G$. The set $\{\hat{\mu}:\mu\in \M(G)\}$, of all Fourier-Stieltjes transforms, is called the Fourier-Stieltjes algebra and denoted by $\mathcal{B}(\hat{G})$.
It is well known (see \cite{BeJM89,Rudi90}) that for non-compact LCA
groups $G$ the inclusions
$$AP(G)\subset \overline{\mathcal{B}(G)}\subset \mathcal{W}(G)$$
hold and are strict.

\begin{pro}\label{Pffacts}
Let $G$ be a locally compact group. The following assertions hold:
\begin{enumerate}
\item If $\mu$ is discrete, then $\hat{\mu}$ is almost periodic,
\item If $\mu$ is absolutely continuous with respect to the Haar measure on $G$, then $\hat{\mu}\in C_0(\hat{G}) \subseteq
\mathcal{W}_0(\hat{G})$
(for $G=\T$ this is the Riemann-Lebesgue Lemma),
\item $m_{\hat{G}}(\hat{\mu})=\mu(\{0_G\})$ for the unique invariant mean $m_{\hat{G}}$ on $\hat{G}$. In particular $\hat{\mu}$ has zero mean-value whenever $\mu$ is atomless.
\end{enumerate}
\end{pro}
\begin{proof} See \cite[Section 1.3]{Rudi90}. \end{proof}

Recall that an LCA group, by Pontryagin's duality theorem, is algebraically and topologically isomorphic to its bi-dual.

\begin{lem}\label{Lftlemma}
Let $G$ be a discrete LCA group and
$(\nu_n)_{n=1}^{\infty}\subseteq \mathcal{M}(\hat{G})$ a bounded sequence of discrete measures. Then the following assertions are equivalent:
\begin{enumerate}
\item The sequence $(\hat{\nu}_n)_{n=1}^{\infty}\subseteq AP(G)$
converges pointwise to a bounded function $f: G\to C$.
\item The sequence $(\nu_n)_{n=1}^{\infty}$ of discrete measures
converges weak-* to a measure $\mu\in \mathcal{M}(\hat G)$.
\end{enumerate}
In this case $f=\hat \mu$, the Fourier-Stieltjes transform of the measure $\mu$.
\end{lem}
\begin{proof}
(1) $\Rightarrow$ (2): Let $f_n:=\hat{\nu}_n$.
By weak-*-compactness of the unit
ball in $C_0(\hat{G})^*=\mathcal{M}(\hat{G})$ we can find a
weak-*-limit-point $\mu$ of the set
$\{\nu_n: n\ge 0\}$.
Due to compactness of $\hat{G}$, for every $x\in G$ the map
$\mu\mapsto \int_{\hat{G}} \chi(x) d\mu(\chi)$ is a weak-*-continuous
functional defined on $\mathcal{M}(\hat{G})$. Thus
for every $x\in G$ and $\eps>0$ there exist
infinitely many $n_k$, $k \in \N$, (depending of course on $x$) such that
\[
|\hat{\mu}(x)-f_{n_k}(x)|=\left|\int_{\hat{G}} \chi(x)
d\mu(\chi)-\int_{\hat{G}} \chi(x) d\nu_{n_k}(\chi)\right| < \eps.
\]

Using that $f_n(x)\to f(x)$ pointwise, we obtain
\[ |\hat{\mu}(x)-f(x)|\le |\hat{\mu}(x)-f_{n_k}(x)|+ |f_{n_k}(x)-f(x)|
                    < 2 \eps.\]
Thus $\lim_{n\to\infty} f_n(x)=\hat{\mu}(x)$. Let $\tilde{\mu}$ be another weak-*-limit point of the set $\{\nu_n: n\in \mathbb{N}\}$. On a compact space weak-*-convergence of measures implies pointwise convergence of their Fourier-Stieltjes transforms. Thus $\hat{\tilde{\mu}}$ and $\hat{\mu}$ coincide. Hence $\mu = \mbox{weak-*-}\lim_{n\to\infty} \nu_n=\tilde\mu$.

(2) $\Rightarrow$ (1): By compactness
of $\hat{G}$, for every $x\in G$ the mapping
$\mu \mapsto \hat{\mu}(x)=\int_{\hat{G}} \chi(x)d\mu(\chi)$ is a weak-*-continuous functional. Thus $f_n:=\hat{\nu}_n$ converges pointwise.
\end{proof}

\subsection{Example}\label{SSexample}

In the following we will investigate the function
\[
f(k) = \prod_{j=1}^{\infty} \cos^2 \left(2 \pi \frac k{3^j}\right)
\]
defined on the group $G = \Z$ of integers and the proof of Theorem
\ref{Tcantor} below. For its formulation we use the singular measure
$\mu_3$ concentrated on the ternary Cantor (middle-third) set in the natural way. To be more precise: Let $\lambda$ be the Lebesgue measure on $[0,1)$. Consider the $\lambda$-almost everywhere uniquely defined mapping
$\phi: [0,1) \to [0,1)$ with
\[
\phi: \sum_{i=1}^{\infty} \frac{a_i}{2^i} \mapsto
\sum_{i=1}^{\infty} \frac{2a_i}{3^i},
\]
$a_i \in \{0,1\}$, and the canonical inclusion $\iota: [0,1) \to \T = \R/\Z$,
$ x \mapsto x+\Z$. Then $\mu_3 = (\iota \circ \phi)\circ \lambda$
(notation as in Proposition \ref{Pcomprepr}).

\begin{thr}
\label{Tcantor}
Let $f: \Z \to [0,1]$ be given by
\[
f(k) = \prod_{j=1}^{\infty} \cos^2 \left(2 \pi \frac k{3^j}\right)
\]
Then the following statements hold:
\begin{enumerate}
\item $f \in (\H_0 \cap \W) \setminus (AP \oplus C_0)$.
\item $m_{\Z}(f)=0$ for $m_{\Z}$ the unique invariant mean on $\H(\Z)$.
\item $f$ can be realized by a Riemann integrable function on the
3-adic compactification $\overline{\Z}^{(3)}$.
\item
$f$ is the Fourier-Stieltjes transform of the singular measure $\mu_3$ corresponding to the ternary Cantor set canonically embedded into $\T$.
\end{enumerate}
\end{thr}
\begin{proof}
Everything will follow from Lemma \ref{finw0},
\ref{fnichtac0} and \ref{f3adisch}.
\end{proof}

We have to fix some notation and then prove the auxiliary statements.
We will construct a function on $\mathbb{Z}$ using
discrete measures on $\hat{\mathbb{Z}}=\mathbb{T}$.
Note that $\mathbb{T}$ is algebraically and topologically isomorphic to the interval $[0,1)$ equipped with addition modulo $1$.

For $\alpha \in [0,1)$ let us denote by $\delta_{\chi_{\alpha}}\in \mathcal{M}(\hat{\Z})$ the probability measure
which is concentrated on the character $\chi_{\alpha}: k \mapsto \exp{(2\pi i k \alpha)}$. We define recursively discrete probability measures $\nu_n\in \mathcal{M}(\Z)$ by
$\nu_0:=\delta_{\chi_{1/2}}$ and
\[\nu_{n}:=\nu_{n-1} * \frac{1}{2}
\left( \delta_{\chi_{-1/3^n}}+\delta_{\chi_{1/3^n}}\right).\]
Note that $\nu_n\to \mu_3$ in the weak-*-topology of $\mathcal{M}(\Z)$.
Using the fact that
$(\nu_n*\nu_{n-1})\,\hat{}=\hat{\nu}_n\hat{\nu}_{n-1}$ and
$\hat\delta_{\chi_{\alpha}}(k)=\chi_{\alpha}(k)=\exp{(2\pi i k \alpha)}$
one easily computes $\hat{\nu}_n(k) =\prod_{j=1}^{n} \cos\left(2\pi \tfrac{k}{3^j}\right)$.

\begin{lem} \label{finw0}
Let
\[ \tilde f_n := \hat{\nu}_n(k) =
\prod_{j=1}^{n} \cos\left(2\pi \frac{k}{3^j}\right).
\] Then $\tilde f_n$ converges pointwise to $\hat{\mu_3}$, the Fourier-Stieltjes transform of the singular measure $\mu_3$ concentrated on the ternary Cantor set. In particular $\lim_{n\to\infty} \tilde f_n$ is weakly almost periodic and has zero-mean value.
\end{lem}
\begin{proof}
Each $\tilde{f}_n$ is a product of finitely many periodic
factors with rational periods, so  $\tilde{f}_n$ is periodic.
We show that the functions $\tilde{f}_n$ converge pointwise.
Observe that $\lim_{j\to\infty}\cos\left(2\pi \tfrac{k}{3^j}\right)= 1$ for fixed $k\in \N$. All terms of this sequence are non-negative provided $j\ge \log_3 (2k)=:j(k)$. Thus $\left(\tfrac{\tilde{f}_{j(k)+n}(k)}{\tilde{f}_{j(k)}(k)}\right)_{n=1}^{\infty}$ is a monotonically decreasing sequence of non-negative real numbers, hence
its limit exists. By Lemma \ref{Lftlemma} we know that
\[ \tilde{f}(k) := \lim_{n\to\infty}\tilde f_n (k)= \prod_{j=1}^{\infty} \cos\left(2\pi \frac{k}{3^j}\right) \]
is the Fourier-Stieltjes transform of the measure $\mu=\mu_3\in\mathcal{M}(\Z)$. Since $\mu_3$
has no atoms, Proposition \ref{Pffacts} implies that
$\tilde f=\widehat{\mu_3}\in\W(\Z)$ and that $m_{\Z}(\tilde f)=0$ for the unique invariant mean $m_{\Z}$ on $\W(\Z)$.

The same considerations apply
to the discrete measures $\nu_n*\nu_n$, the non-negative periodic functions
\[f_n(k):=\tilde{f}_n^2(k)=\prod_{j=1}^{n}
\cos^2\left(2\pi \frac{k}{3^j}\right)\]
and the limit $f=\tilde{f}^2=(\mu_3*\mu_3)\,\hat{}\,$, which is weakly almost periodic with zero-mean value. \end{proof}

\begin{lem} \label{mittelwert}
The periodic functions $f_n: \Z\to [0,1]$ defined via
\[f_n(k):=\tilde{f}_n^2(k)=\prod_{j=1}^{n}
\cos^2\left(2\pi \frac{k}{3^j}\right)\]
have the mean value $m_{\Z}(f_n)=\tfrac{1}{2^n}$, where
$m_{\Z}$ is the unique invariant mean on $AP(\Z)$.
\end{lem}
\begin{proof}
By Proposition \ref{Pffacts} it suffices to compute $(\nu_n*\nu_n)\left(\{0\}\right)$. We leave the elementary calculation to the reader.
\end{proof}

\begin{lem} \label{fnichtac0}
Let $f=\lim_{n\to\infty} f_n$ be as above. Then $f \notin (AP \oplus C_0)$.
\end{lem}
\begin{proof}
$\tilde{f}$ satisfies the functional equation
$\tilde{f}(3k)=\tilde{f}(k)$, $k\in \Z$. This implies \[\tilde{f}(3^k)=\tilde{f}(0)=1.\] Thus both
$\tilde{f},f \notin C_0$.
As $f\ge 0$ but $m_{\Z}(f)= 0$  we have $f\notin AP$ by Corollary \ref{Capnull}.

Suppose there exists a representation $f=f_a+f_0\ge 0$
with non-trivial $f_a\in AP$ and $f_0\in C_0$. Furthermore let
\[ f_a=\underbrace{\max\{f_a,0\}}_{:=f_a^+\ge 0} +
\underbrace{\min\{f_a,0\}}_{:=f_a^-\le 0}.\]
Note that $f_a^+,f_a^-\in AP$ as $AP$ is a lattice.
$m_{\Z}(f_a)=m_{\Z}(f)=0$ implies $m_{\Z}(f_a^+)=-m_{\Z}(f_a^-)$.
As $f_a^-$ is a non-positive almost periodic function $m_{\Z}(f_a^-)<0$. Thus there exists $\eps>0$ such that for all $N \in \N$
\[
\inf_{|k|\ge N} f_a (k) =\inf_{|k|\ge N} f_a^- (k)\le- \eps < 0.
\]
Note that $f_a^-(k_0)\neq 0$ implies $f_a^+(k_0)=0$.
Let $N_0$ be such that $|f_0(k)|<\tfrac{\eps}{2}$ for $|k|\ge N_0$. Thus there
exists $k_0\ge N_0$ such that
\[ f(k_0)=f_a(k_0)+f_0(k_0)=f_a^-(k_0)+f_0(k_0)\le- \eps + \tfrac{\eps}{2} = -\tfrac{\eps}{2}<0. \]
This contradicts $f\ge 0$.
\end{proof}

Consider the compact group of 3-adic integers
$\overline{\Z}^{(3)}$ realized as projective limit of the projective
system of cyclic groups
$C_n:=\mathbb{Z}/3^n\mathbb{Z}$ and mappings $\pi_n$
(reducing $k$ mod $3^{n+1}$ to $k$ mod $3^n$):
\begin{diagram}[small]
     &     &              &    &      &    &
        &    &\overline{\Z}^{(3)}      \\
     &     &              &    &      &    & &
        \ldTo(6,2)^{\kappa_{n-1}}\ldTo_{\kappa_n}    &    \dDots  \\
\{0\}&\lTo_{\pi_0}& \Z/3\Z&\lTo_{\pi_1}&\ldots&\lTo_{\pi_{n-1}}&
   \Z/3^n\Z&\lTo_{\pi_n}&{}
\end{diagram}
The projective limit $\overline{\Z}^{(3)}:=\lim_{\leftarrow} C_n$
of this system can be identified with
a certain closed subgroup of the compact group $\prod_{n=1}^{\infty} C_n$.
Regarding $C_n$ as the set $\{0,1/3^n,\ldots, 1-1/3^n\}$
with addition modulo 1, one
easily checks that for each integer $k\in \Z$ the sequence
$\iota(k):=(k/3^n)_{n=1}^{\infty}$ defines an element of the
projective limit  $\overline{\Z}^{(3)}$.  The mapping
$\iota: \Z \to \overline{\Z}^{(3)}$ is a (continuous) homomorphism. Hence
$(\iota, \overline{\Z}^{(3)})$ is a group compactification of $\Z$, the so
called 3-adic compactification. Note that  each $(C_n,\iota_n)$ is a group compactification
of $\Z$, where $\iota_n$ is reduction modulo $3^n$.
Furthermore $(C_n,\iota_n)\le (C_{n+1},\iota_{n+1})$
via $\pi_n$ and $(C_n,\iota_n) \le (\iota, \overline{\Z}^{(3)})$
via $\kappa_n$ for each $n\in\N$.  By construction every $3^n$-periodic function $f: \Z\to \C$ can be realized by a
continuous function $F: C_n\to \C$.

\begin{lem}\label{f3adisch} Let $f=\lim_{n\to\infty} f_n$ be as above.
Then $f\in\mathcal{H}_0$  and $f$ can be realized in the 3-adic integers.
\end{lem}

\begin{proof}Since every $3^n$-periodic function can be realized by a
continuous function on $C_n$, we can in particular realize
$f_n:=\prod_{j=1}^n \cos^2\left(2\pi \tfrac{k}{3^j}\right)$.
Consequently there exists a unique continuous function $F_n$ on
the 3-adic integers $\overline{\Z}^{(3)}$
such that $f_n=F_n\circ\iota$.

Since for $x\in \iota(\Z)$ the sequence of $(F_n(x))_{n=1}^{\infty}$
is decreasing (note that $0\le \cos^2\left(2\pi \frac{k}{3^j}\right)\le 1$),
$(F_n(x))_{n=1}^{\infty}$
is decreasing for every $x\in \overline{\Z}^{(3)}$ by continuity of $F_n$. In particular the limit $F(x):=\lim_{n\to\infty}F_n(x)$
exists and $F\circ\iota=f$.
We show that $F$ is Riemann integrable on $\overline{\Z}^{(3)}$:
For each $n\in\N$ we have $0 \le F \le F_n$.
Lemma \ref{mittelwert} and
uniqueness of the invariant mean $m_{\Z}$ on $AP$ yield
\[\int_{\overline{\Z}^{(3)}} F_n d\lambda = m_{\Z}(f_n)=\frac{1}{2^n}\]
for the normalized Haar measure $\lambda$ on $\overline{\Z}^{(3)}$.
Thus Proposition \ref{Priemannchar} implies that $F$ is Riemann
integrable on $\overline{\Z}^{(3)}$.

Finally, suppose $F$ has a g.j.d. Then Theorem \ref{Tjump} implies $f\notin \W$ contradicting Lemma \ref{finw0}. Thus
$f\in \H_c$. By Proposition \ref{Pregularization} we can find unique functions $f_a\in AP$ and
$f_0\in\H_0$ such that $f=f_a+f_0$. As $f\ge 0$ we have $f_a\ge 0$. $m_{\Z}(f_a)=m_{\Z}(f)=0$ implies $f_a=0$. So, indeed $f=f_0\in\H_0$.
\end{proof}

\begin{problem} Construct functions $f_1\in \overline{\mathcal{B}}\setminus \H$ and  $f_2\in \H\setminus \overline{\mathcal{B}}$.
\end{problem}

\begin{problem} How are $\overline{\mathcal{B}}$ and $\mathcal{W}\cap\mathcal{H}$ related?
Is there a reasonable condition on functions in $\overline{\mathcal{B}}$
that implies Hartman measurability?
\end{problem} 
\chapter{Summary}\label{CHsum}
The following diagram summarizes some
of our results concerning the space $\H = \H(G)$ of Hartman functions on a topological group $G$. Recall the following
function spaces:

\begin{tabular}{lcl}
        &        &                                                      \\
$AC$    &$\ldots$&  almost convergent functions, see Definition \ref{Dalmostconv}\\
$AP$    &$\ldots$&  almost periodic functions,  see Definition \ref{Dweakap}\\
$\W$    &$\ldots$&  weakly almost periodic functions, see Definition \ref{Dweakap}\\
$\overline{\mathcal{B}}$&$\ldots$& Fourier-Stieltjes algebra, see Section \ref{SSfourier}\\
$\H$    &$\ldots$&  Hartman functions, see Definition \ref{Dhartman} \\
$\H^w$  &$\ldots$&  weak Hartman functions, see Definition \ref{Dweakhartman} \\
$\H_c$  &$\ldots$&  Hartman functions realized without g.j.d.,\\
        &        &  see Definition \ref{defhc}\\
$\H_0$  &$\ldots$&  Hartman functions realized by functions supported \\
        &        &  on a meager null set, see Definition \ref{defh0}\\
$C_0$   &$\ldots$&  continuous functions vanishing at infinity, see pp. \pageref{f0andc0}
\end{tabular}

Inclusions indicated by $\mid$ are proper (at least for certain groups
$G$, e.g.\ for $G=\Z$). For spaces connected by $:$ we did not prove strict inclusions.

\begin{figure}[h]
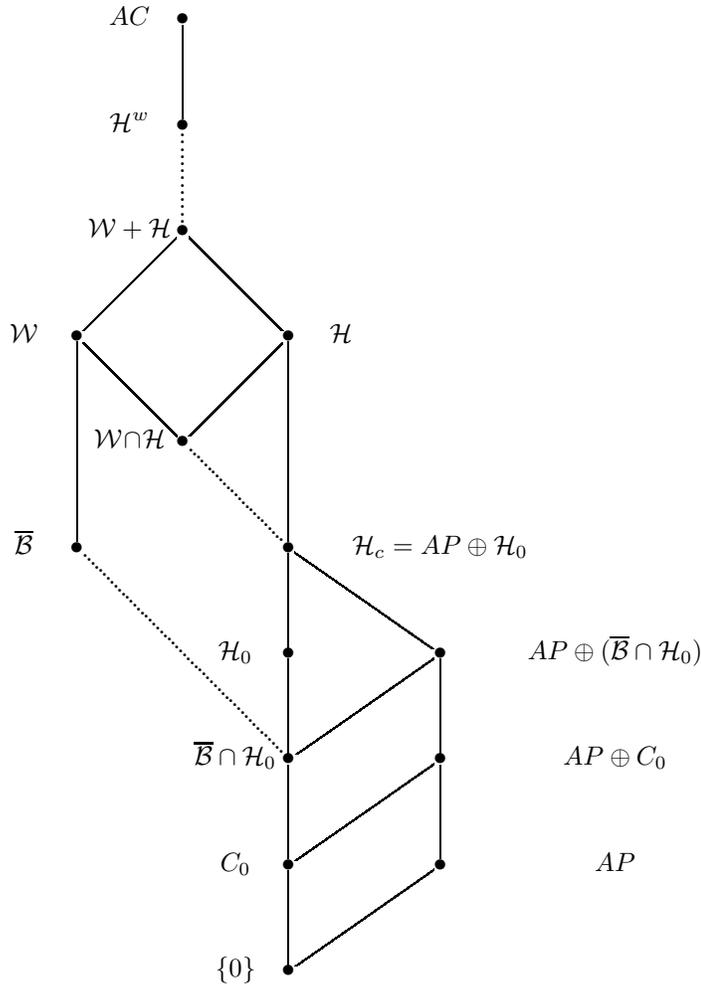
\label{bild}

\begin{center}
\begin{diagram}[small,abut]
&&AC&\bullet&&&                             &&&                 \\
&&&\dLine&&&                              &&&                 \\
&&\H^w&\bullet&&&                         &&&                 \\
&&&\dDots&&&                              &&&                 \\
&&\mathcal{W}+\mathcal{H}&\bullet&&&      &&&                 \\
&&\ldLine&&\rdLine&&                      &&&                 \\
\mathcal{W}&\bullet&&&&\bullet&\mathcal{H}&&&                 \\
&\dLine&\rdLine&& \ldLine&&           &&&                     \\
&&\W\!\cap\!\H&\bullet         \\
&&&&\rdDots&\dLine                                     \\
\overline{\mathcal{B}}&\bullet&&&&\bullet&&\H_c=AP\oplus \H_0&&&              \\
&&\rdDots&&&\dLine&\rdLine  \\
&&&&\H_0&\bullet&&\bullet&AP\oplus (\overline{\mathcal{B}}\cap \H_0)     \\
&&&&\rdDots&\dLine&\ldLine&\dLine         \\
&&&&\overline{\mathcal{B}}\cap \H_0&\bullet&&\bullet &AP\oplus C_0                   \\
&&&&&\dLine&\ldLine&\dLine         \\
&&&&C_0&\bullet&&\bullet&AP     \\
&&&&&\dLine&\ldLine        \\
&&&&\{0\}&\bullet     \\
\end{diagram}
\end{center}
\caption{Spaces of Hartman measurable functions}
\end{figure}

\pagebreak
\bibliographystyle{ams-pln}
\bibliography{Hartman}
\addcontentsline{toc}{chapter}{Bibliography}

\end{document}